\RequirePackage{xcolor}
\documentclass[journal,twoside,web]{ieeecolor}
\usepackage{generic}
\usepackage{cite}
\usepackage{amsmath,amssymb,amsfonts}
\usepackage{graphicx}
\graphicspath{{figures/}}
\usepackage{hyperref}
\usepackage{rotating}

\newtheorem{theorem}{Theorem}

\newtheorem{proposition}[theorem]{Proposition}

\let\labelindent\relax
\usepackage{enumitem}
\setlist[itemize]{label=*}

\usepackage{etoolbox}
\makeatletter
\long\def\RevisionColorHbox#1{\hbox{\color@begingroup#1\color@endgroup}}
\patchcmd{\@makecaption}{\hbox}{\RevisionColorHbox}{}{}
\patchcmd{\@makecaption}{\hbox}{\RevisionColorHbox}{}{}
\makeatother
\AtBeginDocument{\hypersetup{hidelinks}}

\begin{document}
\title{Fischer--Burmeister-based Multiplier Dynamics in Inequality-Constrained Optimization}
\author{Kang Liu, Jianchen Hu, \IEEEmembership{Member, IEEE}
\thanks{This work is supported by the National Key Research and Development Program of China under Grant 2025YFA1017000 (\emph{corresponding author: Jianchen Hu}).}
\thanks{Kang Liu is with the School of Future Technology, Xi'an Jiaotong University, Xi'an, China, 710049 (e-mail: \texttt{kanyo@foxmail.com}).}
\thanks{Jianchen Hu is with the School of Automation Science and Engineering, Xi'an Jiaotong University, Xi'an, China, 710049 (e-mail: \texttt{horace89@xjtu.edu.cn}).}
}

\maketitle

\begin{abstract}
We study proportional–integral multiplier feedback for inequality-constrained optimization using the Fischer–Burmeister complementarity residual. We construct explicit dynamics and characterize gain-dependent stability and convergence rates. For a single active constraint, we derive exact local stability and decay-rate conditions, including a tradeoff between continuous-time convergence and discretization costs. For strongly convex objectives with affine constraints, we establish global convergence and characterize geometric conditions for uniform exponential stability. Numerical experiments compare feedback structures and computational costs.
\end{abstract}


\section{Introduction}
\label{sec:introduction}

Continuous-time primal--dual dynamics solve constrained optimization problems through the interaction between primal variables and Lagrange multipliers, connecting optimization convergence with closed-loop stability~\cite{FeijerPaganini2010}. The controlled multipliers optimization (CMO) framework further treats the multipliers as control inputs and constructs constraint feedback through proportional--integral (PI) control and feedback linearization~\cite{CeroneFossonPirreraRegruto2025TAC, CeroneFossonPirreraRegruto2024CDC}. Within this framework, the choice of the constraint error signal and its feedback structure jointly determine the dynamics of the optimization system.

We consider the inequality-constrained optimization problem $\min_{x\in\mathbb R^n}f(x)$ subject to $g(x)\leq0$, where $f$ and the constraint functions are twice continuously differentiable. For equality constraints, the constraint values can directly serve as error signals to be driven to zero. For inequality constraints, strictly satisfied constraints need not approach their boundaries, while the multipliers must also satisfy nonnegativity and complementarity. The Fischer--Burmeister (FB) function represents these conditions through a residual equation~\cite{Fischer1992}. Specifically, define
\begin{equation}
    \Phi_j(x,\mu)
    :=\sqrt{\mu_j^2+g_j(x)^2}-\mu_j+g_j(x).
    \label{eq:fb-residual}
\end{equation}
Then $\Phi_j(x,\mu)=0$ is equivalent to $g_j(x)\leq0$, $\mu_j\geq0$, and $\mu_jg_j(x)=0$~\cite{Fischer1992,ChenKoPan2010}. The FB residual therefore provides a direct feedback signal for inequality multipliers, without introducing squared slack variables or specifying the active set in advance.
Using this representation, let $L(x,\mu):=f(x)+\mu^\top g(x)$ and consider the system formed by Lagrangian descent and complete PI feedback of the FB residual:
\begin{equation}
\begin{aligned}
    \dot x=-\nabla_xL(x,\mu),~
    \dot\mu=K_p\frac{\mathrm d}{\mathrm dt}\Phi(x,\mu)
    +K_i\Phi(x,\mu),
\end{aligned}
\label{eq:fb-pi-implicit}
\end{equation}
where $K_p\succeq0$ and $K_i\succ0$ are diagonal gain matrices. We refer to this system as the FB-PI dynamics. Its equilibrium conditions correspond to the Karush--Kuhn--Tucker (KKT) conditions of the original problem.

The dependence of $\Phi$ on both $x$ and $\mu$ makes the complete PI feedback implicit, while the ideal residual is nondifferentiable at weak complementarity points~\cite[Lemma~2.1]{ChenKoPan2010}. We first construct an explicit realization and then examine how feedback gains affect local stability, convergence rates, discrete implementation costs, and global convergence. The main results are as follows.

\noindent 1) We construct an explicit realization of complete FB-PI feedback through an invertible integral-coordinate transformation. It remains well-posed where the ideal residual is nondifferentiable, and its equilibria correspond exactly to KKT pairs.
2) We show that integral feedback can stabilize equilibria with an indefinite enhanced Hessian. For a single active constraint, we characterize the exact stability region and show that stabilizing integral gains exist if and only if (iff) the modified normal curvature is positive. We also quantify the gain and discretization costs of approaching the decay-rate ceiling. For multiple active constraints, we give a modal damping criterion for a uniform stability margin at large gains, showing how gain allocation and constraint geometry affect stability.
3) For strongly convex objectives with affine constraints, we establish global convergence conditions. Under these conditions and uniqueness of the KKT pair, the ideal FB-PI equilibrium is globally exponentially stable (GES) iff no nonzero nonnegative combination of constraint normals vanishes. If such cancellation exists, no fixed, finite gains yield a uniform global exponential estimate.

We further examine gain effects on convergence, compare integration costs during constraint switching, and assess behavior under weak complementarity. Numerical experiments support the theoretical findings.

\subsection{Related Work}
\label{subsec:related-work}

\noindent 1.~\textit{Multiplier Feedback and PI Dynamics.}
Projected primal--dual dynamics enforce multiplier nonnegativity \cite{FeijerPaganini2010}, while proximal augmented Lagrangians provide continuously differentiable saddle formulations \cite{DhingraKhongJovanovic2019}. Global exponential stability has been established for augmented primal--dual dynamics under suitable regularity assumptions \cite{QuLi2019}. Under more general constraint assumptions, semi-global exponential guarantees and counterexamples to global exponential stability have also been established \cite{TangQuLi2020}. PI structures also arise in consensus dynamics and projected conic methods \cite{DrogeKawashimaEgerstedt2014,YuElangoTopcuAcikmese2022}. The CMO framework treats multipliers as inputs and constraints as outputs, developing PI and feedback-linearization laws, including nonconvex examples and squared-slack inequality reformulations \cite{CeroneFossonPirreraRegruto2025TAC}.
For inequalities, augmented-Lagrangian PI (AL-PI) feedback combines integral feedback of the augmented dual residual with proportional constraint-velocity feedback, obtaining global exponential convergence under strong convexity and smoothness, full-row-rank affine constraints, and suitable feedback and augmentation parameters \cite{CeroneFossonPirreraRegruto2024CDC}. Antiwindup designs address multiplier saturation \cite{Adegbege2020Antiwindup}. Related proportional--integral--derivative formulations incorporate antiwindup PI feedback for inequalities and admit an equivalent preconditioned augmented-Lagrangian representation \cite{CentorrinoHoteitBaltaLygeros2026SPPID}. Under strict complementarity and matched gains, the AL-PI implementation in \cite{CeroneFossonPirreraRegruto2024CDC} shares the FB-PI active linearization but differs in inactive-channel coupling.

\noindent 2.~\textit{Complementarity-Based Dynamics.}
The FB function provides an equation-based representation of complementarity conditions, supporting Newton-type reformulations \cite{Fischer1992} and perturbed-equation continuation methods \cite{Kanzow1996}. Continuous-time neural models connect these reformulations with merit minimization \cite{LiaoQiQi2001}, including gradient flows of squared generalized FB residuals \cite{ChenKoPan2010}. Here, the FB residual drives complete PI feedback in the multiplier channel, while the primal channel follows ordinary Lagrangian descent. We study the stability and convergence rates of the resulting interconnection.

\section{Explicit Realization and Basic Properties}
\label{sec:fb-dynamics}

This section expresses the FB-PI feedback as an explicit dynamical system and establishes its well-posedness, equilibrium correspondence, and multiplier-channel properties.

Let $J(x)\mathrel{:=}\nabla g(x)\in\mathbb R^{m\times n}$ denote the constraint Jacobian. To treat the ideal dynamics and their smoothed implementation together, we use the standard smoothed FB residual~\cite{Kanzow1996}:
\begin{equation}
    \Phi_j^\varepsilon(x,\mu)
    \mathrel{:=}\sqrt{\mu_j^2+g_j(x)^2+2\varepsilon^2}-\mu_j+g_j(x),
    \label{eq:smoothed-fb-residual}
\end{equation}
where $\varepsilon\geq0$ and $\Phi^0=\Phi$. We replace $\Phi$ in \eqref{eq:fb-pi-implicit} with $\Phi^\varepsilon$ and fix diagonal gains $K_p=\operatorname{diag}(k_{p,j})\succeq0$, $K_i=\operatorname{diag}(k_{i,j})\succ0$.

\subsection{Integral Coordinates and Explicit Realization}\label{subsec:sys}

Where the residual is differentiable, the multiplier feedback equation can be written as
\(
    \frac{\mathrm d}{\mathrm dt}
    \bigl[\mu-K_p\Phi^\varepsilon(x,\mu)\bigr]
    =K_i\Phi^\varepsilon(x,\mu).
\)
The quantity in brackets accumulates the complementarity residual weighted by $K_i$. We therefore define the integral state
\begin{equation}
    \nu:=\mu-K_p\Phi^\varepsilon(x,\mu).
    \label{eq:integral-coordinate}
\end{equation}
Then $\dot\nu=K_i\Phi^\varepsilon(x,\mu)$. Recovering $\mu$ uniquely from the current state $(x,\nu)$ allows the primal and integral-state velocities to be evaluated directly.

The coordinate relation \eqref{eq:integral-coordinate} can be inverted channel by channel. For channel $j$, write $y:=g_j(x)$, $k:=k_{p,j}$, $c:=1+2k$, and $w:=\nu_j+ky$, and define
\[
    M_j^\varepsilon(y,\nu_j)
    :=\left({(1+k)w+k\sqrt{w^2+c(y^2+2\varepsilon^2)}}\right)/{c},
\]
\[
    P_j^\varepsilon(y,\nu_j)
    :=\left({\sqrt{w^2+c(y^2+2\varepsilon^2)}-\nu_j+(1+k)y}\right)/{c}.
\]
Here, $M_j^\varepsilon$ gives the recovered multiplier, and $P_j^\varepsilon$ gives the corresponding complementarity residual. When $k=0$, the formula gives $M_j^\varepsilon(y,\nu_j)=\nu_j$.

For every $\varepsilon\geq0$ and fixed $x$, the transformation \eqref{eq:integral-coordinate} is globally invertible with respect to $\mu$, with its inverse given componentwise by $M_j^\varepsilon$. Moreover, $M_j^\varepsilon$ and $P_j^\varepsilon$ are globally Lipschitz in $(y,\nu_j)$, with Lipschitz constants independent of $\varepsilon$. These properties also hold when $\varepsilon=0$, providing a common basis for realizing the ideal and smoothed systems; see Appendix~\ref{app-real-inversion} for the proof.

Combining the channels into vector maps $M^\varepsilon$ and $P^\varepsilon$ yields
\(
    M^\varepsilon(g(x),\nu)
    =\nu+K_pP^\varepsilon(g(x),\nu),
\)
\(
    P^\varepsilon(g(x),\nu)
    =\Phi^\varepsilon\!\left(x,M^\varepsilon(g(x),\nu)\right).
\)
The resulting explicit FB-PI system is
\begin{equation}
    \dot x=-\nabla f(x)-J(x)^\top M^\varepsilon(g(x),\nu),~
    \dot\nu=K_iP^\varepsilon(g(x),\nu).
    \label{eq:fb-pi-realization}
\end{equation}
Its dynamical state is $(x,\nu)$, and the physical multiplier is recovered as $\mu=M^\varepsilon(g(x),\nu)$. The right-hand side is fully determined by the current state.

Integrating the second equation of \eqref{eq:fb-pi-realization} and using \eqref{eq:integral-coordinate} gives
\begin{equation}
    \mu(t)-K_p\Phi^\varepsilon(t)
    =\mu(0)-K_p\Phi^\varepsilon(0)
    +\int_0^tK_i\Phi^\varepsilon(s)\,\mathrm ds,
    \label{eq:fb-pi-integral-identity}
\end{equation}
where $\Phi^\varepsilon(t):=\Phi^\varepsilon(x(t),\mu(t))$. Together with the primal equation and the coordinate relation \eqref{eq:integral-coordinate}, \eqref{eq:fb-pi-integral-identity} is equivalent to \eqref{eq:fb-pi-realization} for locally absolutely continuous trajectories. This integral formulation preserves complete PI feedback at singular complementarity points of the ideal residual without requiring pointwise differentiation.

\subsection{Equilibria and Channel Properties}
\label{subsec:basic-properties}

The explicit realization above allows the basic properties of the ideal and smoothed systems to be analyzed together.
\begin{proposition}[Basic properties of the explicit realization]
\label{prop:exact-realization}
Under the stated smoothness and gain assumptions, system \eqref{eq:fb-pi-realization} has the following properties for every $\varepsilon\geq0$.

\noindent 1)~\textit{Well-posedness.}
Every initial state generates a unique maximal solution. The recovered $\mu$ and $\Phi^\varepsilon$ are locally absolutely continuous and satisfy \eqref{eq:fb-pi-integral-identity}.

\noindent 2)~\textit{Equilibrium correspondence.}
Equilibria are characterized by
\(
\nabla f(x)+J(x)^\top\mu=0,~
\Phi^\varepsilon(x,\mu)=0,~
\nu=\mu,
\)
and correspond exactly to KKT pairs when $\varepsilon=0$.

\noindent 3)~\textit{Nonnegativity and inactive channels.}
The set $\nu\geq0$ is forward invariant, with $\mu\geq0$. For $\varepsilon=0$, a zero multiplier remains zero while its constraint remains strictly satisfied.
\end{proposition}
Detailed proofs are given in Appendix~\ref{app-real-main}. The Lipschitz bounds in Section~\ref{subsec:sys} and the smoothness of $f$ and $g$ make the vector field in \eqref{eq:fb-pi-realization} locally Lipschitz, giving local existence and uniqueness of solutions~\cite[Sec.~3.1]{Khalil2002}. The recovered trajectories are locally absolutely continuous by composition with $M^\varepsilon$ and $P^\varepsilon$. For affine constraints and a globally Lipschitz objective gradient, the vector field is globally Lipschitz, ensuring existence for all $t\geq0$.

At equilibrium, $K_i\succ0$ implies $\Phi^\varepsilon=0$ and hence $\nu=\mu$. Together with primal stationarity, these conditions are independent of the gains for fixed $\varepsilon$. For $\varepsilon>0$, the zero-residual condition is equivalent to
\(
\mu_j>0,~
g_j(x)<0,~
-\mu_jg_j(x)=\varepsilon^2,
\)
which gives the standard perturbed complementarity conditions~\cite{Kanzow1996}. Figure~\ref{fig:app-fb-residual-contours} illustrates the ideal and perturbed zero sets, and Appendix~\ref{app-approx-main} quantifies the resulting equilibrium displacement and trajectory perturbations.

At $\nu_j=0$, the inverse formulas give $P_j^\varepsilon(g_j(x),0)\geq0$, establishing forward invariance of $\nu\geq0$. They also give $M_j^\varepsilon(g_j(x),\nu_j)\geq0$ whenever $\nu_j\geq0$. For the ideal system, primal--dual feasible initial conditions imply $\Phi^0(x(0),\mu(0))\leq0$ and thus $\nu(0)\geq0$ by \eqref{eq:integral-coordinate}.

For $y<0$, the ideal inverse maps satisfy $M_j^0(y,0)=P_j^0(y,0)=0$. By uniqueness, adding a constraint that remains strictly satisfied along the original trajectory, with zero initial multiplier, leaves that trajectory unchanged when the original data, gains, and initial conditions are retained.

\section{Local Stability and Gain Effects}
\label{sec:local-stability}

This section analyzes the local stabilizing gains, achievable decay rates, and discrete-time implementation costs of the FB-PI system, and then considers integral-gain allocation with multiple active constraints. Throughout this section, we set $\varepsilon=0$ and study \eqref{eq:fb-pi-realization} near a reference KKT pair $(x^\star,\mu^\star)$.

\subsection{Local Structure and Curvature Modification}
\label{subsec:local-structure}

Define the active and inactive sets at the reference KKT pair as
\(
\mathcal A\mathrel{:=}\{j:g_j(x^\star)=0\},\;
\mathcal I\mathrel{:=}\{j:g_j(x^\star)<0\},
\)
and let $J_A$ and $J_I$ collect the corresponding rows of $J(x^\star)$. Write $H\mathrel{:=}\nabla_{xx}^2L(x^\star,\mu^\star)$. Assume strict complementarity, full row rank of $J_A$, and $d^\top Hd>0$ for every nonzero $d\in\ker J_A$. These conditions ensure differentiability of the ideal residual near the reference point and positive tangent curvature, while allowing negative normal curvature and tangent--normal coupling.

By the derivatives of FB~\cite[Lemma~2.1]{ChenKoPan2010}, the first-order residual variations are
\(
\delta\Phi_A=J_A\delta x,~
\delta\Phi_I=-\delta\mu_I.
\)
Let $\eta\mathrel{:=}\delta\mu_A-K_{p,A}J_A\delta x=\delta\nu_A$ and
\(
H_e\mathrel{:=}J_A^\top K_{p,A}J_A+H,~
D_I\mathrel{:=}(I+K_{p,I})^{-1}K_{i,I}.
\)
The linearized system is
{\small
\begin{equation}
\begin{bmatrix}
\delta\dot x\\
\dot\eta\\
\delta\dot\mu_I
\end{bmatrix}
=
\begin{bmatrix}
-H_e&-J_A^\top&-J_I^\top\\
K_{i,A}J_A&0&0\\
0&0&-D_I
\end{bmatrix}
\begin{bmatrix}
\delta x\\
\eta\\
\delta\mu_I
\end{bmatrix}.
\label{eq:local-linearization}
\end{equation}
}
Proportional feedback modifies the normal curvature through $J_A^\top K_{p,A}J_A$ while preserving the quadratic form of $H$ on $\ker J_A$. Integral feedback regulates the primal--multiplier coupling. If $H_e\succ0$, a quadratic energy estimate and full row rank of $J_A$ show that the active subsystem is Hurwitz for every positive diagonal $K_{i,A}$. The linearization theorem~\cite[Ch.~4]{Khalil2002} then gives local exponential stability of the full equilibrium. Under the tangent-space positivity assumption, sufficiently large active proportional gains ensure $H_e\succ0$; an explicit sufficient bound and the proof are given in Appendix~\ref{app-local-linearization}.

The block-triangular structure of \eqref{eq:local-linearization} separates the active spectrum from the inactive eigenvalues $-d_j$, where $d_j\mathrel{:=}k_{i,j}/(1+k_{p,j})$. Increasing an inactive proportional gain at fixed integral gain therefore slows that channel. Inactive multipliers affect the primal transient without changing the active eigenvalues.
The same active structure also appears in AL-PI feedback. For the constraint-velocity implementation in \cite[Eq.~(9)]{CeroneFossonPirreraRegruto2024CDC}, equal scalar integral gains and $k_p^{\mathrm{FB}}=\rho+k_p^{\mathrm{AL}}$ give identical active linearizations in the respective integral coordinates, where $\rho>0$ is the augmentation parameter. The inactive channels differ: FB-PI gives independent multiplier decay, whereas this AL-PI implementation satisfies
\(
\delta\dot\mu_I=-\frac{k_i}{\rho}\delta\mu_I+k_p^{\mathrm{AL}}J_I\delta\dot x,
\)
retaining constraint-velocity forcing; see Appendix~\ref{app-local-al}.

\subsection{Exact Stability Boundary}
\label{subsec:exact-stability}

The condition $H_e\succ0$ guarantees stability independently of the integral-gain magnitude. We now derive exact joint gain conditions for a single active constraint.

Write $J_A=\sigma n_A^\top$, where $\sigma>0$ and $\|n_A\|=1$, and let the columns of $Q$ form an orthonormal basis of the constraint tangent space. Define
\(
T\mathrel{:=}Q^\top HQ\succ0,\;
b\mathrel{:=}Q^\top Hn_A,\;
a\mathrel{:=}n_A^\top Hn_A,
\)
representing tangent curvature, tangent--normal coupling, and normal curvature, respectively. For active gains $k_p$ and $k_i$, define the modified normal curvature $e\mathrel{:=}a+k_p\sigma^2$ and scaled integral gain $q\mathrel{:=}k_i\sigma^2$. If there are no tangent variables, the corresponding blocks are omitted and quadratic forms involving $b$ are interpreted as zero.

\begin{theorem}[Exact stability boundary]
\label{thm:exact-stability}
For every $q>0$, there is a unique $t_q>0$ satisfying
\(
q=t_q\bigl[1+b^\top(T^2+t_qI)^{-1}b\bigr].
\)
Define $\beta(q)\mathrel{:=}b^\top T(T^2+t_qI)^{-1}b$. The reference KKT equilibrium is locally exponentially stable iff
\(
e>\beta(q).
\)

If $e<\beta(q)$, the equilibrium is unstable. If $e=\beta(q)$, the active linearization has a simple purely imaginary pair $\pm i\sqrt{t_q}$, and all its remaining eigenvalues lie in the open left half-plane.
If $b\ne0$, then $\beta(q)$ is strictly decreasing, with
\(
\lim_{q\downarrow0}\beta(q)=b^\top T^{-1}b,~
\lim_{q\to\infty}\beta(q)=0.
\)
If $b=0$, the stability condition reduces to $e>0$.
\end{theorem}

Schur elimination reduces the active characteristic equation to a scalar equation, whose imaginary-axis crossings determine the exact stability boundary. The crossing direction and stability for large proportional gains establish the stable side, while implicit differentiation gives the monotonicity of $\beta$; see Appendix~\ref{app-local-exact-stability}.

By the Schur complement~\cite{BoydVandenberghe2004}, $H_e\succ0$ is equivalent to $e>b^\top T^{-1}b$. For $b\ne0$, Theorem~\ref{thm:exact-stability} gives local exponential stability in the region
\(
\beta(q)<e<b^\top T^{-1}b,
\)
where $H_e$ is indefinite. Integral feedback therefore enlarges the stability region beyond the positive-definite enhanced-Hessian condition.

Since $\beta(q)\geq0$ and $\beta(q)\to0$ as $q\to\infty$, a stabilizing integral gain exists at fixed $k_p$ iff $e>0$. In particular, pure integral feedback ($k_p=0$) achieves local exponential stability for some $k_i>0$ iff $a>0$.
For example, consider
\[
\min_{x_1\le0}
\left\{
-x_1+\frac12x^\top
\begin{bmatrix}
1&2\\
2&1
\end{bmatrix}x
\right\}.
\]
At the KKT pair $x^\star=0$, $\mu^\star=1$, we have $T=1$, $b=2$, and $a=1$. With $k_p=0$, the enhanced Hessian is indefinite, yet local exponential stability holds iff $k_i>6$.
\subsection{Prescribed Decay Rates and Implementation Costs}
\label{subsec:prescribed-rate}

We next determine the gains required for a prescribed decay margin and the gain and discretization costs of approaching its maximum. Retain the single-active-constraint notation, let $M_q$ denote the active subsystem matrix in \eqref{eq:local-linearization}, and define
\(
\gamma(q)\mathrel{:=}-\max_{\lambda\in\operatorname{spec}(M_q)}
\operatorname{Re}\lambda.
\)
For a target rate $0<r<\lambda_{\min}(T)$, set $L_r\mathrel{:=}T-rI$, $B_r\mathrel{:=}b^\top L_r^{-1}b$, and
\[
\begin{aligned}
\beta_r(t)&\mathrel{:=}b^\top T(L_r^2+tI)^{-1}b,\\
h_r(t)&\mathrel{:=}(r^2+t)\bigl[1+b^\top(L_r^2+tI)^{-1}b\bigr].
\end{aligned}
\]
Define $e_r\mathrel{:=}2r+\beta_r(0)$, with the convention $\lambda_{\min}(T)=+\infty$ when there are no tangent variables.

\begin{theorem}[Prescribed decay rates and gain costs]
\label{thm:prescribed-rate}
The following statements hold.

\noindent 1)~If $b\ne0$, then $\gamma(q)>r$ iff $e>2r$ and $q>q_{\min}(r,e)$, where
\begin{equation}
q_{\min}(r,e)\mathrel{:=}
\begin{cases}
h_r(t_{r,e}),&2r<e<e_r,\\
r(e-r-B_r),&e\ge e_r.
\end{cases}
\label{eq:prescribed-rate-threshold}
\end{equation}
Here $t_{r,e}>0$ is the unique solution of $\beta_r(t_{r,e})=e-2r$ in the first branch. If $b=0$, the conditions reduce to $e>2r$ and $q>r(e-r)$.

\noindent 2)~If $b\ne0$ and $0<e<2\lambda_{\min}(T)$ is fixed, then $\sup_{q>0}\gamma(q)=e/2$, while $\gamma(q)<e/2$ for every finite $q>0$. With $c_b\mathrel{:=}b^\top Tb>0$,
\begin{equation}
\begin{aligned}
\gamma(q)=\frac e2-\frac{c_b}{2q}+O(q^{-2}),~
q\to\infty,\\
q_{\min}\!\left(\frac e2-\delta,e\right)
=\frac{c_b}{2\delta}+O(1),~
\delta\downarrow0.
\end{aligned}
\label{eq:rate-gain-cost}
\end{equation}
\end{theorem}

The branches in \eqref{eq:prescribed-rate-threshold} correspond to different rate-limiting modes: the conjugate roots $-r\pm i\sqrt{t_{r,e}}$ when $2r<e<e_r$, and the real root $-r$ when $e>e_r$, with the boundaries meeting at $e=e_r$. In the second branch, the required integral gain increases linearly with $e$, so stronger proportional feedback may require a larger integral gain to maintain the same target rate. The proof determines these thresholds by counting characteristic-root crossings of $\operatorname{Re}s=-r$. Large-gain root expansions give \eqref{eq:rate-gain-cost}, while root counting at the limiting boundary establishes nonattainment of the ceiling at finite gains; see Appendices~\ref{app-rate-high-gain}--\ref{app-rate-costs}.

Under the assumptions of part 2), nonzero tangent--normal coupling imposes the ceiling $e/2$, and reducing the gap to this ceiling to $\delta$ requires an integral gain of order $\delta^{-1}$. When $b=0$ and $e>0$, the normal characteristic polynomial is $s^2+es+q$, whose spectral decay margin reaches $e/2$ for $q\ge e^2/4$. The full nonlinear system also depends on the tangent and inactive modes: if $\gamma(q)>r$ and $d_j>r$ for every $j\in\mathcal I$, a local exponential estimate with rate $r$ holds, with a gain-dependent prefactor and neighborhood.

The same fast modes that determine the large-gain decay margin also restrict explicit discretization. At a locally exponentially stable reference point, the explicit Euler stability condition $|1+\alpha\lambda|<1$~\cite[Sec.~IV.2]{HairerWanner1996} gives Schur stability of the linearization iff $0<\alpha<\alpha_{\max}^{\mathrm{full}}$, where
\begin{equation}
\alpha_{\max}^{\mathrm{full}}\mathrel{:=}
\min\left\{
\min_{\lambda\in\operatorname{spec}(M_q)}
\frac{-2\operatorname{Re}\lambda}{|\lambda|^2},
\;
\min_{j\in\mathcal I}\frac{2}{d_j}
\right\}.
\label{eq:euler-local-bound}
\end{equation}
The second minimum is $+\infty$ when there are no inactive constraints. Under the assumptions of part 2), the fast conjugate pair satisfies $s_\pm(q)=-e/2\pm i\sqrt q+O(q^{-1/2})$, giving $\alpha_{\max}^{\mathcal A}\sim e/q$ for the active subsystem. For $r=e/2-\delta$, choose $q>q_{\min}(r,e)$ with $q/q_{\min}(r,e)\to1$. Then
\(
\alpha_{\max}^{\mathcal A}\sim 2e\delta/c_b
\)
as $\delta\downarrow0$, and fixed inactive-channel gains preserve this scaling for the full system. Thus, approaching the continuous-time rate ceiling requires integral gains growing as $\delta^{-1}$ and admissible Euler step sizes decreasing proportionally to $\delta$. Appendix~\ref{app-approx-euler} gives the nonlinear local Euler argument and distinguishes discretization in physical and integral coordinates.
\subsection{Integral-Gain Allocation with Multiple Active Constraints}
\label{subsec:multiple-active}

With multiple active constraints, integral-gain allocation changes the oscillation directions and hence the curvature experienced by individual modes. Fix the proportional and inactive integral gains, and set $K_{i,A}\mathrel{:=}\kappa K_0$, where $K_0\succ0$ is fixed and diagonal and $\kappa>0$ is the overall gain scale. Let $U$ and $N$ have orthonormal columns spanning $\ker J_A$ and its orthogonal complement, respectively, and define
\[
R\mathrel{:=}J_AN,\;
T\mathrel{:=}U^\top HU,\;
E\mathrel{:=}N^\top H_eN,\;
Q_0\mathrel{:=}R^\top K_0R.
\]
Full row rank of $J_A$ makes $R$ nonsingular and $Q_0\succ0$. Here $T$ and $E$ are the tangent and modified normal curvature blocks, while $Q_0$ combines constraint geometry with integral-gain allocation.

Let $\lambda_\ell$ be the distinct eigenvalues of $Q_0$, with orthogonal spectral projectors $\Pi_\ell$, and define
\(
E_\ell\mathrel{:=}
\left.\Pi_\ell E\Pi_\ell\right|_{\operatorname{range}\Pi_\ell}.
\)
As $\kappa\to\infty$, the active eigenvalues satisfy, counting algebraic multiplicity,
\[
\begin{aligned}
s_j^0(\kappa)
&=-\lambda_j(T)+O(\kappa^{-1/2}),\\
s_{\ell,j}^{\pm}(\kappa)
&=\pm i\sqrt{\kappa\lambda_\ell}
-\frac12\lambda_j(E_\ell)+O(\kappa^{-1/2}).
\end{aligned}
\]
The first family is omitted when there are no tangent variables. Scaling $\eta$ by $\kappa^{-1/2}K_0^{-1/2}$ makes the leading system matrix skew-symmetric. Semisimple eigenvalue perturbation theory on its frequency subspaces~\cite[Ch.~II, \S~2.3]{Kato1995} yields these expansions, including repeated frequencies. Related stabilization mechanisms involving rotation and indefinite damping are studied in \cite{CrauelDammIlchmann2007,Kirillov2012}.

The gain scale $\kappa$ determines the fast frequencies, while $E_\ell$ determines the limiting modal damping. A strictly positive spectral decay margin, uniform for all sufficiently large $\kappa$, exists iff $E_\ell\succ0$ for every $\ell$. A negative eigenvalue of any $E_\ell$ implies instability for all sufficiently large $\kappa$. If every $E_\ell\succeq0$ and at least one is singular, no uniform positive margin exists, and eventual stability depends on higher-order terms.

For simple eigenvalues of $Q_0$ with unit eigenvectors $v_\ell$, the criterion becomes $v_\ell^\top Ev_\ell>0$ for every $\ell$, which can hold even when $E$ is indefinite. If $Q_0$ is a positive scalar multiple of the identity, the criterion instead requires $E\succ0$. Integral-gain allocation therefore changes the damping condition by separating or merging modal frequencies. The attainable modal directions are restricted by $Q_0=R^\top K_0R$ with diagonal $K_0$. The complete proof, a sufficient finite-gain condition, and a two-channel example are given in Appendix~\ref{app-modal-analysis}.

\section{Nonlinear Dissipation and Convergence}
\label{sec:global-convergence}

The local analysis characterizes stability and decay rates near a KKT pair. We now study convergence from arbitrary initial conditions and determine when fixed gains admit a uniform global exponential estimate. Consider affine constraints $g(x)=Cx-d$ and an $\ell$-strongly convex objective $f$ with a globally Lipschitz gradient. For fixed gains $K_p\succeq0$ and $K_i\succ0$, we combine an FB-residual dissipation estimate with conditions on the constraint geometry to control the full state.

\subsection{Dissipation of the FB Residual}
\label{subsec:global-dissipation}

First consider the scalar residual $\psi\mathrel{:=}\sqrt{\mu^2+g^2+2\varepsilon^2}-\mu+g$, where $\varepsilon\geq0$. Let the reference values satisfy $\bar\mu\geq0$, $\bar g\leq0$, and $-\bar\mu\bar g=\varepsilon^2$. For arbitrary current values $(\mu,g)$,
\begin{equation}
(\mu-\bar\mu)\bigl[\psi-(g-\bar g)\bigr]
\leq(g-\bar g)\psi-\frac12\psi^2.
\label{eq:shifted-fb-inequality}
\end{equation}
This inequality follows from the algebraic identity defining the residual and the Cauchy--Schwarz inequality. It requires no differentiation of the residual and therefore also applies at infeasible states, negative multipliers, and singular complementarity points of the ideal residual. The derivation is given in Appendix~\ref{app-global-scalar}.

Now consider an equilibrium $(\bar x,\bar\mu)$. Since $\bar\nu=\bar\mu$, define $u\mathrel{:=}x-\bar x$, $v\mathrel{:=}\nu-\bar\mu$, and $\psi\mathrel{:=}P^\varepsilon(Cx-d,\nu)$. The exact realization gives $\mu-\bar\mu=v+K_p\psi$. For the energy function $V\mathrel{:=}\|u\|^2/2+v^\top K_i^{-1}v/2$, summing \eqref{eq:shifted-fb-inequality} over the channels and using strong convexity yields
\begin{equation}
\dot V\leq-\ell\|u\|^2+(Cu)^\top\psi
-\psi^\top\left(K_p+\frac12I\right)\psi.
\label{eq:global-dissipation}
\end{equation}
If the proportional gains satisfy
\begin{equation}
\ell I-C^\top\operatorname{diag}
\left(\frac{1}{4k_{p,j}+2}\right)C\succ0,
\label{eq:global-gain-condition}
\end{equation}
the Schur complement gives a constant $c_0>0$ such that
\[
\dot V\leq-c_0\bigl(\|u\|^2+\|\psi\|^2\bigr).
\]
For a uniform proportional gain $K_p=kI$, a direct sufficient choice is $k\geq0$ and $k>\|C\|^2/(4\ell)-1/2$.

Equation \eqref{eq:global-dissipation} shows that proportional feedback increases dissipation through the squared residual. This estimate directly controls the primal deviation and complementarity residual, but not yet the entire integral state. To exploit the coupling between the primal and integral states, choose a sufficiently small $\eta>0$ and define $W\mathrel{:=}V+\eta u^\top C^\top v$. Using the Lipschitz continuity of the objective gradient to absorb the cross terms makes $W$ equivalent to $\|u\|^2+\|v\|^2$ and gives
\[
\dot W\leq-c_1\bigl(\|u\|^2+\|\psi\|^2+\|C^\top v\|^2\bigr)
\]
for some $c_1>0$. The comparison and dissipation constants can all be chosen independently of $\varepsilon$. The construction is detailed in Appendix~\ref{app-global-modified-energy}.
Thus, global exponential convergence depends on whether the primal deviation, complementarity residual, and net constraint feedback $C^\top v$ jointly control the entire integral-state deviation.

\subsection{Uniform Exponential Convergence}
\label{subsec:global-uniform}

If $C$ has full row rank, the net constraint feedback directly controls the entire integral state through $\|C^\top v\|\geq\sigma_{\min}(C)\|v\|$. Substitution into the modified energy estimate gives exponential decay of the full state.

Specifically, under \eqref{eq:global-gain-condition}, the system has a unique equilibrium $(x_\varepsilon^\star,\mu_\varepsilon^\star)$ for every $\varepsilon\geq0$. Define the integral-coordinate state $w\mathrel{:=}(x,\nu)$ and equilibrium $w_\varepsilon^\star\mathrel{:=}(x_\varepsilon^\star,\mu_\varepsilon^\star)$. There exist constants $A\geq1$ and $c>0$ such that
\begin{equation}
\|w(t)-w_\varepsilon^\star\|
\leq Ae^{-ct}\|w(0)-w_\varepsilon^\star\|,
\;t\geq0,
\label{eq:smoothing-uniform-ges}
\end{equation}
where $A$ and $c$ can be chosen independently of the smoothing parameter and initial condition.

Equilibrium existence also follows from full row rank and strong convexity. Full row rank makes $Cx=d-\mathbf1$ solvable, so the constraints are strictly feasible. For the ideal problem, strong convexity gives a unique primal optimizer, Slater's condition ensures the existence of a KKT multiplier~\cite[Sec.~5.2.3]{BoydVandenberghe2004}, and injectivity of $C^\top$ gives multiplier uniqueness. For $\varepsilon>0$, the perturbed KKT equations correspond to a strongly convex logarithmic-barrier problem on the strictly feasible region~\cite[Sec.~11.2.2]{BoydVandenberghe2004}, whose unique solution gives the corresponding equilibrium. Existence and estimate~\eqref{eq:smoothing-uniform-ges} are proved in Appendices~\ref{app-global-equilibrium}--\ref{app-global-uniform}.

Estimate \eqref{eq:smoothing-uniform-ges} holds for arbitrary initial conditions and does not require strict complementarity. The uniform Lipschitz properties of the exact coordinate transformation and its inverse give the same exponential rate in physical coordinates $(x,\mu)$, with an adjusted prefactor. Uniformity is with respect to the smoothing parameter and initial condition; the constants may depend on the problem data and fixed feedback gains. Appendix~\ref{app-approx-practical} gives a complementary practical estimate relative to the ideal equilibrium.

Thus, full row rank allows the net constraint feedback to control every integral direction, yielding a uniform global exponential convergence estimate for both the ideal and smoothed systems.

\subsection{Constraint Geometry of the Ideal System}
\label{subsec:global-geometry}

For the ideal residual, the complementarity residual itself also controls some integral directions in $\ker C^\top$, allowing the full-row-rank condition to be weakened. We now set $\varepsilon=0$, remove the full-row-rank requirement, and assume that the problem has a unique KKT pair $(x^\star,\mu^\star)$. Strong convexity guarantees uniqueness of the primal optimizer; multiplier uniqueness is assumed separately.

\begin{theorem}[Global convergence and constraint geometry]
\label{thm:global-geometry}
Suppose that $f$ is strongly convex with a globally Lipschitz gradient, the constraints are $Cx-d\leq0$, and a unique KKT pair exists. For fixed diagonal gains $K_i\succ0$ and $K_p\succeq0$ satisfying \eqref{eq:global-gain-condition}, the ideal FB-PI equilibrium is globally asymptotically stable. Moreover, it is GES iff
\begin{equation}
\ker C^\top\cap\mathbb R_+^m=\{0\}.
\label{eq:nonnegative-kernel-condition}
\end{equation}
If \eqref{eq:nonnegative-kernel-condition} fails, no fixed, finite diagonal gains $K_i\succ0$ and $K_p\succeq0$ can make the equilibrium GES. This impossibility conclusion does not require \eqref{eq:global-gain-condition}.
\end{theorem}

Condition \eqref{eq:nonnegative-kernel-condition} allows linearly dependent constraint rows and more constraints than primal variables. It excludes cancellation of the constraint normals along a nonzero nonnegative multiplier direction: there is no $a\geq0$, $a\ne0$, such that $C^\top a=0$. By the theorem of alternatives for strict linear inequalities~\cite[Example~2.21]{BoydVandenberghe2004}, this condition is equivalent to the existence of a direction $h$ satisfying $Ch<0$, along which all constraint values strictly decrease. It is stronger than Slater's condition, which only requires a strictly feasible point.

Global asymptotic stability follows from the dissipation estimate and LaSalle's invariance principle~\cite[Sec.~4.2]{Khalil2002}. Condition \eqref{eq:global-gain-condition} confines trajectories to bounded energy sublevel sets. A trajectory remaining in the zero-dissipation set must have its primal state at the optimizer and its complementarity residual equal to zero; maintaining that state also requires stationarity. Uniqueness of the KKT pair therefore makes the reference equilibrium the largest invariant subset of the zero-dissipation set.

Exponential convergence further requires extending this dissipation to the entire integral state. Define $u\mathrel{:=}x-x^\star$, $v\mathrel{:=}\nu-\mu^\star$, and $\psi\mathrel{:=}P^0(Cx-d,\nu)$. Under uniqueness of the KKT pair and \eqref{eq:nonnegative-kernel-condition}, for every fixed $K_p\succeq0$ there exists $K>0$ such that every state satisfies the global error bound
\begin{equation}
\|v\|\leq K\bigl(\|u\|+\|\psi\|+\|C^\top v\|\bigr).
\label{eq:global-residual-bound}
\end{equation}
This bound is established by a normalization and contradiction argument. Violating sequences at finite or vanishing scales produce another KKT multiplier, whereas sequences at unbounded scales produce a nonzero nonnegative kernel vector through positive homogeneity of the ideal FB residual. Multiplier uniqueness and \eqref{eq:nonnegative-kernel-condition} exclude these possibilities, respectively. Substituting \eqref{eq:global-residual-bound} into the modified energy estimate gives exponential decay of the full state. This argument requires neither independence of the active constraint gradients nor strict complementarity. The error bound and convergence proof are given in Appendices~\ref{app-global-error-bound}--\ref{app-global-sufficiency}.

Necessity of \eqref{eq:nonnegative-kernel-condition} reveals a convergence limitation at large multiplier deviations. Suppose that a nonzero $a\geq0$ satisfies $C^\top a=0$, and consider the initial-state family
\[
x(0)=x^\star,\;
\nu(0)=\mu^\star+Ma,\;
M\to\infty.
\]
The initial state error grows linearly with $M$, while the additional integral state $Ma$ has zero net constraint effect. Both the FB residual and the vector field evaluated at these initial states remain bounded independently of $M$. Global Lipschitz continuity of the realized vector field and Gronwall's inequality~\cite{Khalil2002} then bound the displacement from the initial state over any fixed time interval independently of $M$. Consequently, the system cannot reduce an arbitrarily large initial error by a uniform fraction within a fixed time, ruling out an exponential estimate valid for all initial conditions. This argument uses only affine constraints, a reference KKT pair, and a globally Lipschitz objective gradient, so the same obstruction also applies beyond strongly convex objectives; see Appendix~\ref{app-global-necessity}.

These initial states lie in the nonnegative integral-state cone, and their recovered physical multipliers are also nonnegative. The obstruction arises from cancellation of the constraint normals and the residual behavior at large multipliers. It is compatible with asymptotic convergence of every trajectory and does not exclude local exponential stability near the equilibrium.
We provide two examples to illustrate the distinction between full row rank and \eqref{eq:nonnegative-kernel-condition}. First, consider
\[
\min_{x\in\mathbb R^2}~\frac12\|x-(1,1)^\top\|^2~
\text{s.t.}~x_1\leq0,\;x_2\leq0,\;x_1+x_2\leq-1.
\]
Its unique KKT pair is $x^\star=(-1/2,-1/2)^\top$, $\mu^\star=(0,0,3/2)^\top$. The constraint matrix has rows $(1,0)$, $(0,1)$, and $(1,1)$, so it does not have full row rank. Nevertheless, $h=(-1,-1)^\top$ satisfies $Ch<0$. With $K_p=I$, the matrix in \eqref{eq:global-gain-condition} has minimum eigenvalue $1/2$. Hence, the system is GES for every fixed positive diagonal integral gain $K_i$.

Next, consider the interval constrained problem $\min \frac12x^2$
\(
\text{s.t.}~{-1\leq x\leq1}.
\)
Here, $C=(1,-1)^\top$, the unique KKT pair is $x^\star=0$, $\mu^\star=0$, and $(1,1)^\top\in\ker C^\top\cap\mathbb R_+^2$. With $K_p=I$, condition \eqref{eq:global-gain-condition} holds and the system is globally asymptotically stable. However, no fixed, finite feedback gains can make it GES.
More generally, a nonempty bounded polyhedron admits no direction satisfying $Ch<0$, which would generate an unbounded feasible ray. For such problems with a unique KKT pair, \eqref{eq:global-gain-condition} guarantees global asymptotic stability of the ideal FB-PI system, but no fixed, finite diagonal gains can make it GES. This geometric obstruction persists even when the equilibrium is locally exponentially stable; the geometric equivalence is proved in Appendix~\ref{app-global-polyhedron}.

\section{Numerical Experiments}
\label{sec:numerical-experiments}

We examine how feedback gains shape local dynamics, how the FB channel structure affects numerical implementation, and how inequality representations behave near weak complementarity. Accuracy is measured by the original KKT residual
\[
\begin{aligned}
R_{\mathrm{KKT}}(x,\mu)\mathrel{:=}\bigl(
&\|\nabla f(x)+J(x)^\top\mu\|^2+\|[g(x)]_+\|^2\\
&+\|[-\mu]_+\|^2+\|\mu\odot g(x)\|^2
\bigr)^{1/2},
\end{aligned}
\]
together with the physical-state error
\(E(x,\mu)\mathrel{:=}\|(x,\mu)-(x^\star,\mu^\star)\|\).
The local nonconvex experiment uses the ideal FB residual; the other FB tests in this section use $\varepsilon=10^{-6}$. Full experimental settings, parameter derivations, and accuracy checks are provided in Appendix~\ref{app:experiments}.

\subsection{Integral Feedback, Stability, and Rate Costs}
\label{subsec:exp-local}

Consider the nonconvex problem
\[
\min_{x_1\leq0}\;-x_1+\frac12x^\top
\begin{bmatrix}1&2\\2&1\end{bmatrix}x,
\]
with reference KKT point $x^\star=0$, $\mu^\star=1$. The Hessian is positive definite on the constraint tangent space but indefinite on the full space. Setting $k_p=0$ leaves the enhanced Hessian indefinite. Nevertheless, with $q\mathrel{:=}k_i$, Theorem~\ref{thm:exact-stability} predicts local exponential stability iff $q>6$.

\begin{figure*}[t]
    \centering
    \includegraphics[width=0.8\textwidth]{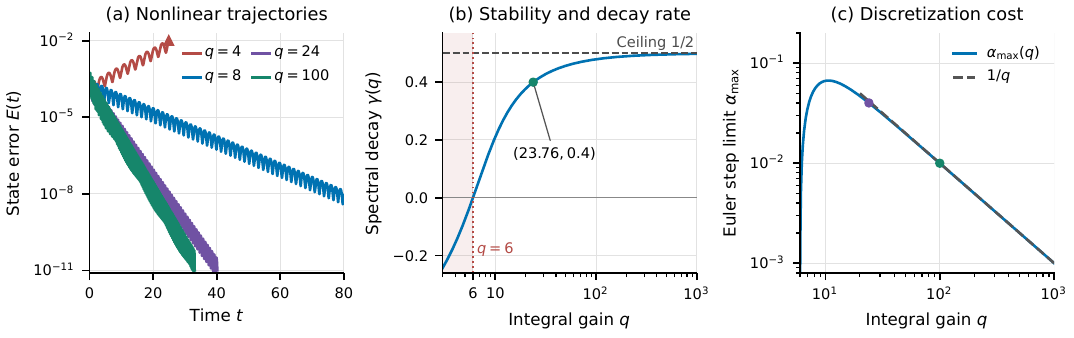}
    \caption{Gain-dependent local dynamics: (a) nonlinear state errors, (b) spectral decay margins, and (c) Euler stability limits.}
    \label{fig:note-local-dynamics}
\end{figure*}

Figure~\ref{fig:note-local-dynamics}(a) shows the corresponding nonlinear behavior from a common nearby initial condition. The error grows for $q=4$ and decays for $q=8,24,100$. Integral feedback therefore stabilizes the coupled primal--multiplier dynamics without restoring positive definiteness of the enhanced Hessian. A separate comparison with positive proportional gains and matched active linearizations gives similar local convergence for FB-PI and AL-PI, consistent with their shared active structure; see Appendix~\ref{app-exp-local}.

Figure~\ref{fig:note-local-dynamics}(b)--(c) illustrates the cost of increasing the decay margin. The exact threshold for the target margin $r=0.4$ is $q=23.76$. Increasing $q$ from $24$ to $100$ raises the margin from $0.4012$ to $0.4792$, while reducing the Euler stability limit from $0.04004$ to $0.009983$. As predicted by Theorem~\ref{thm:prescribed-rate}, the margin approaches $1/2$, whereas the admissible step size decreases as $q^{-1}$. Nonlinear Euler tests on either side of the predicted step-size boundary likewise exhibit error decay and departure from the local neighborhood. Thus, approaching the continuous-time rate ceiling introduces increasingly rapid oscillations, restricting the step size available for explicit implementation.

\subsection{FB Channel Structure and Numerical Implementation}
\label{subsec:exp-channels}

\paragraph{Inactive channels and discrete convergence}
The local analysis identifies two distinct effects of complete FB proportional feedback: it increases normal curvature in active channels and increases the multiplier time constant in inactive channels. Since a uniform explicit step size must accommodate all modes, these effects jointly influence discrete convergence.

Consider the strongly convex quadratic program
\[
\begin{aligned}
\min_{x\in\mathbb R^3}\;&
\frac12x_1^2+x_2^2+\frac12x_3^2-\frac25x_1-x_2,\\
\text{s.t.}\;&
\frac25x_1\leq0,\quad x_2\leq0,\quad x_3-1\leq0.
\end{aligned}
\]
Its unique KKT point is $x^\star=0$, $\mu^\star=(1,1,0)^\top$. For this problem, FB-PI and AL-PI have identical active characteristic polynomials when their integral gains agree and
\(k_p^{\mathrm{FB}}=k_p^{\mathrm{AL}}+\rho\).
Their inactive multiplier eigenvalues, however, are
\(-k_i/(1+k_p^{\mathrm{FB}})\) and \(-k_i/\rho\), respectively. The inactive channel therefore imposes different restrictions on the full discrete spectrum.

Let $A$ denote the ideal Jacobian and
\(r_{\mathrm E}\mathrel{:=}\operatorname{spr}(I+\alpha A)\).
We independently minimize this spectral factor for each method within the class of shared scalar gains and a uniform Euler step size in physical coordinates, including optimization of $\rho$ for AL-PI. The Schur conditions for the modal characteristic polynomials give the exact minima
\(
r_{\mathrm{FB}}^\star=\frac37,
r_{\mathrm{AL}}^\star=\frac{1+2\sqrt2}{7}.
\)
The derivation and attaining parameters are given in Appendix~\ref{app:exp-discrete}.
Using the resulting parameters in the full nonlinear dynamics, $20$ common random directions at each perturbation radius give median iteration counts of $21$ versus $26$ at radius $10^{-2}$ and $23.5$ versus $30$ at radius $10^{-1}$, with stopping condition $R_{\mathrm{KKT}}\leq10^{-8}$. In this example, the FB inactive-channel time scale permits a more favorable balance of discrete modal contraction, and the nonlinear iterations follow the spectral prediction.

\paragraph{Feedback smoothness and integration work}
The feedback structures also differ during constraint transitions. The ideal FB residual is smooth away from $(\mu_j,g_j)=(0,0)$, whereas the AL positive-part feedback changes branch on $\mu_j+\rho g_j=0$. Constraint switching can therefore occur without encountering a nonsmooth FB point.

We examine this effect on sparse strongly convex quadratic programs
\[
\min_x\;\frac12x^\top Hx-(C^\top\mu^\star)^\top x,
\quad \text{s.t.}\quad Cx-d\leq0,
\]
where $C=[I_m\;0]$ and $H$ combines a positive diagonal matrix with a weighted ring-graph Laplacian. Three independent instances use $n=1000$, $m=500$, and half the constraints active at $x^\star=0$. Here, the ideal active linearizations are matched. Both methods use DOP853 with identical integration tolerances and stopping condition $R_{\mathrm{KKT}}\leq10^{-6}$; AL-PI uses the lowest evaluation count among four matching parameter choices for each instance.

\begin{figure*}[t]
    \centering
    \includegraphics[width=0.8\textwidth]{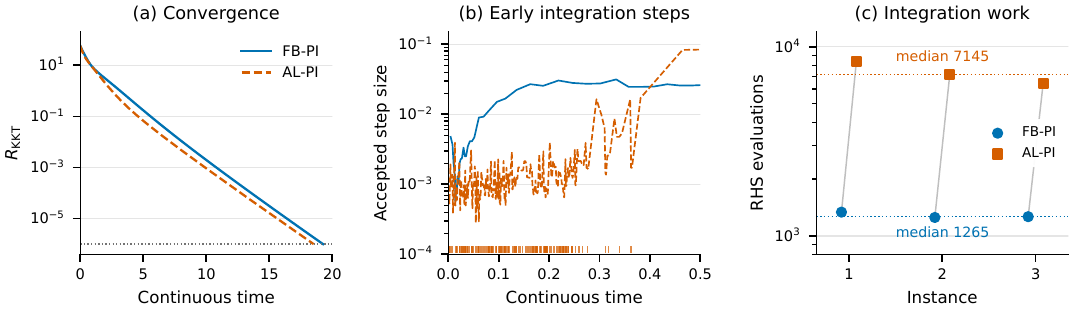}
    \caption{Constraint switching: (a) KKT residuals, (b) accepted step sizes and detected AL branch changes, and (c) function evaluations across three instances. Panels (a)--(b) show the first instance.}
    \label{fig:note-switch-cost}
\end{figure*}

At relative tolerance $10^{-9}$, the median function-evaluation counts are $1265$ for FB-PI and $7145$ for AL-PI. Their continuous stopping times remain close: $18.70$--$19.27$ for FB-PI and $18.38$--$18.52$ for AL-PI. Figure~\ref{fig:note-switch-cost} shows how the work difference develops. During the first $0.5$ time units of the first instance, FB-PI takes $42$ accepted steps and AL-PI takes $262$. Small AL steps concentrate near detected feedback branch changes, while the sampled FB complementarity pairs remain away from the ideal residual's singular point. These observations are consistent with smooth feedback facilitating adaptive integration. The benefit here is reduced integration work at comparable continuous convergence times. Appendix~\ref{app-exp-switching} gives the parameter choices and tighter-tolerance accuracy checks.

\subsection{Weak Complementarity and Inequality Forms}
\label{subsec:exp-weak}

We extend the sparse quadratic-program comparison to $n=5000$ and $m=2500$, including projected primal--dual dynamics (PD) and the squared-slack feedback-linearization implementation of CMO (FL-CMO). Multiplying the active reference multipliers and the corresponding objective linear term by $\theta\in\{1,0\}$ gives strict and weak complementarity while preserving the Hessian, constraints, and primal solution.

Both PI methods use modal parameter design, while PD and FL-CMO are tuned on the same independent smaller training instances. All parameters are frozen before testing; the design and training procedures are specified in Appendices~\ref{app:exp-modal-design} and~\ref{app:exp-large-weak}. The methods share primal initial conditions; FB-PI, AL-PI, and PD also share multiplier initial conditions, whereas FL-CMO determines its multiplier algebraically from the primal and slack states. All methods use explicit Euler updates in their respective coordinates, with multiplier projection for PD.

Table~\ref{tab:large-scale-complementarity} reports median iteration counts over three test instances at the common threshold $R_{\mathrm{KKT}}\leq10^{-6}$. FL-CMO requires the fewest iterations under strict complementarity, but its iteration count increases substantially when the active multipliers vanish. Its slack dynamics,
\(\dot s=-2\mu\odot s\),
lose linear damping at weakly complementary active constraints, consistent with this slowdown. FB-PI and AL-PI introduce no additional slack state and maintain short, similar iteration counts across both settings. Their iteration counts are also below those of PD in these tests.

\begin{table}[t]
    \centering
    \caption{Large-scale convergence under strict and weak complementarity.}
    \label{tab:large-scale-complementarity}
    \footnotesize
    \setlength{\tabcolsep}{2pt}
    \begin{tabular}{lrrrr}
        \hline
        Setting & FB-PI & AL-PI & PD & FL-CMO\\
        \hline
        Iterations, $\theta=1$ & 21 & 20 & 63 & 18\\
        Iterations, $\theta=0$ & 19 & 18 & 63 & 3891\\
        $\max E$, $\theta=0$
        & $5.23{\times}10^{-5}$
        & $9.23{\times}10^{-7}$
        & $2.81{\times}10^{-7}$
        & $4.16{\times}10^{-3}$\\
        \hline
    \end{tabular}
\end{table}

Smoothing also introduces a quantifiable equilibrium displacement in the weakly complementary case. The $1250$ active constraints satisfy $g_j(x)=x_j$ and $x_j^\star=\mu_j^\star=0$. At a smoothed equilibrium, perturbed complementarity gives
\(-\mu_j^\varepsilon x_j^\varepsilon=\varepsilon^2\).
Consequently, its physical-state error $E_\varepsilon\mathrel{:=}E(x^\varepsilon,\mu^\varepsilon)$ satisfies
\(
E_\varepsilon\geq\sqrt{2|\mathcal A|}\,\varepsilon
=5\times10^{-5}.
\)
Independent equilibrium solves place the terminal states within $6.84\times10^{-7}$ of the smoothed equilibria, showing that equilibrium displacement accounts for the observed state-error scale in this test.
Appendices~\ref{app:exp-weak}--\ref{app:exp-robustness} provide the scalar weak-complementarity analysis, smoothed-equilibrium verification, and additional tests across initial perturbations, problem sizes, and parameter choices. Appendix~\ref{app:exp-original-benchmarks} additionally reports constrained Rosenbrock and Griewank tests.

\section{Conclusion}
\label{sec:conclusion}

This paper studied the stability and convergence of PI multiplier dynamics driven by the Fischer--Burmeister residual. Integral feedback can achieve local stability with an indefinite enhanced Hessian through the dynamic coupling between primal variables and multipliers. Beyond stability, gain selection shapes decay rates and oscillation frequencies, leading to a tradeoff between continuous-time performance and admissible explicit Euler step sizes. For strongly convex objectives with affine constraints and a unique KKT pair, suitable gains provide global asymptotic convergence through residual dissipation, while the geometry of the constraint normals determines whether a global exponential estimate holds uniformly over all initial conditions. The exact gain characterization primarily concerns a single active constraint, and a general finite-gain stability characterization with multiple active constraints requires further study. Extending the global convergence analysis to nonconvex objectives and nonlinear constraints is another direction for future work.

\section*{References}
\bibliographystyle{IEEEtran}
\bibliography{reference}
\ifdefined\useRomanappendicesfalse\useRomanappendicesfalse\fi
\appendices
\numberwithin{equation}{section}
\numberwithin{figure}{section}
\numberwithin{table}{section}
\renewcommand{\thesubsection}{\thesection.\arabic{subsection}}
\ifdefined\thesubsectiondis
\renewcommand{\thesubsectiondis}{\thesubsection.}
\fi

\section{Exact Integral Realization}
\label{app-real-main}

This appendix establishes the coordinate properties in Section~\ref{subsec:sys} and proves Proposition~\ref{prop:exact-realization}, including the statements at singular pairs of the ideal residual. All gains are fixed, with $K_p\succeq0$ and $K_i\succ0$ diagonal. The symbols $y$ and $\nu$ in the scalar calculations denote one constraint value and its integral state.

\subsection{Global Inversion and Uniform Lipschitz Bounds}
\label{app-real-inversion}

For $k\geq0$ and $\varepsilon\geq0$, define
\[
    T_y^\varepsilon(\mu)
    \mathrel{:=}(1+k)\mu-k\sqrt{\mu^2+y^2+2\varepsilon^2}-ky.
\]
The reverse triangle inequality gives, whenever $\mu>\widetilde\mu$,
\[
    \mu-\widetilde\mu
    \leq T_y^\varepsilon(\mu)-T_y^\varepsilon(\widetilde\mu)
    \leq(1+2k)(\mu-\widetilde\mu).
\]
Thus $T_y^\varepsilon$ is continuous, strictly increasing, and onto $\mathbb R$. Its inverse exists globally, including when $\varepsilon=y=0$.

Set $c\mathrel{:=}1+2k$ and $w\mathrel{:=}\nu+ky$. The equation $T_y^\varepsilon(\mu)=\nu$ is equivalent to
\[
    (1+k)\mu-w=k\sqrt{\mu^2+y^2+2\varepsilon^2}.
\]
For $k>0$, squaring yields $c\mu^2-2(1+k)w\mu+w^2-k^2(y^2+2\varepsilon^2)=0$. The plus root is
\[
    M^\varepsilon(y,\nu)
    =\frac{(1+k)w+k\sqrt{w^2+c(y^2+2\varepsilon^2)}}{c}.
\]
It satisfies the unsquared equation: writing $R\mathrel{:=}\sqrt{w^2+c(y^2+2\varepsilon^2)}$, its left side is $k[kw+(1+k)R]/c\geq0$, since $R\geq|w|$. Strict monotonicity of $T_y^\varepsilon$ then establishes uniqueness. For $k=0$, the same formula gives $M^\varepsilon(y,\nu)=\nu$ directly. Substitution into the scalar residual gives
\[
    P^\varepsilon(y,\nu)
    =\frac{\sqrt{w^2+c(y^2+2\varepsilon^2)}-\nu+(1+k)y}{c}.
\]
In particular, $M^\varepsilon=\nu+kP^\varepsilon$, and substitution verifies that $P^\varepsilon$ is the FB residual evaluated at $M^\varepsilon$.

The lower slope bound for $T_y^\varepsilon$ and the estimate $|T_y^\varepsilon(\mu)-T_{\widetilde y}^\varepsilon(\mu)|\leq2k|y-\widetilde y|$ imply
\[
    |M^\varepsilon(y,\nu)-M^\varepsilon(\widetilde y,\widetilde\nu)|
    \leq|\nu-\widetilde\nu|+2k|y-\widetilde y|.
\]
For the residual map, regard $R$ as the Euclidean norm of $(\nu+ky,\sqrt c\,y,\sqrt{2c}\,\varepsilon)$. The reverse triangle inequality gives $|R-\widetilde R|\leq|\nu-\widetilde\nu|+(k+\sqrt c)|y-\widetilde y|$ and hence
\[
\begin{aligned}
    |P^\varepsilon(y,\nu)-P^\varepsilon(\widetilde y,\widetilde\nu)|
    \leq\frac{2}{c}|\nu-\widetilde\nu|
    +\left(1+\frac1{\sqrt c}\right)|y-\widetilde y|.
\end{aligned}
\]
These bounds hold globally and are independent of $\varepsilon$. Applying them componentwise proves the corresponding vector-map assertions.

\subsection{Solutions and the Integral Feedback Identity}
\label{app-real-solutions}

Since $f,g\in C^2$, the functions $g$, $J$, and $\nabla f$ are locally Lipschitz. The bounds above show that the right side of \eqref{eq:fb-pi-realization} is locally Lipschitz in $(x,\nu)$. Therefore every initial state has a unique maximal classical solution.

On each compact time interval contained in the maximal interval, the solution $(x,\nu)$ is continuously differentiable. The recovered functions $\mu=M^\varepsilon(g(x),\nu)$ and $\Phi^\varepsilon=P^\varepsilon(g(x),\nu)$ are locally absolutely continuous, by composition with the Lipschitz maps established above. Integrating $\dot\nu=K_iP^\varepsilon$ and using $\nu=\mu-K_p\Phi^\varepsilon$ gives \eqref{eq:fb-pi-integral-identity}. Consequently, the implicit feedback equality holds almost everywhere as an equality of derivatives of absolutely continuous functions, even when the ideal residual encounters a singular pair. A classical Jacobian formula for the residual is needed only to express this equality through partial derivatives.

Conversely, let a locally absolutely continuous physical trajectory satisfy the primal equation almost everywhere and \eqref{eq:fb-pi-integral-identity}. Defining $\nu=\mu-K_p\Phi^\varepsilon$ gives $\dot\nu=K_i\Phi^\varepsilon$ almost everywhere. Global inversion then gives \eqref{eq:fb-pi-realization}. Its continuous right side makes $(x,\nu)$ a classical solution, so uniqueness also applies to this integral interpretation of complete PI feedback.

If $g(x)=Cx-d$ and $\nabla f$ is globally Lipschitz, the realized vector field is globally Lipschitz: $J=C$ is constant and both inverse maps are globally Lipschitz. For each fixed $\varepsilon$, it therefore has at most linear growth in $(x,\nu)$, and all maximal forward solutions are defined for every $t\geq0$.

\subsection{Equilibrium Correspondence}
\label{app-real-equilibria}

At an equilibrium, $K_i\succ0$ gives $P^\varepsilon=\Phi^\varepsilon=0$, and the primal equation gives $\nabla_xL=0$. To identify the scalar residual zero set, write $b=-g_j(x)$. The equality $\sqrt{\mu_j^2+b^2+2\varepsilon^2}=\mu_j+b$ implies $\mu_j+b\geq0$ and $\mu_jb=\varepsilon^2$.

For $\varepsilon=0$, these conditions imply $\mu_j\geq0$, $b\geq0$, and $\mu_jb=0$: the product is zero, and the nonzero entry, if any, must be nonnegative. For $\varepsilon>0$, the entries have the same sign and a positive sum, so both are strictly positive. Conversely, these sign and product conditions imply the residual equation. Thus the realized equilibria correspond exactly to the ideal or perturbed KKT equations stated in Proposition~\ref{prop:exact-realization}. At every equilibrium, $\nu=\mu$, and none of these equations involves the feedback gains.

\begin{figure*}[t]
\centering
\includegraphics[width=0.8\textwidth]{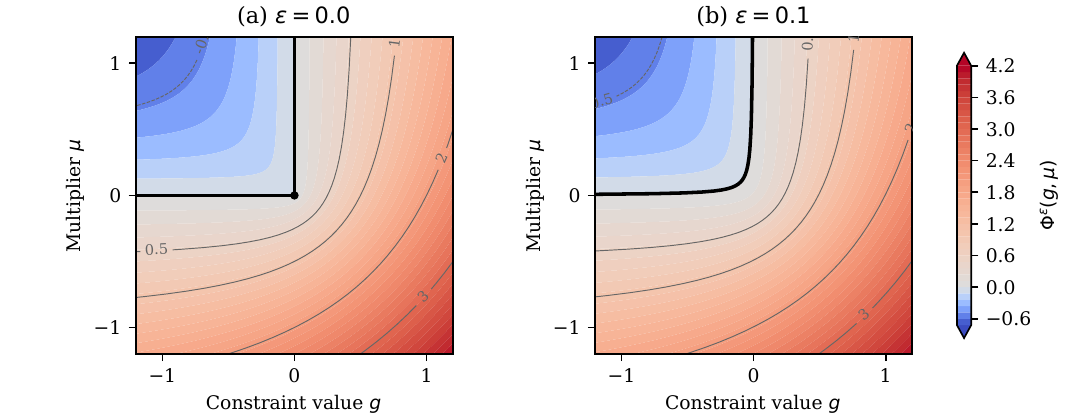}
\caption{Ideal and smoothed FB residuals. Black curves indicate their zero sets; panel (b) uses $\varepsilon=0.1$.}
\label{fig:app-fb-residual-contours}
\end{figure*}

Figure~\ref{fig:app-fb-residual-contours} displays the scalar residual
\(\Phi^\varepsilon(y,\mu)\mathrel{:=}\sqrt{\mu^2+y^2+2\varepsilon^2}-\mu+y\).
The ideal zero set consists of the rays $y=0$, $\mu\geq0$ and $\mu=0$, $y\leq0$, meeting at the singular pair $(y,\mu)=(0,0)$. Smoothing replaces these rays by the branch $y<0$, $\mu>0$, $-\mu y=\varepsilon^2$. The plotted value $\varepsilon=0.1$ makes this displacement visible; the main-text smoothed experiments use $10^{-6}$ and the Rosenbrock and Griewank tests in Appendix~\ref{app:exp-original-benchmarks} use $10^{-5}$.

\subsection{Cone Invariance and Zero Inactive Channels}
\label{app-real-invariance}

For every fixed $y$,
\[
    T_y^\varepsilon(0)=-k\bigl(\sqrt{y^2+2\varepsilon^2}+y\bigr)\leq0.
\]
Strict monotonicity of $T_y^\varepsilon$ therefore implies $M^\varepsilon(y,\nu)\geq0$ whenever $\nu\geq0$. On a boundary face $\nu_j=0$, if $k_{p,j}>0$, the identity $\mu_j=k_{p,j}P_j^\varepsilon$ gives $\dot\nu_j=k_{i,j}\mu_j/k_{p,j}\geq0$. If $k_{p,j}=0$, then $\mu_j=0$ and $\dot\nu_j=k_{i,j}(\sqrt{g_j(x)^2+2\varepsilon^2}+g_j(x))\geq0$. These boundary inequalities, together with local Lipschitz continuity, are the tangent-cone condition for forward invariance of $\mathbb R^n\times\mathbb R_+^m$. They prove invariance over the solution's interval of existence and nonnegativity of the recovered multiplier on this set.

For an ideal primal-dual feasible pair, $\mu_j\geq0$ and $b_j\mathrel{:=}-g_j(x)\geq0$ imply $\sqrt{\mu_j^2+b_j^2}\leq\mu_j+b_j$. Hence $\Phi_j^0\leq0$ and $\nu_j=\mu_j-k_{p,j}\Phi_j^0\geq0$. Thus primal-dual feasibility initially places the trajectory in the invariant integral-state cone.

The zero-channel assertion concerns the ideal system. Suppose $\varepsilon=0$, $\mu_j(t_0)=0$, and $g_j(x(t))<0$ on an interval starting at $t_0$. At its initial time, $\Phi_j^0=0$ and $\nu_j(t_0)=0$. For any $y<0$, direct substitution gives $M_j^0(y,0)=P_j^0(y,0)=0$. Treating the continuous function $g_j(x(t))$ along the solution as a prescribed input, the scalar equation $\dot\nu_j=k_{i,j}P_j^0(g_j(x(t)),\nu_j)$ is Lipschitz in $\nu_j$ and has the zero solution. Uniqueness implies $\nu_j(t)=\mu_j(t)=0$ throughout the interval.

If a constraint remains strictly satisfied along an original ideal trajectory, append its zero integral state and zero multiplier to that trajectory. The enlarged equations are satisfied because the new residual and its contribution to the primal equation vanish. Uniqueness of the enlarged realization proves that adding this channel leaves the original trajectory unchanged. For $\varepsilon>0$, the residual at a zero multiplier is strictly positive, so the zero-channel property is specific to the ideal system.

\section{Local Stability and Prescribed-Rate Proofs}
\label{app-local-proofs}

This appendix derives the linearization in Section~\ref{subsec:local-structure}, proves Theorems~\ref{thm:exact-stability} and~\ref{thm:prescribed-rate}, and details the active and inactive AL-PI comparison. The local assumptions and notation are those of Section~\ref{sec:local-stability}.

\subsection{Linearization and a Sufficient Proportional Gain}
\label{app-local-linearization}

At an active strictly complementary pair, the derivatives of
$\Phi_j^0(x,\mu)$ with respect to $\mu_j$ and $g_j$ are $0$ and $1$.
At an inactive pair they are $-1$ and $0$.
The residual is continuously differentiable in a neighborhood of the
reference point because none of these complementarity pairs is singular.
Linearizing the stationarity equation and complete PI feedback therefore gives
\[
\begin{aligned}
\delta\dot x&=-H\delta x-J_A^\top\delta\mu_A-J_I^\top\delta\mu_I,\\
\delta\dot\mu_A&=K_{p,A}J_A\delta\dot x+K_{i,A}J_A\delta x,\\
(I+K_{p,I})\delta\dot\mu_I&=-K_{i,I}\delta\mu_I.
\end{aligned}
\]
Terms differentiating residual derivative factors multiply zero equilibrium
velocities and thus do not contribute at first order.
Substitution of $\eta=\delta\mu_A-K_{p,A}J_A\delta x$ yields
\eqref{eq:local-linearization}.
Its block triangular structure shows that the inactive eigenvalues are
$-d_j=-k_{i,j}/(1+k_{p,j})$, regardless of their forcing into the primal equation.

For the active subsystem, let
$P\mathrel{:=}\operatorname{diag}(I,K_{i,A}^{-1})$ and denote its matrix by $M$.
For a complex eigenvector $(v,w)\ne0$ with eigenvalue $s$,
\[
\operatorname{Re}s\bigl(\|v\|^2+w^*K_{i,A}^{-1}w\bigr)
=-v^*H_ev.
\]
If $H_e\succ0$, this identity excludes $\operatorname{Re}s>0$.
Equality forces $v=0$, and the primal equation then implies
$J_A^\top w=0$. Full row rank of $J_A$ forces $w=0$, a contradiction.
Thus $M$ is Hurwitz for every $K_{i,A}\succ0$.
The inactive block is also Hurwitz, so the full equilibrium is locally
exponentially stable by the continuously differentiable linearization principle.

To obtain an explicit sufficient gain, choose orthonormal tangent and normal
bases $U,N$ and write
\[
\begin{aligned}
T_0&\mathrel{:=}U^\top HU,\; B_0\mathrel{:=}U^\top HN,\\
C_0&\mathrel{:=}N^\top HN,\; R_0\mathrel{:=}J_AN.
\end{aligned}
\]
Here $T_0\succ0$ and $R_0$ is nonsingular whenever both subspaces are nonempty.
Orthogonal congruence and the Schur complement give
\[
\begin{gathered}
H_e\succ0\ \Longleftrightarrow\\
C_0+R_0^\top K_{p,A}R_0-B_0^\top T_0^{-1}B_0\succ0.
\end{gathered}
\]
Consequently, it suffices to choose $K_{p,A}\succeq\kappa I$ with
\[
\kappa>
\max\left\{0,
\frac{\|B_0\|_2^2/\lambda_{\min}(T_0)-\lambda_{\min}(C_0)}
{\sigma_{\min}(R_0)^2}\right\}.
\]
If there are no active constraints, the tangent condition is $H\succ0$ and
no active gain is needed. If every primal direction is normal, it suffices that
$\kappa>\max\{0,-\lambda_{\min}(H)/\sigma_{\min}(J_A)^2\}$.
These sufficient bounds establish the existence of stabilizing proportional gains.

\subsection{Proof of Theorem~\ref{thm:exact-stability}}
\label{app-local-exact-stability}

Write $\delta x=n_Ay+Qw$. The active equations are
\[
\begin{aligned}
\dot y&=-ey-b^\top w-\sigma\eta,\\
\dot w&=-by-Tw,\; \dot\eta=k_i\sigma y.
\end{aligned}
\]
Let $d\mathrel{:=}\dim T$ and interpret determinants of empty matrices as one.
Eliminating $w$ and $\eta$ gives the characteristic polynomial
\begin{equation}
\begin{aligned}
p_q(s)&=\det(sI+T)(s^2+es+q)\\
&\phantom{=}-s b^\top\operatorname{adj}(sI+T)b,\\
p_q(s)&=s\det(sI+T)D(s,q),\\
D(s,q)&\mathrel{:=}s+e+q/s-b^\top(sI+T)^{-1}b.
\end{aligned}
\label{app-local-characteristic}
\end{equation}
The second identity is used away from its displayed poles and extends after
cancellation. In particular, $p_q(0)=q\det T>0$, so zero is never an eigenvalue.

Diagonalize $T$ orthogonally, writing its positive eigenvalues as $\tau_\ell$
and the coordinates of $b$ as $b_\ell$.
At $s=i\omega$, $\omega>0$, the two parts of $D=0$ are
\[
\begin{aligned}
e&=\sum_{\ell=1}^d\frac{\tau_\ell b_\ell^2}{\tau_\ell^2+\omega^2},\\
q&=\omega^2\left(1+\sum_{\ell=1}^d
\frac{b_\ell^2}{\tau_\ell^2+\omega^2}\right).
\end{aligned}
\]
The function
$h(t)\mathrel{:=}t[1+b^\top(T^2+tI)^{-1}b]$ satisfies
\[
\begin{gathered}
h(0)=0,\; \lim_{t\to\infty}h(t)=\infty,\\
h'(t)=1+b^\top T^2(T^2+tI)^{-2}b>0.
\end{gathered}
\]
Thus $t_q$ exists uniquely, and $e=\beta(q)$ is the only possible
imaginary-axis crossing as $e$ varies at fixed $q>0$.
At that crossing, direct differentiation gives
\[
\operatorname{Re}D_s(i\sqrt{t_q},q)
=2\bigl[1+b^\top T^2(T^2+t_qI)^{-2}b\bigr]>0.
\]
The conjugate pair is simple. Its continuation satisfies
\[
\begin{aligned}
\frac{\mathrm ds}{\mathrm de}&=-\frac{1}{D_s(s,q)},\\
\frac{\mathrm d\operatorname{Re}s}{\mathrm de}
&=-\frac{\operatorname{Re}D_s(s,q)}{|D_s(s,q)|^2}<0
\end{aligned}
\]
at the crossing.

For sufficiently large $e$, the Schur complement gives
$H_e\succ0$, hence stability by the preceding argument.
The eigenvalues vary continuously as a multiset, and their number in the
open right half-plane can change on a bounded parameter interval only at an
imaginary-axis crossing. There is exactly one such value of $e$.
Therefore all $e>\beta(q)$ are stable. Decreasing $e$ through $\beta(q)$
moves exactly two roots into the right half-plane, and no further crossing
can restore stability. At equality the unique imaginary pair is simple,
and all other roots remain strictly in the left half-plane.

Implicit differentiation of $h(t_q)=q$ gives
\[
\beta'(q)=
-\frac{b^\top T(T^2+t_qI)^{-2}b}
{1+b^\top T^2(T^2+t_qI)^{-2}b}.
\]
For $b\ne0$ this derivative is strictly negative. The limits follow from
$t_q\downarrow0$ as $q\downarrow0$ and $t_q\to\infty$ as $q\to\infty$.
Moreover, $\beta(q)>0$ for every finite $q$, so $e\le0$ is unstable when
$b\ne0$. If $b=0$, the polynomial factors as
$\det(sI+T)(s^2+es+q)$: $e>0$ is necessary and sufficient for stability,
$e<0$ is unstable, and $e=0$ has a purely imaginary normal pair.
The last case does not imply nonlinear instability merely from linearization.

The inactive eigenvalues remain strictly negative throughout.
A Hurwitz Jacobian is also necessary for local exponential stability of a
continuously differentiable flow: differentiating its local exponential
estimate with respect to the initial state at the equilibrium gives the
same estimate for the linearized flow. A right-half-plane eigenvalue implies
nonlinear instability. These facts establish the nonlinear statements of
Theorem~\ref{thm:exact-stability} and distinguish them from the boundary case.
Finally, $H_e\succ0$ is equivalent to $e>b^\top T^{-1}b$ by its Schur complement.
This proves the stated strict enlargement beyond the positive-definite region.
For the example $T=1$, $b=2$, $e=1$, the characteristic polynomial is
$s^3+2s^2+(q-3)s+q$. The cubic Hurwitz conditions reduce to $q>6$,
and at $q=6$ it factors as $(s+2)(s^2+3)$.

\subsection{High-Gain Root Structure}
\label{app-rate-high-gain}

We first establish the root expansion needed in the prescribed-rate proof.
All parameters except $q$ are fixed in this subsection.
On compact sets, $p_q(s)/q\to\det(sI+T)$.
The argument principle on contours avoiding $-\operatorname{spec}(T)$
therefore locates exactly $d$ bounded roots, counted with multiplicity,
converging to $-\operatorname{spec}(T)$.
Scaling $s=\sqrt q\,z$ in \eqref{app-local-characteristic} gives two remaining
roots near $z=\pm i$. These simple scaled roots have expansions
\begin{equation}
s_\pm(q)=\pm i\sqrt q-e/2+O(q^{-1/2}).
\label{app-rate-fast-roots}
\end{equation}
For example, substituting $s=\pm i\sqrt q+c+O(q^{-1/2})$ into $D=0$
gives $2c+e=0$.

A contour calculation supplies a real-part expansion without assuming that
the eigenvalues of $T$ are simple. Let $\Gamma$ be a positively oriented
contour enclosing all bounded roots and $-\operatorname{spec}(T)$ but no
escaping roots, and define
$f_T(s)\mathrel{:=}b^\top(sI+T)^{-1}b$ and $F(s)\mathrel{:=}s^2+es-sf_T(s)$.
Uniformly on $\Gamma$,
\[
\begin{aligned}
\frac{p_q'(s)}{p_q(s)}
&=\frac{\mathrm d}{\mathrm ds}\log\det(sI+T)\\
&\phantom{=}+q^{-1}F'(s)+O(q^{-2}).
\end{aligned}
\]
Multiplying by $s$ and integrating gives the sum of the bounded roots:
\[
\begin{aligned}
\sum_{\mathrm{bounded}}s
&=-\operatorname{tr}T
+\frac{1}{2\pi i q}\oint_\Gamma sf_T(s)\,\mathrm ds\\
&\phantom{=}+O(q^{-2})\\
&=-\operatorname{tr}T-\frac{b^\top Tb}{q}+O(q^{-2}).
\end{aligned}
\]
Here integration by parts removes $F'$, and the residues of $sf_T(s)$ at
$-\tau_\ell$ sum to $-\sum_\ell\tau_\ell b_\ell^2$.
The total sum of all roots is $-\operatorname{tr}T-e$.
The escaping roots form a conjugate pair for sufficiently large $q$, hence
\begin{equation}
\operatorname{Re}s_\pm(q)
=-e/2+\frac{c_b}{2q}+O(q^{-2}),\; c_b\mathrel{:=}b^\top Tb.
\label{app-rate-real-expansion}
\end{equation}
In particular, when $b\ne0$ and $0<r<\lambda_{\min}(T)$, all roots lie
strictly to the left of $-r$ for sufficiently large $q$ if $e>2r$.
If $e\le2r$, exactly two roots instead lie strictly to its right for
sufficiently large $q$. The equality case uses $c_b>0$ in
\eqref{app-rate-real-expansion}.

\subsection{Proof of the Prescribed-Rate Condition}
\label{app-rate-prescribed}

Assume first that $b\ne0$, and fix $0<r<\lambda_{\min}(T)$.
Neither $s$ nor $\det(sI+T)$ vanishes on $\operatorname{Re}s=-r$.
At a nonreal root $s=-r+i\omega$, separating $D=0$ into real and imaginary
parts yields
\[
q=h_r(\omega^2),\; e-2r=\beta_r(\omega^2).
\]
Since
$\beta_r'(t)=-b^\top T(L_r^2+tI)^{-2}b<0$ and
$\beta_r(t)\to0$, a nonreal crossing exists exactly when
$2r<e<e_r$. It occurs at the unique gain
$q_H\mathrel{:=}h_r(t_{r,e})$. The only possible real crossing is $s=-r$, at
\begin{equation}
q_R\mathrel{:=}r(e-r-B_r).
\label{app-rate-real-boundary}
\end{equation}
When the nonreal crossing exists, substitution of
$e-2r=\beta_r(t_{r,e})$ gives
\[
q_H-q_R=t_{r,e}\left[
1+b^\top T L_r^{-1}(L_r^2+t_{r,e}I)^{-1}b\right]>0.
\]

At a nonreal crossing, implicit differentiation gives
$\mathrm ds/\mathrm dq=-1/(sD_s)$ and
\[
\operatorname{Re}(sD_s)
=2\omega^2 b^\top T(L_r^2+\omega^2I)^{-2}b>0.
\]
Thus the pair is simple and crosses to the left as $q$ increases.
At the real crossing,
\[
D_s(-r,q_R)=\frac{e_r-e}{r},\;
\frac{\mathrm ds}{\mathrm dq}=\frac{1}{e_r-e}
\]
when $e\ne e_r$. This root crosses to the right if $e<e_r$ and to the
left if $e>e_r$.

At the meeting point $e=e_r$, the identities
$\beta_r(0)=B_r+r b^\top L_r^{-2}b$ and
\eqref{app-rate-real-boundary} give
\[
q_R=r^2(1+b^\top L_r^{-2}b)>0,\; D_s(-r,q_R)=0,
\]
whereas
\[
D_{ss}(-r,q_R)
=-2q_R/r^3-2b^\top L_r^{-3}b<0.
\]
Hence the root has multiplicity exactly two; an exceptional triple root is
excluded. To determine the nearby real splitting, set $G(s,q)\mathrel{:=}sD(s,q)$.
Then $G_q=1$ and $G_{ss}(-r,q_R)=-rD_{ss}(-r,q_R)>0$.
For $q<q_R$ sufficiently close to $q_R$, the two nearby roots are
\[
s=-r\pm\sqrt{\frac{2(q_R-q)}{G_{ss}(-r,q_R)}}+O(q_R-q).
\]
They are real and lie on opposite sides of the target line.

We can now count roots by decreasing $q$ from a sufficiently large value.
If $2r<e<e_r$, all roots initially lie to the left of the target line.
The unique nonreal crossing at $q_H$ introduces two right-hand roots.
The only additional crossing is $q_R<q_H$; if it is positive, it removes
only one right-hand root as $q$ decreases. Thus the strict target rate holds
exactly for $q>q_H$.
If $e>e_r$, there is no nonreal crossing and $q_R>0$.
The unique real crossing introduces one right-hand root as $q$ decreases,
so the strict target rate holds exactly for $q>q_R$.
If $e=e_r$, the double-root calculation gives one right-hand root immediately
below $q_R$, and there is no further positive crossing. The same condition
$q>q_R$ therefore holds. At every limiting threshold, all other roots lie
strictly to the left of the target line.

Finally, if $e\le2r$, there are already two right-hand roots at large $q$
and no nonreal crossing is available. Any positive real crossing has
$e<e_r$ and removes only one right-hand root as $q$ decreases.
At least one root remains strictly to the right for every finite $q>0$,
including any real-crossing value itself. No gain attains the strict target.
This proves the first part of Theorem~\ref{thm:prescribed-rate} for $b\ne0$.
If $b=0$, the normal factor shifted by $s=z-r$ is
\[
z^2+(e-2r)z+q-r(e-r).
\]
Its Hurwitz conditions are precisely $e>2r$ and $q>r(e-r)$;
the tangent roots already satisfy the target inequality.

\subsection{Rate Ceiling, Gain Cost, and Euler Step Size}
\label{app-rate-costs}

Suppose $b\ne0$ and $0<e<2\lambda_{\min}(T)$.
Apply the preceding root count with $r=e/2$.
It leaves a root strictly to the right of $-e/2$ for every finite $q>0$,
so $\gamma(q)<e/2$.
For sufficiently large $q$, the tangent roots have a strictly larger decay
margin than the escaping pair. Equation~\eqref{app-rate-real-expansion}
therefore gives
\[
\begin{gathered}
\gamma(q)=e/2-c_b/(2q)+O(q^{-2}),\\
\sup_{q>0}\gamma(q)=e/2.
\end{gathered}
\]
Let $r=e/2-\delta$ with $\delta\downarrow0$.
The nonreal branch applies for all sufficiently small $\delta$, and its
frequency parameter tends to infinity. Uniformly for $r$ near $e/2$,
\[
2\delta=\beta_r(t)=c_b/t+O(t^{-2}).
\]
It follows that $t=c_b/(2\delta)+O(1)$ and
$h_r(t)=t+O(1)$, proving the asserted gain cost.
When $b=0$, the explicit roots of $s^2+es+q$ show that the normal spectral
margin equals $e/2$ for $q\ge e^2/4$.
This is a statement about spectral margin: at $q=e^2/4$ the double normal
root can be defective, so an exponential estimate at exactly that limiting
rate does not follow solely from its value.

For a Hurwitz matrix, its explicit Euler linearization is Schur exactly when
\[
0<\alpha<\min_{s\in\operatorname{spec}(M_q)}
\frac{-2\operatorname{Re}s}{|s|^2}.
\]
The bounded roots have positive limiting step bounds, since their limits
are the negative eigenvalues of $T$.
For the fast pair, \eqref{app-rate-fast-roots} and
\eqref{app-rate-real-expansion} give
$-2\operatorname{Re}s_\pm=e+O(q^{-1})$ and $|s_\pm|^2=q+O(1)$.
They therefore determine the minimum for all sufficiently large $q$, yielding
$\alpha_{\max}^{\mathcal A}=e/q[1+o(1)]$.
Along the stated near-threshold path,
$q/q_{\min}(e/2-\delta,e)\to1$ and hence
$\alpha_{\max}^{\mathcal A}\sim2e\delta/c_b$.
Each inactive eigenvalue $-d_j$ contributes the bound
$2/d_j=2(1+k_{p,j})/k_{i,j}$.
Taking the minimum proves the full-system step bound in the main text.
The physical and integral linear coordinates are related by a nonsingular
linear map, so they have the same Euler eigenvalue criterion.
Their nonlinear Euler updates need not coincide.

To transfer a strict spectral rate to the nonlinear system, include the
inactive block and assume $r<\min\{\gamma(q),\min_j d_j\}$.
Choose $r_0$ strictly between $r$ and this minimum.
For the full Jacobian $A$, there is a constant $C$ such that
$\|e^{At}\|_2\le C e^{-r_0t}$.
The continuously differentiable vector field has remainder
$R(z)=o(\|z\|)$ in equilibrium-centered coordinates.
Variation of constants and Gronwall's inequality, with
$\|R(z)\|\le\epsilon_0\|z\|$ in a sufficiently small neighborhood and
$C\epsilon_0<r_0-r$, give a local bound with rate at least $r$.
A smaller initial neighborhood keeps the trajectory inside that neighborhood.
The estimate constant and the admissible neighborhood may depend on the gains.

\subsection{The Active Linearization of Augmented-Lagrangian PI Feedback}
\label{app-local-al}

For $\rho>0$, define $p_\rho\mathrel{:=}[\mu+\rho g(x)]_+$ and
$q_\rho\mathrel{:=}(p_\rho-\mu)/\rho$. The AL-PI implementation with
proportional feedback of the constraint velocity is
\begin{equation}
\begin{aligned}
\dot x&=-\nabla f(x)-J(x)^\top p_\rho,\\
\dot\mu&=k_iq_\rho+k_pJ(x)\dot x,
\end{aligned}
\label{eq:al-pi-feedback}
\end{equation}
where $k_p\geq0$ and $k_i>0$; see
\cite[Eq.~(9)]{CeroneFossonPirreraRegruto2024CDC}.
At a strictly complementary point, its branches are fixed locally.
Suppressing the subscript $\rho$ in the variations gives
\[
\begin{aligned}
\delta p_A&=\delta\mu_A+\rho J_A\delta x,\; \delta p_I=0,\\
\delta q_A&=J_A\delta x,\; \delta q_I=-\delta\mu_I/\rho.
\end{aligned}
\]
For affine constraints, the linearization is
\[
\begin{aligned}
\delta\dot x&=-(H+\rho J_A^\top J_A)\delta x-J_A^\top\delta\mu_A,\\
\delta\dot\mu_A&=k_iJ_A\delta x+k_pJ_A\delta\dot x,\\
\delta\dot\mu_I&=-k_i\delta\mu_I/\rho+k_pJ_I\delta\dot x.
\end{aligned}
\]
For the nonlinear-constraint extension, differentiating $J(x)^\top p_\rho$
adds the constraint Hessian terms weighted by the reference multipliers.
These are exactly the terms already included in the Lagrangian Hessian $H$,
so the same displayed linearization applies at the reference point.
The change $\eta_{\mathrm{AL}}\mathrel{:=}\delta\mu_A-k_pJ_A\delta x$ gives the active
matrix in \eqref{eq:local-linearization}, with
$H_e=H+(\rho+k_p)J_A^\top J_A$ and $K_{i,A}=k_iI$.
Thus matching $k_p^{\mathrm{FB}}=\rho+k_p^{\mathrm{AL}}$ and $k_i$ matches
the active matrices in their respective integral coordinates.

The AL inactive subsystem is driven by the active subsystem but does not
feed back into it, whereas the FB linearization has the opposite triangular
orientation. Both orientations preserve the active eigenvalues.
If the complete residual derivative $k_p\dot q_\rho$ is used instead of
$k_pJ\dot x$, its inactive equation becomes
$(1+k_p/\rho)\delta\dot\mu_I=-(k_i/\rho)\delta\mu_I$.
Its decay rate is therefore $k_i/(\rho+k_p)$.
On an interval where that channel remains strictly inactive, an initial
zero multiplier is also preserved.

\section{Multiple Active Constraints and Modal Damping}
\label{app-modal-analysis}

This appendix proves the high-integral-gain statements in
Section~\ref{subsec:multiple-active} and gives a sufficient finite-gain
certificate. The reference ideal KKT point satisfies the regularity
assumptions of Section~\ref{subsec:local-structure}. The proportional gains and inactive integral
gains are fixed. The active integral gain is
$K_{i,A}=\kappa K_0$, where $\kappa>0$ and $K_0$ is a fixed positive
diagonal matrix. Matrix norms without subscripts are Euclidean operator
norms, and $\|\cdot\|_F$ denotes the Frobenius norm.

\subsection{Tangent--Normal Reduction}
\label{app-modal-reduction}

Let $m_A\geq1$ be the number of active constraints. Choose matrices $U$ and
$N$ with orthonormal columns spanning $\ker J_A$ and its orthogonal
complement, respectively, so that $[U\ N]$ is orthogonal. Set
$R\mathrel{:=}J_AN$. The full row rank of $J_A$ implies that the square matrix
$R$ is nonsingular. In addition to the matrices used in the main text,
define $B\mathrel{:=}U^\top HN$. Then
\[
\begin{aligned}
T&\mathrel{:=}U^\top HU,\; Q_0\mathrel{:=}R^\top K_0R,\\
E&\mathrel{:=}N^\top HN+R^\top K_{p,A}R.
\end{aligned}
\]
Here $T\succ0$, $E=E^\top$, and $Q_0\succ0$. If $m_A=n$, all
blocks involving $U$, $T$, and $B$ are omitted.

For the active subsystem of \eqref{eq:local-linearization}, introduce
$u\mathrel{:=}U^\top\delta x$, $v\mathrel{:=}N^\top\delta x$, and $p\mathrel{:=}R^\top\eta$.
Using $J_AU=0$ gives
\begin{equation}
\frac{\mathrm d}{\mathrm dt}
\begin{bmatrix}u\\v\\p\end{bmatrix}
=
\begin{bmatrix}
-T&-B&0\\
-B^\top&-E&-I\\
0&\kappa Q_0&0
\end{bmatrix}
\begin{bmatrix}u\\v\\p\end{bmatrix}.
\label{app-modal-active-matrix}
\end{equation}
The matrix $Q_0=R^\top K_0R$ is determined jointly by the constraint
geometry and the diagonal entries of $K_0$.

Put $\omega\mathrel{:=}\sqrt\kappa$, $S\mathrel{:=}Q_0^{1/2}$, and
$z\mathrel{:=}S^{-1}p/\omega$. This coordinate is related to the weighted integral coordinate
$\zeta\mathrel{:=}\kappa^{-1/2}K_0^{-1/2}\eta$ by
\[
z=O\zeta,\; O\mathrel{:=}Q_0^{-1/2}R^\top K_0^{1/2},\; OO^\top=I.
\]
Thus the two integral coordinates differ by an orthogonal transformation.
The scaled state $y\mathrel{:=}(u,v,z)$ satisfies
$\dot y=(\omega A_0+A_1)y$, where
\[
\begin{aligned}
A_0&\mathrel{:=}\begin{bmatrix}0&0&0\\0&0&-S\\0&S&0\end{bmatrix},\\
A_1&\mathrel{:=}\begin{bmatrix}-T&-B&0\\-B^\top&-E&0\\0&0&0\end{bmatrix}.
\end{aligned}
\]
The fast operator is skew-symmetric, $A_0^\top=-A_0$, whereas
$A_1^\top=A_1$. This separates rapid conservative motion from the
curvature and damping terms, including an indefinite normal block.

\subsection{Spectral Expansion and the Uniform Margin Criterion}
\label{app-modal-spectrum}

Let $\lambda_1,\ldots,\lambda_h$ be the distinct eigenvalues of $Q_0$,
and let $\Pi_j$ be the corresponding orthogonal projectors. Define
$E_j\mathrel{:=}\Pi_jE\Pi_j|_{\operatorname{range}\Pi_j}$. Counting algebraic
multiplicity, the active eigenvalues have the expansions
\begin{equation}
\begin{aligned}
s_{\ell}^{0}(\kappa)
&=-\lambda_\ell(T)+O(\kappa^{-1/2}),\\
s_{j,\ell}^{\pm}(\kappa)
&=\pm\mathrm i\sqrt{\kappa\lambda_j}
-\frac12\lambda_\ell(E_j)+O(\kappa^{-1/2}).
\end{aligned}
\label{app-modal-expansions}
\end{equation}
There are $n-m_A$ slow eigenvalues in the first line and
$2\operatorname{rank}\Pi_j$ fast eigenvalues associated with each $j$.
The constants in the remainders may depend on the fixed problem data
and on $K_0$. The expansions hold as multiset statements; no analytic
labeling of individual eigenvalues through repeated roots is assumed.

\paragraph{Derivation of the expansion.}
The zero eigenspace of $A_0$, when present, consists of $(u,0,0)$,
and the compression of $A_1$ to that space is $-T$. For
$r\in\operatorname{range}\Pi_j$, the isometries
\[
V_j^\pm r\mathrel{:=}2^{-1/2}(0,r,\mp\mathrm i r)
\]
span the eigenspaces of $A_0$ at $\pm\mathrm i\sqrt{\lambda_j}$.
Their compressed first-order matrices are
$(V_j^\pm)^*A_1V_j^\pm=-E_j/2$, where $*$ denotes conjugate transpose.
In particular, all compressed first-order matrices are Hermitian,
even when the frequencies or their damping eigenvalues are repeated.

To justify the error order in this setting, put $\tau\mathrel{:=}\omega^{-1}$
and consider $A_0+\tau A_1$. For a distinct eigenvalue $\mathrm i\nu$
of the normal matrix $A_0$, a sufficiently small fixed contour enclosing
only that eigenvalue defines an analytic spectral projector for all
small $\tau$. An analytic change of basis identifies its range with
the unperturbed eigenspace. On that fixed space, the reduced matrix is
\[
\mathrm i\nu I+\tau C_\nu+\tau^2F_\nu(\tau),
\]
where $C_\nu$ is the compression of $A_1$ and $F_\nu$ is bounded.
After subtracting $\mathrm i\nu I$ and dividing by $\tau$, the
remaining matrix is $C_\nu+O(\tau)$. Since $C_\nu$ is Hermitian,
its resolvent outside its spectrum has norm equal to the inverse
distance to that spectrum. The Neumann-series resolvent test therefore
places every perturbed eigenvalue within $O(\tau)$ of
$\operatorname{spec}(C_\nu)$. Spectral projection counts around the
distinct eigenvalues of $C_\nu$ preserve their multiplicities.
Consequently, the reduced eigenvalues are
$\mathrm i\nu+\tau\lambda_\ell(C_\nu)+O(\tau^2)$, including repeated
values of $\lambda_\ell(C_\nu)$. Multiplication by $\omega$ gives
\eqref{app-modal-expansions}.

\paragraph{Consequences for stability.}
Define the active limiting margin
\[
\gamma_*\mathrel{:=}
\min\left\{\lambda_{\min}(T),
\frac12\min_j\lambda_{\min}(E_j)\right\},
\]
where the tangent term is omitted when there are no tangent variables.
The active spectral margin converges to $\gamma_*$. For the complete
linearization, the block-triangular form in
\eqref{eq:local-linearization} adds the fixed inactive eigenvalues
$-d_j$, with $d_j=k_{i,j}/(1+k_{p,j})>0$. Hence its limiting margin is
$\min\{\gamma_*,\min_{j\in\mathcal I}d_j\}$, with the second minimum
omitted if there are no inactive constraints.

It follows that a positive spectral margin, uniform for all sufficiently
large $\kappa$, exists if and only if every $E_j\succ0$. Sufficiency
follows by choosing any positive number strictly below the limiting
margin and applying \eqref{app-modal-expansions}. For necessity, if an
$E_j$ has a negative eigenvalue, a fast eigenvalue has a positive limiting
real part and is therefore unstable for every sufficiently large gain.
If an $E_j$ has a zero eigenvalue, the real part of a corresponding fast
eigenvalue tends to zero, which precludes a uniform strict margin.
When no compression has a negative eigenvalue but at least one has a
zero eigenvalue, higher-order terms determine eventual Hurwitz
stability; the leading expansion alone does not decide it.

Strict complementarity makes the ideal realized vector field
continuously differentiable in a neighborhood of the reference point.
Thus a Hurwitz linearization implies local exponential stability for
each fixed sufficiently large gain, and a positive-real-part eigenvalue
implies nonlinear instability. The admissible nonlinear neighborhood may depend on the gain;
the marginal cases require higher-order analysis.

\subsection{A Constructive Finite-Gain Certificate}
\label{app-modal-finite-gain}

Suppose every $E_j\succ0$, and define
\[
\begin{aligned}
\overline E&\mathrel{:=}\sum_{j=1}^h\Pi_jE\Pi_j,\\
\overline A&\mathrel{:=}\operatorname{diag}(-T,-\overline E/2,-\overline E/2).
\end{aligned}
\]
The rapid-motion average is
\begin{equation}
\overline A=
\lim_{b\to\infty}\frac1b\int_0^b
 e^{-tA_0}A_1e^{tA_0}\,\mathrm dt.
\label{app-modal-average}
\end{equation}
Indeed, in the complex orthogonal spectral decomposition of $A_0$,
conjugation multiplies the block between frequencies $\nu$ and $\nu'$
by $e^{\mathrm i(\nu'-\nu)t}$. Its average vanishes unless
$\nu=\nu'$. Retaining the equal-frequency blocks computed above gives
\eqref{app-modal-average}. In particular,
$\overline A\preceq-\gamma_*I$ with $\gamma_*>0$.

There is a real symmetric matrix $P_1$, with zero component commuting
with $A_0$, that solves the homological equation
\begin{equation}
P_1A_0-A_0P_1=-2(A_1-\overline A).
\label{app-modal-homological}
\end{equation}
For an explicit construction, let $\Pi_\nu$ denote the complex
orthogonal spectral projector of $A_0$ at $\mathrm i\nu$. Set
\[
\Pi_\nu P_1\Pi_{\nu'}\mathrel{:=}
\begin{cases}
-2\Pi_\nu A_1\Pi_{\nu'}/[\mathrm i(\nu'-\nu)],&\nu\ne\nu',\\
0,&\nu=\nu'.
\end{cases}
\]
The Hermitian and complex-conjugacy relations between these blocks
make $P_1$ real symmetric. Substitution verifies
\eqref{app-modal-homological} block by block.

Put $R_1\mathrel{:=}A_1P_1+P_1A_1$. A sufficient finite-gain condition is
\begin{equation}
\sqrt\kappa\geq
\max\left\{2\|P_1\|,\frac{\|R_1\|}{\gamma_*}\right\}.
\label{app-modal-certificate}
\end{equation}
To verify it, define $P_\kappa\mathrel{:=}I+P_1/\sqrt\kappa$. Then
$\tfrac12I\preceq P_\kappa\preceq\tfrac32I$, and direct expansion yields
\[
\begin{aligned}
&(\omega A_0+A_1)^\top P_\kappa+P_\kappa(\omega A_0+A_1)\\
&=2\overline A+\frac{R_1}{\omega}\preceq-\gamma_*I.
\end{aligned}
\]
Thus $V(y)\mathrel{:=}y^\top P_\kappa y$ satisfies
$\dot V\leq-\gamma_*\|y\|^2\leq-(2\gamma_*/3)V$, proving the scaled
active trajectory bound
\begin{equation}
\|y(t)\|\leq\sqrt3\,e^{-\gamma_*t/3}\|y(0)\|.
\label{app-modal-scaled-bound}
\end{equation}
Condition \eqref{app-modal-certificate} certifies the active block;
the fixed negative inactive eigenvalues then certify the full spectrum.

The certificate also displays the cost of nearby distinct frequencies.
Let $\Omega$ contain the distinct values $\pm\sqrt{\lambda_j}$,
together with zero when a tangent block is present, and define
$\delta\mathrel{:=}\min_{\nu\ne\nu',\,\nu,\nu'\in\Omega}|\nu-\nu'|>0$.
Orthogonality of the spectral blocks in Frobenius norm gives
\[
\begin{aligned}
\|P_1\|&\leq\|P_1\|_F\leq\frac{2\|A_1-\overline A\|_F}{\delta},\\
\|R_1\|&\leq2\|A_1\|\|P_1\|.
\end{aligned}
\]
Therefore, the explicit sufficient condition
\begin{equation}
\kappa\geq
\frac{16\|A_1-\overline A\|_F^2}{\delta^2}
\max\left\{1,\frac{\|A_1\|^2}{\gamma_*^2}\right\}
\label{app-modal-gap-certificate}
\end{equation}
implies \eqref{app-modal-certificate}. At fixed dimension and bounded
matrix norms and $\gamma_*^{-1}$, this is an $O(\delta^{-2})$
sufficient gain as distinct frequencies approach one another.

The norm in \eqref{app-modal-scaled-bound} involves
$z=Q_0^{-1/2}R^\top\eta/\sqrt\kappa$. Its equivalence constants with
an unscaled physical-state norm can depend on $\kappa$. For example,
with one normal variable, $E=e>0$, and $Q_0=1$, the initial condition
$v(0)=1$, $p(0)=0$ gives
\[
p(t)=\frac{\kappa}{\sqrt{\kappa-e^2/4}}
 e^{-et/2}\sin\bigl(\sqrt{\kappa-e^2/4}\,t\bigr)
\]
for $\kappa>e^2/4$. Its magnitude near time $\pi/(2\sqrt{\kappa-e^2/4})$
is of order $\sqrt\kappa$. Hence a uniform strict spectral margin
and the scaled estimate do not imply a gain-independent transient
prefactor in the physical multiplier coordinates.

\subsection{An Exact Two-Channel Example}
\label{app-modal-two-channel}

For $B=0$, differentiation of the normal equation in
\eqref{app-modal-active-matrix} gives
$\ddot v+E\dot v+\kappa Q_0v=0$. Take
\[
E=\begin{bmatrix}e_1&c\\c&e_2\end{bmatrix},\;
Q_0=\operatorname{diag}(q_1,q_2),
\]
with $e_1,e_2,q_1,q_2>0$ and $c^2>e_1e_2$. If $q_1=q_2$, an
orthogonal diagonalization of $E$ yields a scalar oscillator with
negative damping, so no positive $\kappa$ stabilizes the system.
If $q_1\ne q_2$, its exact Hurwitz condition is
\begin{equation}
\kappa>
\frac{(e_1+e_2)(e_1q_2+e_2q_1)(c^2-e_1e_2)}
{e_1e_2(q_1-q_2)^2}.
\label{app-modal-two-threshold}
\end{equation}
To check this, write $a\mathrel{:=}e_1+e_2$, $b_q\mathrel{:=}q_1+q_2$, and
$c_q\mathrel{:=}e_1q_2+e_2q_1$. The characteristic polynomial is
\[
s^4+as^3+(\det E+\kappa b_q)s^2+
\kappa c_qs+\kappa^2q_1q_2.
\]
Its third Hurwitz determinant is
$\Delta_3=\kappa[ac_q\det E+\kappa e_1e_2(q_1-q_2)^2]$.
Since $a$, $\kappa c_q$, and $\kappa^2q_1q_2$ are positive,
$\Delta_3>0$ also forces the second Hurwitz determinant and the
quadratic coefficient to be positive. The quartic Hurwitz criterion
therefore gives precisely \eqref{app-modal-two-threshold}.

This calculation is realizable with diagonal multiplier gains. Consider
\[
\begin{aligned}
\min_{y,v_1,v_2}\;&\frac12y^2+2v_1v_2-v_1-v_2,\\
\text{s.t.}\;&v_1\leq0,\;v_2\leq0.
\end{aligned}
\]
The unique global minimizer is the origin, with multipliers $(1,1)$.
The active gradients are independent, strict complementarity holds,
and the tangent curvature is $T=1$. Choosing $K_{p,A}=I$ and
$K_{i,A}=\kappa\operatorname{diag}(1,r)$, $r>0$, gives
$E=\left[\begin{smallmatrix}1&2\\2&1\end{smallmatrix}\right]$
and $Q_0=\operatorname{diag}(1,r)$. Equal gains never stabilize the
linearization. For $r\ne1$, local exponential stability holds exactly
when $\kappa>6(1+r)/(r-1)^2$. For example, $r=4$ gives threshold
$10/3$. As $r\to1$, this threshold is asymptotic to $12/(r-1)^2$,
while the smallest distinct-frequency gap is asymptotic to $|r-1|/2$.
The exact required gain is consequently $\Theta(\delta^{-2})$ along
this family, matching the order of \eqref{app-modal-gap-certificate}.
The proportional correction remains fixed; only the relative integral
gains separate or merge the normal frequencies.

\section{Nonlinear Dissipation and Global Convergence}
\label{app-global-proofs}

This appendix proves the dissipation statements in
Section~\ref{sec:global-convergence}, the smoothing-uniform estimate
\eqref{eq:smoothing-uniform-ges}, the global error bound
\eqref{eq:global-residual-bound},
and Theorem~\ref{thm:global-geometry}. Throughout, $g(x)=Cx-d$,
$f$ is $\ell$-strongly convex, and $\nabla f$ has global Lipschitz
constant $L_f$. All gains are fixed diagonal matrices with
$K_i\succ0$ and $K_p\succeq0$.

\subsection{The Shifted Scalar Inequality}
\label{app-global-scalar}

Let $\psi\mathrel{:=}\sqrt{\mu^2+g^2+2\varepsilon^2}-\mu+g$, and suppose that
$\bar\mu\geq0$, $\bar g\leq0$, and $-\bar\mu\bar g=\varepsilon^2$.
Squaring the defining identity of $\psi$ gives
\[
    \mu g=\mu\psi-g\psi+\tfrac12\psi^2-\varepsilon^2.
\]
Consequently, the difference between the left and right sides of
\eqref{eq:shifted-fb-inequality} is exactly
\[
    B\mathrel{:=}\bar g(\mu+\psi)+\bar\mu(g-\psi)+2\varepsilon^2.
\]
Write $s\mathrel{:=}-g$, $\bar s\mathrel{:=}-\bar g$, and
$r\mathrel{:=}\sqrt{\mu^2+s^2+2\varepsilon^2}$. Since
$\bar\mu\bar s=\varepsilon^2$, we have
\[
    \| (\bar\mu,\bar s,\sqrt{2}\varepsilon)\|
    =\bar\mu+\bar s.
\]
Using $\mu+\psi=r-s$ and $g-\psi=\mu-r$ therefore yields
\[
\begin{aligned}
    B
    &=\bar\mu\mu+\bar s s+2\varepsilon^2
      -(\bar\mu+\bar s)r\leq0,
\end{aligned}
\]
by the Cauchy--Schwarz inequality. This proves
\eqref{eq:shifted-fb-inequality}. The argument places no sign
restriction on the current $\mu$ or $g$ and uses no derivative of the
residual. It therefore also covers infeasible states, negative
multipliers, and the ideal singular pair $(\mu,g)=(0,0)$.

\subsection{Energy Dissipation}
\label{app-global-energy}

Fix $\varepsilon\geq0$ and an equilibrium $(\bar x,\bar\mu)$.
Set $u\mathrel{:=}x-\bar x$, $v\mathrel{:=}\nu-\bar\mu$,
$y\mathrel{:=}Cx-d$, $\bar y\mathrel{:=}C\bar x-d$,
$\psi\mathrel{:=}P^\varepsilon(y,\nu)$, and
$h\mathrel{:=}\nabla f(x)-\nabla f(\bar x)$. The equilibrium residual vanishes,
and the exact realization gives
\[
\begin{aligned}
    \mu-\bar\mu&=v+K_p\psi,\\
    \dot u&=-h-C^\top v-C^\top K_p\psi,\\
    \dot v&=K_i\psi.
\end{aligned}
\]
For $V\mathrel{:=}\|u\|^2/2+v^\top K_i^{-1}v/2$, differentiation yields
\[
    \dot V=-u^\top h-(Cu)^\top(\mu-\bar\mu)+v^\top\psi.
\]
Applying the preceding scalar identity componentwise gives the exact
expression
\[
\begin{aligned}
    \dot V={}&-u^\top h+(Cu)^\top\psi\\
    &-\psi^\top(K_p+I/2)\psi+\sum_{j=1}^m B_j,
\end{aligned}
\]
where
$B_j\mathrel{:=}\bar y_j(\mu_j+\psi_j)+\bar\mu_j(y_j-\psi_j)+2\varepsilon^2
\leq0$. Strong convexity gives $u^\top h\geq\ell\|u\|^2$ and hence
\eqref{eq:global-dissipation}.

Define
\[
    G\mathrel{:=}
    \begin{bmatrix}
        \ell I&-C^\top/2\\
        -C/2&K_p+I/2
    \end{bmatrix}.
\]
The Schur complement of its lower diagonal block is
$\ell I-C^\top\operatorname{diag}((4k_{p,j}+2)^{-1})C$.
Thus \eqref{eq:global-gain-condition} implies $G\succ0$.
With $\gamma\mathrel{:=}\lambda_{\min}(G)>0$, we obtain
\begin{equation}
    \dot V\leq-\gamma(\|u\|^2+\|\psi\|^2).
    \label{app-global-Vdot}
\end{equation}
Equivalently, completing the square gives
\[
\begin{aligned}
    \dot V\leq{}&-\ell\|u\|^2
    +u^\top C^\top\operatorname{diag}\!\left(\frac1{4k_{p,j}+2}\right)Cu\\
    &-\sum_{j=1}^m(k_{p,j}+1/2)
    \left(\psi_j-\frac{(Cu)_j}{2k_{p,j}+1}\right)^2.
\end{aligned}
\]
For $K_p=kI$, the sufficient scalar condition in the main text follows
from $C^\top C\preceq\|C\|^2I$.

\subsection{A Modified Energy with Integral-State Dissipation}
\label{app-global-modified-energy}

The following construction does not require any rank condition on $C$.
Let $Z\mathrel{:=}u^\top C^\top v$. Along the realized system,
\[
\begin{aligned}
    \dot Z={}&-h^\top C^\top v-\|C^\top v\|^2\\
    &-(C^\top v)^\top C^\top K_p\psi
    +u^\top C^\top K_i\psi.
\end{aligned}
\]
Since $\|h\|\leq L_f\|u\|$, Young's inequality gives
\[
\begin{aligned}
    |h^\top C^\top v|
    &\leq\tfrac14\|C^\top v\|^2+L_f^2\|u\|^2,\\
    |(C^\top v)^\top C^\top K_p\psi|
    &\leq\tfrac14\|C^\top v\|^2
      +\|C\|^2\|K_p\|^2\|\psi\|^2,\\
    |u^\top C^\top K_i\psi|
    &\leq\tfrac12\|C\|\|K_i\|(\|u\|^2+\|\psi\|^2).
\end{aligned}
\]
Thus, writing
$a_0\mathrel{:=}L_f^2+\|C\|\|K_i\|/2$ and
$b_0\mathrel{:=}\|C\|^2\|K_p\|^2+\|C\|\|K_i\|/2$, we have
\[
    \dot Z\leq-\tfrac12\|C^\top v\|^2
    +a_0\|u\|^2+b_0\|\psi\|^2.
\]
For explicit comparison constants, put
\[
\begin{aligned}
    v_-&\mathrel{:=}\tfrac12\min\{1,\lambda_{\min}(K_i^{-1})\},\\
    v_+&\mathrel{:=}\tfrac12\max\{1,\|K_i^{-1}\|\}.
\end{aligned}
\]
With $N^2\mathrel{:=}\|u\|^2+\|v\|^2$, we have
$v_-N^2\leq V\leq v_+N^2$ and $|Z|\leq\|C\|N^2/2$.
Choose $\eta>0$ such that
\[
    \eta\|C\|\leq v_-,
    \;\eta\max\{a_0,b_0\}\leq\gamma/2,
\]
and define $W\mathrel{:=}V+\eta Z$. Such a choice exists because
$v_->0$, $\gamma>0$, and all problem and gain constants are finite.
Then
\begin{equation}
\begin{aligned}
    w_-N^2&\leq W\leq w_+N^2,\\
    \dot W&\leq-\tfrac\gamma2(\|u\|^2+\|\psi\|^2)
    -\tfrac\eta2\|C^\top v\|^2,
\end{aligned}
\label{app-global-Wbounds}
\end{equation}
where $w_-\mathrel{:=}v_-/2$ and $w_+\mathrel{:=}v_++\eta\|C\|/2$.
In particular, the derivative bound asserted in the main text holds
with $c_1=\min\{\gamma,\eta\}/2$.

All constants so far are independent of $\varepsilon$ and of the
chosen equilibrium. Only the smooth functions $V$ and $Z$ of the
integral coordinates have been differentiated. The realized vector
field is continuous and globally Lipschitz, so these computations
remain valid at ideal singular complementarity pairs.

\subsection{Equilibrium Existence Under Full Row Rank}
\label{app-global-equilibrium}

Suppose that $C$ has full row rank. The equation
$Cx=d-\boldsymbol{1}$ has a solution, which is strictly feasible.
For $\varepsilon=0$, the feasible set is nonempty and closed, and
strong convexity makes $f$ coercive. Hence a minimizer exists and is
unique. Strict feasibility gives a KKT multiplier, and injectivity of
$C^\top$ makes that multiplier unique. Thus the ideal realized system
has a unique equilibrium.

For $\varepsilon>0$, consider the open nonempty convex domain
$\mathcal D\mathrel{:=}\{x:Cx-d<0\}$ and the barrier objective
\[
    F_\varepsilon(x)\mathrel{:=}f(x)-\varepsilon^2
    \sum_{j=1}^m\log(d_j-C_jx),
    \; x\in\mathcal D.
\]
The barrier terms are convex, so $F_\varepsilon$ is strongly convex.
To verify existence, choose $b\geq1$ with
$d_j-C_jx\leq b(1+\|x\|)$ on $\mathcal D$ for every $j$.
Strong convexity gives constants $b_1,b_2\geq0$ such that
\[
\begin{aligned}
    F_\varepsilon(x)\geq{}&
    \tfrac\ell2\|x\|^2-b_1\|x\|-b_2\\
    &-m\varepsilon^2\log\bigl(b(1+\|x\|)\bigr).
\end{aligned}
\]
For each fixed $\varepsilon$, this lower bound tends to infinity as
$\|x\|\to\infty$. Along a bounded sequence approaching the boundary
of $\mathcal D$, at least one logarithmic term tends to infinity,
while all other terms have a common finite lower bound. Consequently,
barrier sublevel sets are bounded and bounded away from the boundary.
A minimizing sequence therefore has a limit in $\mathcal D$, giving
a unique minimizer $x_\varepsilon^\star$.

The stationarity equation for $F_\varepsilon$ is the original
stationarity equation with
$\mu_{\varepsilon,j}^\star=\varepsilon^2/(d_j-C_jx_\varepsilon^\star)>0$.
This multiplier satisfies the perturbed complementarity equations.
Conversely, any perturbed KKT pair gives a stationary point of the
strongly convex barrier objective and hence this same minimizer and
multiplier. The equilibrium is therefore unique for every
$\varepsilon\geq0$, and its integral state equals its multiplier.

\subsection{The Smoothing-Uniform Exponential Estimate}
\label{app-global-uniform}

With full row rank, $\|C^\top v\|\geq\sigma_{\min}(C)\|v\|$.
Combining this fact with \eqref{app-global-Wbounds} yields
\[
\begin{aligned}
    \dot W&\leq-d_0N^2\leq-\omega W,\\
    d_0&\mathrel{:=}\tfrac12\min\{\gamma,\eta\sigma_{\min}(C)^2\},
    \;\omega\mathrel{:=}d_0/w_+>0.
\end{aligned}
\]
Integration and quadratic comparison give
\[
    \|w(t)-w_\varepsilon^\star\|
    \leq\sqrt{w_+/w_-}\,e^{-\omega t/2}
    \|w(0)-w_\varepsilon^\star\|.
\]
This proves \eqref{eq:smoothing-uniform-ges} with
$A=\sqrt{w_+/w_-}$ and $c=\omega/2$.
The affine constraint map, globally Lipschitz objective gradient,
and globally Lipschitz inverse residual realization give a globally
Lipschitz vector field. Hence solutions exist for all $t\geq0$.
The displayed constants depend on the problem data and fixed gains,
but not on $\varepsilon$, the equilibrium location, or the initial state.

For completeness, this estimate also transfers uniformly to physical
coordinates. Write $q_\mu\mathrel{:=}\mu-\bar\mu$ and $\kappa\mathrel{:=}\|K_p\|$.
For fixed constraint values, the inverse map is $1$-Lipschitz in $\nu$.
For fixed $\nu$, it is at most $2\kappa$-Lipschitz in the constraint
values. Thus
\[
    \|q_\mu\|\leq\|v\|+2\kappa\|C\|\|u\|.
\]
The scalar residual is $2$-Lipschitz in each of its two arguments,
uniformly in $\varepsilon$. Since its equilibrium value is zero,
$\|\psi\|\leq2\|q_\mu\|+2\|C\|\|u\|$. Consequently,
\[
    \|v\|\leq(1+2\kappa)\|q_\mu\|+2\kappa\|C\|\|u\|.
\]
These inequalities make the two state errors globally equivalent with
constants independent of $\varepsilon$. They preserve the exponential
rate and alter only the prefactor. No strict complementarity is used.

\subsection{Proof of the Global Residual Error Bound}
\label{app-global-error-bound}

Set $\varepsilon=0$ and assume a unique KKT pair
$(x^\star,\mu^\star)$ and \eqref{eq:nonnegative-kernel-condition}.
For the fixed gains, write $u\mathrel{:=}x-x^\star$,
$v\mathrel{:=}\nu-\mu^\star$, and
$\psi\mathrel{:=}P^0(Cx-d,\nu)$. We prove
\eqref{eq:global-residual-bound} by contradiction. Suppose the bound fails.
There is then a sequence of states with $t_k\mathrel{:=}\|v_k\|>0$ such that
\[
    \frac{\|u_k\|+\|\psi_k\|+\|C^\top v_k\|}{t_k}\longrightarrow0.
\]
After taking a subsequence, $v_k/t_k\to a$ with $\|a\|=1$, and
$t_k$ tends either to zero, to a positive finite limit, or to infinity.
In every case $C^\top a=0$, while
\[
    \mu_k=\mu^\star+v_k+K_p\psi_k,
    \;g(x_k)=g(x^\star)+Cu_k.
\]

If $t_k\to\infty$, these identities imply
$\mu_k/t_k\to a$ and $g(x_k)/t_k\to0$. Positive homogeneity and
continuity of the ideal scalar FB residual give
\[
    0=\lim_{k\to\infty}\frac{\psi_{k,j}}{t_k}
    =|a_j|-a_j.
\]
Thus $a\geq0$. Together with $C^\top a=0$ and $\|a\|=1$, this
contradicts \eqref{eq:nonnegative-kernel-condition}.

If $t_k\to t\in(0,\infty)$, then $u_k\to0$, $\psi_k\to0$, and
$v_k\to ta$. At the limit the vanishing residual implies that
$\mu^\star+ta$ is nonnegative and complementary to $g(x^\star)$.
Moreover,
\[
    \nabla f(x^\star)+C^\top(\mu^\star+ta)=0.
\]
This is a second KKT multiplier at the same optimizer, contradicting
uniqueness.

It remains to consider $t_k\to0$. Partition the constraint indices as
\[
\begin{aligned}
    \mathcal I&\mathrel{:=}\{j:g_j(x^\star)<0\},\\
    \mathcal S&\mathrel{:=}\{j:\mu_j^\star>0\},\\
    \mathcal W&\mathrel{:=}\{j:g_j(x^\star)=\mu_j^\star=0\}.
\end{aligned}
\]
For $j\in\mathcal I$, the residual is differentiable at the reference
pair, with derivative $-1$ in $\mu_j$ and zero in $g_j$. Since
$\mu_k-\mu^\star=t_ka+o(t_k)$ and $g(x_k)-g(x^\star)=o(t_k)$,
we have $\psi_{k,j}/t_k\to-a_j$, and hence $a_j=0$.
For $j\in\mathcal W$, homogeneity gives
$\psi_{k,j}/t_k\to|a_j|-a_j$, so $a_j\geq0$.
No sign conclusion is needed on $\mathcal S$.

Choose $\theta>0$ sufficiently small that
$\mu_j^\star+\theta a_j\geq0$ for every $j\in\mathcal S$.
On $\mathcal I$, the perturbed multiplier remains zero, while on
$\mathcal W$ it is nonnegative. Thus $\mu^\star+\theta a$ is
complementary to $g(x^\star)$; stationarity is preserved because
$C^\top a=0$. Since $a\neq0$, this again gives a distinct KKT
multiplier. All possible subsequences lead to contradictions, proving
\eqref{eq:global-residual-bound}. This proof requires neither LICQ
nor strict complementarity and makes no assumption about uniqueness
of complementary multipliers at points other than $x^\star$.

\subsection{Global Asymptotic Stability and Exponential Sufficiency}
\label{app-global-sufficiency}

We first prove the asymptotic assertion of
Theorem~\ref{thm:global-geometry}, without imposing
\eqref{eq:nonnegative-kernel-condition}. Set $\varepsilon=0$ and take
the unique KKT pair as the reference equilibrium. The realized field
is globally Lipschitz, and \eqref{app-global-Vdot} holds without a
rank condition. The positive-definite proper function $V$ has compact
positively invariant sublevel sets. A trajectory contained in
$\{\dot V=0\}$ must have $u=0$ and $\psi=0$ throughout.
It follows that $\dot v=0$, $\mu=\nu$, and complementarity holds at
$x^\star$. Maintaining $u=0$ also requires $C^\top v=0$, which gives
stationarity. Uniqueness of the KKT pair therefore makes the largest
invariant subset of $\{\dot V=0\}$ the reference equilibrium.
LaSalle's invariance principle gives global attraction, and the
quadratic bounds on $V$ give Lyapunov stability.

Now impose \eqref{eq:nonnegative-kernel-condition}. Let $K$ be the
constant in \eqref{eq:global-residual-bound}. Squaring that bound
and using $(a+b+c)^2\leq3(a^2+b^2+c^2)$ gives
\[
    N^2\leq(1+3K^2)
    (\|u\|^2+\|\psi\|^2+\|C^\top v\|^2).
\]
Together with \eqref{app-global-Wbounds}, this implies
\[
    \dot W\leq-\omega_0W,
    \;\omega_0\mathrel{:=}\frac{\min\{\gamma,\eta\}}
    {2(1+3K^2)w_+}>0.
\]
Quadratic comparison proves global exponential stability in integral
coordinates, with prefactor $\sqrt{w_+/w_-}$ and rate $\omega_0/2$.
The global coordinate bounds established above transfer both the
asymptotic and exponential conclusions to physical coordinates.

\subsection{Necessity of the Geometric Condition}
\label{app-global-necessity}

Suppose that a vector $a\geq0$ satisfies $\|a\|=1$ and $C^\top a=0$.
Fix any finite admissible gains, without assuming
\eqref{eq:global-gain-condition}. Denote the globally Lipschitz
realized ideal field by $F$ and consider integral-coordinate initial
states $w_M\mathrel{:=}(x^\star,\mu^\star+Ma)$, where $M>0$.
These states lie in the nonnegative integral-state cone, so their
realized multipliers are nonnegative. For $g_j^\star\mathrel{:=}g_j(x^\star)
\leq0$ and any $\mu_j\geq0$,
\[
    g_j^\star\leq
    \sqrt{\mu_j^2+(g_j^\star)^2}-\mu_j+g_j^\star\leq0.
\]
Consequently, $\|\psi(w_M)\|\leq\|g(x^\star)\|$, independently of
$M$. Stationarity and $C^\top a=0$ imply
\[
    F(w_M)=
    \begin{bmatrix}
        -C^\top K_p\psi(w_M)\\
        K_i\psi(w_M)
    \end{bmatrix}.
\]
Thus $\|F(w_M)\|\leq B$ for all $M$, where the finite constant
$B\mathrel{:=}\sqrt{\|C^\top K_p\|^2+\|K_i\|^2}\,\|g(x^\star)\|$ depends
only on the fixed problem and gains.

Choose a positive global Lipschitz constant $L_F$ for $F$.
The integral equation and the bound
$\|F(w)-F(w_M)\|\leq L_F\|w-w_M\|$ give, by Gronwall's inequality,
\[
    \|w(t;w_M)-w_M\|\leq\frac{B}{L_F}(e^{L_Ft}-1).
\]
The right side is independent of $M$ for every fixed time $t$.
If global exponential stability held, constants $A\geq1$ and $c>0$
would satisfy $\|w(t;w_M)-w^\star\|\leq Ae^{-ct}M$ for all $M,t$.
Choose $t_0$ so that $Ae^{-ct_0}<1/2$. The reverse triangle inequality
instead gives
\[
    \|w(t_0;w_M)-w^\star\|
    \geq M-\frac{B}{L_F}(e^{L_Ft_0}-1)>M/2
\]
for sufficiently large $M$, a contradiction.
The global equivalence of physical and integral coordinates implies
that exponential stability in physical coordinates is also impossible.
This proves necessity for every fixed finite admissible gain choice,
regardless of whether the dissipation gain condition holds.

\subsection{Equivalent Geometry and Bounded Feasible Sets}
\label{app-global-polyhedron}

For completeness, the equivalence between
\eqref{eq:nonnegative-kernel-condition} and the existence of $h$ with
$Ch<0$ follows directly from finite-dimensional separation.
If $Ch<0$ and $a\geq0$ is nonzero, then $a^\top Ch<0$, excluding
$C^\top a=0$. Conversely, under
\eqref{eq:nonnegative-kernel-condition}, the origin does not belong
to the convex hull of the row vectors of $C$: convex coefficients
representing zero would provide a nonzero nonnegative kernel vector.
Let $p\neq0$ be the point of minimum Euclidean norm in this compact
convex hull. Optimality gives $p^\top C_j^\top\geq\|p\|^2$ for every
row $C_j$. Hence $h=-p$ satisfies
$C_jh\leq-\|p\|^2<0$ for all $j$.

If the feasible polyhedron is nonempty and $Ch<0$, any feasible
$x_0$ generates the feasible ray $x_0+th$, $t\geq0$.
Since $h\neq0$, the feasible set is then unbounded. A nonempty
bounded polyhedron must therefore violate
\eqref{eq:nonnegative-kernel-condition}. Under the unique-KKT-pair
and dissipation gain assumptions, its ideal FB-PI equilibrium is
globally asymptotically stable by the preceding argument, while the
necessity argument excludes a global exponential estimate uniform
over all initial states for every fixed finite admissible gain choice.

\section{Smoothing and Discrete Implementation}
\label{app-approx-main}

This appendix supplements Section~\ref{sec:fb-dynamics} with equilibrium sensitivity and trajectory perturbation estimates for smoothing, and Section~\ref{subsec:prescribed-rate} with the local Euler stability argument. The equilibrium sensitivity results concern the physical state $z\mathrel{:=}(x,\mu)$, while the global and finite-horizon estimates are first established in integral coordinates $w\mathrel{:=}(x,\nu)$.

\subsection{Equilibrium Sensitivity at a Regular KKT Pair}
\label{app-approx-regular}

Let $z^\star$ satisfy the local assumptions of Section~\ref{subsec:local-structure}. For sufficiently small $\varepsilon\geq0$, there is a unique equilibrium near $z^\star$, denoted by $z_\varepsilon^\star$, and $z_\varepsilon^\star-z^\star=O(\varepsilon^2)$. Its existence and sensitivity do not require stability of the ideal feedback. If the fixed gains make the ideal linearization Hurwitz, then the nearby smoothed equilibrium is also locally exponentially stable for all sufficiently small $\varepsilon$.

To prove these assertions, introduce $\tau\mathrel{:=}\varepsilon^2$ and define
\[
    \mathcal H(x,\mu,\tau)
    \mathrel{:=}\begin{bmatrix}
        \nabla_xL(x,\mu)\\
        \operatorname{diag}(\mu)g(x)+\tau\mathbf1
    \end{bmatrix}.
\]
This map is $C^1$ near $(z^\star,0)$. Order the multiplier coordinates by the active and inactive sets. Its state Jacobian at the reference point is
\[
    \mathcal B\mathrel{:=}
    \begin{bmatrix}
        H&J_A^\top&J_I^\top\\
        \operatorname{diag}(\mu_A^\star)J_A&0&0\\
        0&0&\operatorname{diag}(g_I(x^\star))
    \end{bmatrix}.
\]
If $\mathcal B(u,v_A,v_I)^\top=0$, strict complementarity gives $v_I=0$ and $J_Au=0$. Taking the inner product of $Hu+J_A^\top v_A=0$ with $u$ and using the tangent-space second-order condition yields $u=0$. Full row rank of $J_A$ then gives $v_A=0$. Thus $\mathcal B$ is nonsingular.

The implicit function theorem provides a unique $C^1$ solution branch $z(\tau)$, with
\begin{equation}
    z_\varepsilon^\star
    =z^\star-\varepsilon^2\mathcal B^{-1}
    \begin{bmatrix}0\\\mathbf1\end{bmatrix}
    +o(\varepsilon^2).
    \label{app-approx-regular-expansion}
\end{equation}
For small positive $\tau$, continuity preserves $\mu_A(\tau)>0$ and $g_I(x(\tau))<0$. The product equations then force $g_A(x(\tau))<0$ and $\mu_I(\tau)>0$. Hence this branch lies on the correct perturbed complementarity branch and consists of equilibria of the smoothed FB-PI system. The inverse function argument also gives local uniqueness among these equilibria.

Writing $z'(0)=(u,v_A,v_I)$ gives the more explicit sensitivity equations
\[
\begin{aligned}
    J_Au&=-\operatorname{diag}(\mu_A^\star)^{-1}\mathbf1,\\
    v_I&=-\operatorname{diag}(g_I(x^\star))^{-1}\mathbf1,\\
    Hu+J_A^\top v_A+J_I^\top v_I&=0.
\end{aligned}
\]
In particular, $g_j(x_\varepsilon^\star)=-\varepsilon^2/\mu_j^\star+o(\varepsilon^2)$ for active channels, and $\mu_{\varepsilon,j}^\star=\varepsilon^2/(-g_j(x^\star))+o(\varepsilon^2)$ for inactive channels. The little-$o$ remainder follows from the $C^1$ dependence of the equilibrium branch on $\tau$.

Strict complementarity keeps all radicands away from zero near the ideal pair. The physical and integral vector fields and their state Jacobians therefore vary continuously with $\tau$ along the branch. A Hurwitz ideal Jacobian remains Hurwitz for small $\tau$, proving the stability assertion for the fixed gains. The equilibrium branch itself is independent of these gains.

\subsection{Bias Without Strict Complementarity}
\label{app-approx-weak}

Suppose $f$ is $\ell$-strongly convex and $g(x)=Cx-d$. Let $(x^\star,\mu^\star)$ be an ideal KKT pair and $(x_\varepsilon^\star,\mu_\varepsilon^\star)$ a stationary pair satisfying the positive perturbed complementarity equations. Then
\begin{equation}
    \|x_\varepsilon^\star-x^\star\|
    \leq\sqrt{m/\ell}\,\varepsilon.
    \label{app-approx-primal-bias}
\end{equation}
Indeed, set $\Delta x\mathrel{:=}x_\varepsilon^\star-x^\star$ and $\Delta\mu\mathrel{:=}\mu_\varepsilon^\star-\mu^\star$. Subtracting stationarity and using strong monotonicity yields
\[
\begin{aligned}
    \ell\|\Delta x\|^2
    &\leq-(C\Delta x)^\top\Delta\mu\\
    &=m\varepsilon^2+g(x_\varepsilon^\star)^\top\mu^\star
      +g(x^\star)^\top\mu_\varepsilon^\star
    \leq m\varepsilon^2.
\end{aligned}
\]
Both cross terms are nonpositive by feasibility and multiplier nonnegativity. This proves \eqref{app-approx-primal-bias} without a constraint rank or strict-complementarity assumption.

If $C$ also has full row rank and $\nabla f$ is globally $L_f$-Lipschitz, stationarity gives
\[
    \|\mu_\varepsilon^\star-\mu^\star\|
    \leq\frac{L_f}{\sigma_{\min}(C)}\|x_\varepsilon^\star-x^\star\|.
\]
This follows from $\sigma_{\min}(C)\|\Delta\mu\|\leq\|C^\top\Delta\mu\|=\|\nabla f(x_\varepsilon^\star)-\nabla f(x^\star)\|$.

The first-order bias is attained by the weakly complementary problem $\min_{x\leq0}hx^2/2$, $h>0$. Its ideal KKT pair is $(0,0)$, and its smoothed pair is $(-\varepsilon/\sqrt h,\varepsilon\sqrt h)$. In contrast, for $\min_{x\leq0}(x-1)^2/2$, the strictly complementary reference is $(0,1)$ and
\[
    x_\varepsilon^\star=-\frac{2\varepsilon^2}{\sqrt{1+4\varepsilon^2}+1},
    \; \mu_\varepsilon^\star=1-x_\varepsilon^\star.
\]
These formulas distinguish a first-order weakly complementary bias from the second-order local branch in \eqref{app-approx-regular-expansion}. For any trajectory, the error decomposition $\|z(t)-z^\star\|\leq\|z(t)-z_\varepsilon^\star\|+\|z_\varepsilon^\star-z^\star\|$ separates transient error from equilibrium displacement.

\subsection{Finite-Horizon Consistency Through Singular Pairs}
\label{app-approx-finite-horizon}

Fix a finite interval $[0,T]$ on which an ideal solution exists. If its smoothed counterpart has an initial discrepancy $O(\varepsilon)$ in the complete integral-coordinate state $w=(x,\nu)$, then, for all sufficiently small $\varepsilon$, the smoothed solution exists on $[0,T]$ and differs by $O(\varepsilon)$ uniformly on that interval. If the ideal path stays uniformly separated from every singular pair $(\mu_j,g_j)=(0,0)$ and the initial discrepancy is $O(\varepsilon^2)$, the uniform error improves to $O(\varepsilon^2)$. Both conclusions also hold in physical coordinates.

For the proof, write $k_j\mathrel{:=}k_{p,j}$, $c_j\mathrel{:=}1+2k_j$, $a_j\mathrel{:=}(\nu_j+k_jy_j)^2+c_jy_j^2$. The inverse formulas give
\[
\begin{aligned}
    M_j^\varepsilon-M_j^0
    &=\frac{k_j}{c_j}\bigl(\sqrt{a_j+2c_j\varepsilon^2}-\sqrt{a_j}\bigr),\\
    P_j^\varepsilon-P_j^0
    &=\frac1{c_j}\bigl(\sqrt{a_j+2c_j\varepsilon^2}-\sqrt{a_j}\bigr).
\end{aligned}
\]
Consequently,
\begin{equation}
\begin{aligned}
    \sup_{y_j,\nu_j}|M_j^\varepsilon-M_j^0|
    &=\frac{\sqrt2\,k_j}{\sqrt{c_j}}\varepsilon,\\
    \sup_{y_j,\nu_j}|P_j^\varepsilon-P_j^0|
    &=\frac{\sqrt2}{\sqrt{c_j}}\varepsilon.
\end{aligned}
    \label{app-approx-map-perturbation}
\end{equation}
These equalities follow from $0\leq\sqrt{a+b}-\sqrt a\leq\sqrt b$; the extrema are attained at $y_j=\nu_j=0$, with the first difference identically zero if $k_j=0$.

Let $\mathcal F_\varepsilon(w)$ denote the integral-coordinate vector field. On a compact neighborhood of the ideal path, the Lipschitz bounds in Appendix~\ref{app-real-inversion}, the $C^2$ data, and boundedness of the inverse maps give a state Lipschitz constant $L_0$ uniform over $0\leq\varepsilon\leq\bar\varepsilon$. If $J_0$ bounds $\|J(x)\|$ there, set
\[
\begin{gathered}
    a_M^2\mathrel{:=}\sum_j\frac{2k_j^2}{c_j},\;
    a_P^2\mathrel{:=}\sum_j\frac2{c_j},\\
    B_0\mathrel{:=}\sqrt{J_0^2a_M^2+\|K_i\|^2a_P^2}.
\end{gathered}
\]
Equation~\eqref{app-approx-map-perturbation} implies $\|\mathcal F_\varepsilon(w)-\mathcal F_0(w)\|\leq B_0\varepsilon$. While both solutions remain in this neighborhood, Gronwall's inequality gives
\[
\begin{aligned}
    \|w_\varepsilon(t)-w_0(t)\|
    \leq{}& e^{L_0t}\|w_\varepsilon(0)-w_0(0)\|\\
    &+B_0\varepsilon\frac{e^{L_0t}-1}{L_0},
\end{aligned}
\]
where the quotient is interpreted as $t$ if $L_0=0$. Choose a compact tube around the ideal path strictly inside the neighborhood. For sufficiently small initial discrepancy and $\varepsilon$, the right side is smaller than the tube radius on $[0,T]$. A first-exit argument and continuation then justify existence and the estimate throughout $[0,T]$.

If the ideal physical path is separated from the singular pairs, every $a_j$ has a positive lower bound on a sufficiently small tube: $a_j=0$ is equivalent to $y_j=\nu_j=0$, which is equivalent to $y_j=\mu_j=0$ for the ideal inverse. Rationalizing the square-root difference gives $2c_j\varepsilon^2/(\sqrt{a_j+2c_j\varepsilon^2}+\sqrt{a_j})$. Thus the vector-field perturbation is $O(\varepsilon^2)$ on that tube, and the same argument proves the sharper estimate.

Uniform Lipschitz continuity of $M^\varepsilon$ transfers these estimates to physical coordinates. Conversely, identical physical initial states satisfy $\nu_\varepsilon(0)-\nu_0(0)=-K_p(\Phi^\varepsilon(z(0))-\Phi^0(z(0)))$. This initial discrepancy is $O(\varepsilon)$ in general and $O(\varepsilon^2)$ at nonsingular initial pairs. Hence the same conclusions apply to initialization at the same physical state. All finite-horizon constants may depend on $T$ and the reference path.

\subsection{A Uniform Practical Bound Under Full Row Rank}
\label{app-approx-practical}

Under the full-row-rank assumptions in Section~\ref{sec:global-convergence} and gain condition~\eqref{eq:global-gain-condition}, there are constants $A_0\geq1$, $c_0>0$, and $B_0>0$, independent of $\varepsilon\geq0$ and the initial state, such that
\begin{equation}
    \|w_\varepsilon(t)-w^\star\|
    \leq A_0e^{-c_0t}\|w_\varepsilon(0)-w^\star\|+B_0\varepsilon,
    \label{app-approx-practical-bound}
\end{equation}
where $w^\star\mathrel{:=}(x^\star,\mu^\star)$ is the ideal equilibrium. Thus the smoothed trajectory approaches an $O(\varepsilon)$ neighborhood of the ideal equilibrium with a uniform exponential transient bound.

To prove the bound, put $u\mathrel{:=}x-x^\star$, $v\mathrel{:=}\nu-\mu^\star$, $\psi\mathrel{:=}P^\varepsilon(Cx-d,\nu)$, $h\mathrel{:=}\nabla f(x)-\nabla f(x^\star)$, and $N^2\mathrel{:=}\|u\|^2+\|v\|^2$. For one channel, writing $y\mathrel{:=}g_j(x)$, $y^\star\mathrel{:=}g_j(x^\star)$, and using $\mu^\star$ for $\mu_j^\star$, the FB identity yields
\[
\begin{aligned}
    &(\mu-\mu^\star)[\psi-(y-y^\star)]\\
    &\;- (y-y^\star)\psi+\tfrac12\psi^2\\
    &=\varepsilon^2+y^\star(\mu+\psi)+\mu^\star(y-\psi)\leq\varepsilon^2.
\end{aligned}
\]
Indeed, $\mu+\psi=\sqrt{\mu^2+y^2+2\varepsilon^2}+y\geq0$ and $y-\psi=\mu-\sqrt{\mu^2+y^2+2\varepsilon^2}\leq0$, while $y^\star\leq0$ and $\mu^\star\geq0$. No signs of the current multiplier or integral state are needed.

For $V\mathrel{:=}(\|u\|^2+v^\top K_i^{-1}v)/2$, summing this identity gives
\[
\begin{aligned}
    \dot V&\leq-\ell\|u\|^2+(Cu)^\top\psi-\psi^\top D\psi+m\varepsilon^2,\\
    D&\mathrel{:=}K_p+\tfrac12I.
\end{aligned}
\]
Condition~\eqref{eq:global-gain-condition} makes the symmetric block matrix $S\mathrel{:=}\left[\begin{smallmatrix}\ell I&-C^\top/2\\-C/2&D\end{smallmatrix}\right]$ positive definite. With $\gamma\mathrel{:=}\lambda_{\min}(S)>0$, this gives $\dot V\leq-\gamma(\|u\|^2+\|\psi\|^2)+m\varepsilon^2$.

Add $Z\mathrel{:=}u^\top C^\top v$. Since $\dot u=-h-C^\top(v+K_p\psi)$ and $\dot v=K_i\psi$, differentiation and Young's inequality yield
\[
    \dot Z\leq-\tfrac12\|C^\top v\|^2+A_u\|u\|^2+A_\psi\|\psi\|^2,
\]
where one may take $A_u=L_f^2+\|C^\top K_i\|/2$ and $A_\psi=\|C^\top K_p\|^2+\|C^\top K_i\|/2$. Define $v_-\mathrel{:=}\min\{1,\lambda_{\max}(K_i)^{-1}\}/2$ and $v_+\mathrel{:=}\max\{1,\lambda_{\min}(K_i)^{-1}\}/2$. Choose $\eta>0$ so that $\eta\|C\|\leq v_-$ and $\eta\max\{A_u,A_\psi\}\leq\gamma/2$. Then $W\mathrel{:=}V+\eta Z$ satisfies
\[
    w_-N^2\leq W\leq w_+N^2,
    \; w_-=v_-/2,
    \; w_+=v_++v_-/2,
\]
and
\[
    \dot W\leq-\tfrac\gamma2\|u\|^2-\tfrac\eta2\|C^\top v\|^2+m\varepsilon^2
    \leq-\omega W+m\varepsilon^2,
\]
with $\omega=\min\{\gamma,\eta\sigma_{\min}(C)^2\}/(2w_+)>0$. Integration and the quadratic bounds give \eqref{app-approx-practical-bound} with $A_0=\sqrt{w_+/w_-}$, $c_0=\omega/2$, and $B_0=\sqrt{m/(\omega w_-)}$.

The corresponding physical-coordinate estimate follows from the uniform coordinate Lipschitz bounds and \eqref{app-approx-map-perturbation}. At the reference state, their additional offset is $O(\varepsilon)$; the same is true when expressing the initial integral error in terms of the initial physical error. Thus only the constants in \eqref{app-approx-practical-bound} change.

\subsection{Explicit Physical Dynamics and Euler Stability}
\label{app-approx-euler}

Where all residual components are differentiable, define $r_j\mathrel{:=}\sqrt{\mu_j^2+g_j(x)^2+2\varepsilon^2}$ and
\[
\begin{aligned}
    D_\mu&\mathrel{:=}\operatorname{diag}(\mu_j/r_j-1),\\
    D_g&\mathrel{:=}\operatorname{diag}(1+g_j(x)/r_j).
\end{aligned}
\]
The chain rule in \eqref{eq:fb-pi-implicit} gives the physical vector field
\begin{equation}
\begin{aligned}
    \dot x&=-\nabla_xL(x,\mu),\\
    \dot\mu&=(I-K_pD_\mu)^{-1}
    [K_i\Phi^\varepsilon(x,\mu)+K_pD_gJ(x)\dot x].
\end{aligned}
\label{app-approx-physical-field}
\end{equation}
Each diagonal denominator lies in $[1,1+2k_{p,j}]$. Hence this formula uses componentwise divisions and Jacobian products rather than a coupled linear solve. It is defined globally for $\varepsilon>0$ and in a neighborhood of every strictly complementary ideal equilibrium. At an ideal singular pair, a channel with $k_{p,j}=0$ is directly defined by $\dot\mu_j=k_{i,j}\Phi_j^0$. For nonzero proportional gains, the integral realization handles the singularity without assigning values to undefined derivative ratios.

Let $F_\varepsilon$ denote \eqref{app-approx-physical-field} in a smooth equilibrium neighborhood and $A_\varepsilon\mathrel{:=}DF_\varepsilon(z_\varepsilon^\star)$. Suppose $A_\varepsilon$ is Hurwitz. Simultaneous explicit Euler gives $z_{k+1}=z_k+\alpha F_\varepsilon(z_k)$, with fixed-point Jacobian $B_\alpha\mathrel{:=}I+\alpha A_\varepsilon$. For each eigenvalue $\lambda$, $|1+\alpha\lambda|^2=1+2\alpha\operatorname{Re}\lambda+\alpha^2|\lambda|^2$. Therefore $B_\alpha$ is Schur if and only if
\begin{equation}
    0<\alpha<\min_{\lambda\in\sigma(A_\varepsilon)}
    \frac{-2\operatorname{Re}\lambda}{|\lambda|^2}.
    \label{app-approx-euler-threshold}
\end{equation}
For any $\theta$ strictly between $\rho(B_\alpha)$ and $1$, choose an intermediate number $\widehat\theta\in(\rho(B_\alpha),\theta)$. There is an equivalent vector norm in which the induced norm of $B_\alpha$ is below $\widehat\theta$. Continuity of the derivative of the Euler map makes that map a contraction toward the fixed point with factor at most $\theta$ on a sufficiently small ball in this norm. Norm equivalence yields
\[
    \|z_k-z_\varepsilon^\star\|
    \leq C_\theta\theta^k\|z_0-z_\varepsilon^\star\|
\]
for sufficiently close initial states. This proof covers nondiagonalizable linearizations. Equality in \eqref{app-approx-euler-threshold} gives a unit-modulus eigenvalue and is not included in the exponential-stability conclusion; larger steps produce a linearly unstable fixed point.

Finally, let $\mathcal T_\varepsilon(x,\mu)\mathrel{:=}(x,\mu-K_p\Phi^\varepsilon(x,\mu))$. In a smooth equilibrium neighborhood it is a local diffeomorphism. If $A_w$ is the integral-coordinate linearization and $P\mathrel{:=}D\mathcal T_\varepsilon(z_\varepsilon^\star)$, then $A_w=PA_\varepsilon P^{-1}$, because the vector field vanishes at the equilibrium. Thus Euler linearizations in the two coordinates are also similar and have the same local spectral step-size threshold. Their nonlinear updates generally differ: applying Euler and then a nonlinear coordinate transformation does not give the same map as transforming the continuous-time field and then applying Euler. No such smooth-coordinate argument is used at singular ideal pairs.

\section{Experimental Details and Supplementary Results}
\label{app:experiments}
\label{app:feedback-experiments}

This appendix gives the implementations, parameter calculations, and supplementary tests for Section~\ref{sec:numerical-experiments}. We first specify the feedback laws and accuracy measures, then examine local stabilization, discrete modal design, integration through constraint transitions, and weak complementarity. Additional tests vary the problem size, initial perturbation, and parameter-selection procedure.

\subsection{Implementations and Accuracy Measures}
\label{app:experimental-implementations}

\paragraph{FB-PI}
For scalar gains shared across channels, define
\(
r_j\mathrel{:=}\sqrt{\mu_j^2+g_j(x)^2+2\varepsilon^2},\;
\phi_j\mathrel{:=}r_j-\mu_j+g_j(x).
\)
Where $r_j>0$, the complete FB-PI feedback has the physical-coordinate realization
\begin{equation}
\begin{aligned}
\dot x&=-\nabla f(x)-J(x)^\top\mu,\\
\dot\mu_j&=
\frac{k_i\phi_j+k_p\bigl(1+g_j(x)/r_j\bigr)[J(x)\dot x]_j}
{1+k_p(1-\mu_j/r_j)}.
\end{aligned}
\label{eq:app-fb-physical}
\end{equation}
The denominator lies in $[1,1+2k_p]$. To avoid cancellation when $\mu_j-g_j(x)>0$, the residual is evaluated as
\(
\phi_j=2\bigl(\mu_jg_j(x)+\varepsilon^2\bigr)/(r_j+\mu_j-g_j(x)).
\)
The sparse and three-dimensional Euler experiments use $\varepsilon=10^{-6}$ and update $(x,\mu)$ directly. The ideal local experiments use $\varepsilon=0$; their strictly complementary trajectories remain in the differentiable region.

The scalar complementarity experiments include the singular pair $(\mu,g)=(0,0)$ and use the exact integral realization. For $y=g_j(x)$, define
\(
\nu_j\mathrel{:=}\mu_j-k_p\phi_j,\;
c\mathrel{:=}1+2k_p,\;
w\mathrel{:=}\nu_j+k_py.
\)
The recovered multiplier and residual are
\begin{equation}
\begin{aligned}
M^\varepsilon_j(y,\nu_j)
&\mathrel{:=}\frac{(1+k_p)w+k_p\sqrt{w^2+c(y^2+2\varepsilon^2)}}{c},\\
P^\varepsilon_j(y,\nu_j)
&\mathrel{:=}\frac{\sqrt{w^2+c(y^2+2\varepsilon^2)}-\nu_j+(1+k_p)y}{c}.
\end{aligned}
\label{eq:app-integral-recovery}
\end{equation}
We integrate $\dot x=-\nabla f-J^\top M^\varepsilon$ and $\dot\nu=k_iP^\varepsilon$, initializing $\nu$ from the prescribed physical state. This implementation remains defined at singular ideal complementarity pairs. An Euler update in $(x,\mu)$ and an Euler update in $(x,\nu)$ followed by nonlinear recovery generally give different discrete maps; the discrete comparisons below consistently use the stated coordinates.

\paragraph{AL-PI and projected primal--dual dynamics}
The AL implementation uses proportional feedback of the constraint velocity:
\begin{equation}
\begin{aligned}
p_\rho&\mathrel{:=}[\mu+\rho g(x)]_+,\\
\dot x&=-\nabla f(x)-J(x)^\top p_\rho,\\
\dot\mu&=\frac{k_i}{\rho}(p_\rho-\mu)+k_pJ(x)\dot x,
~~ \rho>0.
\end{aligned}
\label{eq:app-al-physical}
\end{equation}
Euler updates are applied to $(x,\mu)$. Projected primal--dual dynamics (PD) use the same ordinary Lagrangian primal velocity and the componentwise dual velocity $k_i g_j(x)$ when $\mu_j>0$ or $g_j(x)>0$, and zero otherwise. Their discrete implementation is
\begin{equation}
\begin{aligned}
x_{k+1}&=x_k-\alpha\bigl(\nabla f(x_k)+J(x_k)^\top\mu_k\bigr),\\
\mu_{k+1}&=[\mu_k+\alpha k_i g(x_k)]_+.
\end{aligned}
\label{eq:app-pd-euler}
\end{equation}
No multiplier projection is applied to FB-PI or AL-PI during integration or Euler updates.

\paragraph{Squared-slack feedback linearization}
For FL-CMO, let $z\mathrel{:=}(x,s)$, $F(z)\mathrel{:=}f(x)$, $h(z)\mathrel{:=}g(x)+s\odot s$, and $J_h(z)\mathrel{:=}[J(x)\ \ 2\operatorname{diag}(s)]$. With output gain $\kappa>0$, the feedback is
\begin{equation}
\begin{aligned}
\mu&=(J_hJ_h^\top)^{-1}\bigl(\kappa h-J(x)\nabla f(x)\bigr),\\
\dot x&=-\nabla f(x)-J(x)^\top\mu,~~
\dot s=-2\mu\odot s.
\end{aligned}
\label{eq:app-fl-feedback}
\end{equation}
It gives $\dot h=-\kappa h$ in continuous time. Euler updates are applied to $(x,s)$, and $\mu$ is recovered algebraically after each update. The constraint Jacobians in the tested families have disjoint row supports, so the Gram solve is diagonal. In particular, for $g(x)=[I_m\;0]x-d$, its $j$th denominator is $1+4s_j^2$. All four methods exploit the same sparse objective and constraint operations.

\paragraph{Accuracy and termination}
All methods are evaluated on the original problem using
\begin{equation}
\begin{aligned}
R_{\mathrm{KKT}}^2(x,\mu)\mathrel{:=}
&\|\nabla f(x)+J(x)^\top\mu\|^2+\|[g(x)]_+\|^2\\
&+\|[-\mu]_+\|^2+\|\mu\odot g(x)\|^2,\\
E(x,\mu)\mathrel{:=}&\|(x,\mu)-(x^\star,\mu^\star)\|.
\end{aligned}
\label{eq:app-experimental-errors}
\end{equation}
The norms are Euclidean and are not normalized by dimension. For FL-CMO, both quantities use the recovered multiplier. Sparse Euler runs stop at $R_{\mathrm{KKT}}\le10^{-6}$; reference recovery is additionally recorded by $E\le10^{-4}$. The three-dimensional comparison uses the tighter residual threshold $10^{-8}$. A reported iteration count is the first qualifying discrete iterate. Adaptive integrations use a downward residual-threshold event; the accuracy study below separately checks whether an earlier crossing occurs within an accepted step. Counts of right-hand-side evaluations measure integration work, whereas iteration counts measure discrete updates; neither is identified with wall-clock time.

\subsection{Local Stability, Matched Feedback, and Euler Implementation}
\label{app-exp-local}

\paragraph{Exact boundaries and nonlinear trajectories.}
Consider
\[
\min_{x_1\leq0}\;-x_1+\frac12x^\top Hx,~~
H\mathrel{:=}\begin{bmatrix}1&2\\2&1\end{bmatrix}.
\]
Write $z\mathrel{:=}(x_1,x_2,\mu)^\top$ and $z^\star\mathrel{:=}(0,0,1)^\top$.
The reference point is a strict local constrained minimum. Its tangent curvature is positive, whereas $H$ has eigenvalues $-1$ and $3$. The experiment concerns stabilization near this reference point; the objective is unbounded below on the feasible set.

For the ideal FB system with $k_p=0$ and $q\mathrel{:=}k_i$, the physical-coordinate linearization and its characteristic polynomial are
\[
A_q\mathrel{:=}\begin{bmatrix}-1&-2&-1\\-2&-1&0\\q&0&0\end{bmatrix},~~
p_q(s)\mathrel{:=}s^3+2s^2+(q-3)s+q.
\]
For $q>0$, the Routh--Hurwitz conditions reduce to $q>6$. At $q=6$, $p_q(s)=(s+2)(s^2+3)$. For a target decay margin $0<r<1/2$, applying the same criterion to $p_q(s-r)$ gives
\begin{equation}
\gamma(q)>r~\Longleftrightarrow~
q>q_r\mathrel{:=}2(1-r)\left[\frac{4}{1-2r}-(1-2r)\right].
\label{eq:app-local-rate-boundary}
\end{equation}
Indeed, the shifted cubic has coefficients
\(
2-3r,\;q-3-4r+3r^2,\;(1-r)q+3r+2r^2-r^3,
\)
and its remaining Routh--Hurwitz inequality is $(1-2r)(q-q_r)>0$. The boundaries for $r=0.25,0.4,0.49$ are $q_r=11.25,23.76,203.9796$, respectively. Direct eigenvalue evaluation at these gains recovers the prescribed margins to numerical precision.

The nonlinear trajectories use
\(
z(0)=z^\star+10^{-4}v_0,\;
v_0\mathrel{:=}(1,-2,1)^\top/\sqrt6.
\)
DOP853 is run with relative tolerance $10^{-11}$, absolute tolerance $10^{-14}$, and horizon $80$. A trajectory is stopped when $E=10^{-2}$ or $E=10^{-11}$, whichever occurs first. Table~\ref{tab:app-local-spectra} records the spectrum and trajectory endpoints. The unstable case reaches the boundary of the displayed neighborhood at $t=24.9074$; the stable cases decay toward the reference point. The spectral margins are computed from $A_q$, independently of the nonlinear trajectories.

\begin{table}[t]
\centering\small
\caption{Local spectra and nonlinear trajectory endpoints.}
\label{tab:app-local-spectra}
\setlength{\tabcolsep}{5pt}
\begin{tabular}{c r c r c}
\hline
$q$&$\gamma(q)$&$\alpha_{\max}$&End time&End condition\\
\hline
$4$&$-0.157298$&---&$24.9074$&$E=10^{-2}$\\
$8$&$0.123914$&$0.0542797$&$80$&Horizon\\
$24$&$0.401170$&$0.0400388$&$40.2016$&$E=10^{-11}$\\
$100$&$0.479176$&$0.00998265$&$33.0797$&$E=10^{-11}$\\
\hline
\end{tabular}
\end{table}

Here,
\(
\alpha_{\max}\mathrel{:=}\min_{\lambda\in\operatorname{spec}(A_q)}
[-2\operatorname{Re}\lambda/|\lambda|^2]
\)
is defined for the stable linearizations. The large-gain expansions give $\gamma(q)=1/2-2/q+O(q^{-2})$ and $\alpha_{\max}\sim1/q$. Across the spectral grid used in the main figure, the maximum characteristic-polynomial residual, scaled by $1+q$, is $6.09\times10^{-14}$.

\paragraph{Matched active linearizations.}
To examine the local behavior shared by the two feedback structures, fix
\[
k_p^{\mathrm{FB}}=\frac12,~~
k_p^{\mathrm{AL}}=\rho=\frac14,~~ k_i=8.
\]
The common enhanced Hessian is
\(
H_e=\left[\begin{smallmatrix}1.5&2\\2&1\end{smallmatrix}\right],
\)
with eigenvalues $-0.765564$ and $3.265564$. In the respective integral perturbations
\(
\eta_{\mathrm{FB}}=\delta\mu-\tfrac12\delta x_1
\)
and
\(
\eta_{\mathrm{AL}}=\delta\mu-\tfrac14\delta x_1,
\)
both active linearizations are
\[
\begin{bmatrix}\delta\dot x_1\\\delta\dot x_2\\\dot\eta\end{bmatrix}
=
\begin{bmatrix}-1.5&-2&-1\\-2&-1&0\\8&0&0\end{bmatrix}
\begin{bmatrix}\delta x_1\\\delta x_2\\\eta\end{bmatrix}.
\]
Their common characteristic polynomial is $s^3+2.5s^2+5.5s+8$, with roots $-0.321206\pm2.050241i$ and $-1.857587$. Thus, both feedback laws stabilize the reference point with an indefinite enhanced Hessian.

Generate $20$ independent Gaussian vectors in $\mathbb R^3$ using NumPy's \texttt{default\_rng(20261023)} and normalize each to unit length. The same directions are used for both methods at radii $\delta\in\{10^{-4},10^{-3},10^{-2}\}$, giving $120$ trajectories with $z(0)=z^\star+\delta v$. The integrations use physical coordinates, DOP853, tolerances $10^{-11}$ and $10^{-14}$, and horizon $90$. Define the first-crossing times
\[
\begin{aligned}
t_R&\mathrel{:=}\inf\{t:R_{\mathrm{KKT}}(z(t))\leq10^{-8}\},\\
t_E&\mathrel{:=}\inf\{t:E(z(t))\leq10^{-8}\}.
\end{aligned}
\]
Both thresholds are reached in every run. Table~\ref{tab:app-matched-local} reports their medians. Matching the active spectrum gives comparable local decay. The coordinate transformations differ, so a common physical initial state does not give a common integral initial state or an identical physical transient.

\begin{table}[t]
\centering\small
\caption{Median first-crossing times over $20$ shared directions.}
\label{tab:app-matched-local}
\begin{tabular}{c rr rr}
\hline
&\multicolumn{2}{c}{$t_R$}&\multicolumn{2}{c}{$t_E$}\\
$\delta$&FB-PI&AL-PI&FB-PI&AL-PI\\
\hline
$10^{-4}$&$30.2041$&$30.1524$&$28.5639$&$28.5707$\\
$10^{-3}$&$37.1891$&$37.1162$&$35.6289$&$35.7911$\\
$10^{-2}$&$44.3607$&$44.2798$&$42.6652$&$42.6564$\\
\hline
\end{tabular}
\end{table}

For the representative direction $v_0$, each nonlinear trajectory is compared with its own physical-coordinate linearization over $[0,30]$. At the three radii, the maximum FB discrepancies are $4.70155\times10^{-9}$, $4.70137\times10^{-7}$, and $4.69997\times10^{-5}$, respectively, giving approximately $0.470\delta^2$. The AL trajectories show no detected branch change and agree with their affine active dynamics to within $1.57\times10^{-12}$ on the same sampling grid. A separate check initializes both methods at the same exact integral perturbation $\delta v_0$ with $\delta=10^{-4}$, recovering the physical multipliers through the respective coordinate maps. The maximum difference in their integral-coordinate trajectories over the same grid is $4.31\times10^{-9}$. This check uses different physical multiplier initializations and is separate from Table~\ref{tab:app-matched-local}.

First-crossing times are computed from dense output on a grid of spacing $0.005$, followed by Brent refinement of the first detected downward crossing. Repeating this procedure with spacing $0.0025$ changes $t_R$ and $t_E$ by at most $2.41\times10^{-8}$ and $2.92\times10^{-9}$, respectively, over all $120$ runs. This check resolves the brief threshold crossings that can occur between accepted adaptive nodes. For the representative radius $10^{-4}$, tighter integrations with relative and absolute tolerances $2.5\times10^{-13}$ and $2.5\times10^{-16}$ change sampled states by at most $1.93\times10^{-13}$ and $t_R$ by at most $3.36\times10^{-6}$. These comparisons assess numerical consistency of the reported trajectories and crossing times.

\paragraph{Actual Euler trajectories.}
The discrete tests retain the original pure-integral FB system, with $q=24$ or $100$ and $z(0)=z^\star+10^{-6}v_0$. For each gain, physical-coordinate Euler updates use $\alpha=0.9\alpha_{\max}$ or $1.1\alpha_{\max}$ and stop at $E\leq10^{-10}$, $E\geq10^{-2}$, or elapsed time $500$. The outcomes in Table~\ref{tab:app-euler-local} agree with the local Schur stability boundary: steps below the threshold give decay, while steps above it produce departure from the neighborhood. At $0.9\alpha_{\max}$, the larger integral gain reaches the error threshold in less elapsed time but requires more iterations.

\begin{table}[t]
\centering\small
\caption{Euler trajectories below and above the local step limit.}
\label{tab:app-euler-local}
\begin{tabular}{c c c r c}
\hline
$q$&$\alpha/\alpha_{\max}$&$\alpha$&Iterations&End condition\\
\hline
$24$&$0.9$&$0.0360349$&$5952$&$E\leq10^{-10}$\\
$24$&$1.1$&$0.0440426$&$4628$&$E\geq10^{-2}$\\
$100$&$0.9$&$0.00898439$&$19550$&$E\leq10^{-10}$\\
$100$&$1.1$&$0.0109809$&$14672$&$E\geq10^{-2}$\\
\hline
\end{tabular}
\end{table}

\begin{figure*}[t]
\centering
\includegraphics[width=0.8\textwidth]{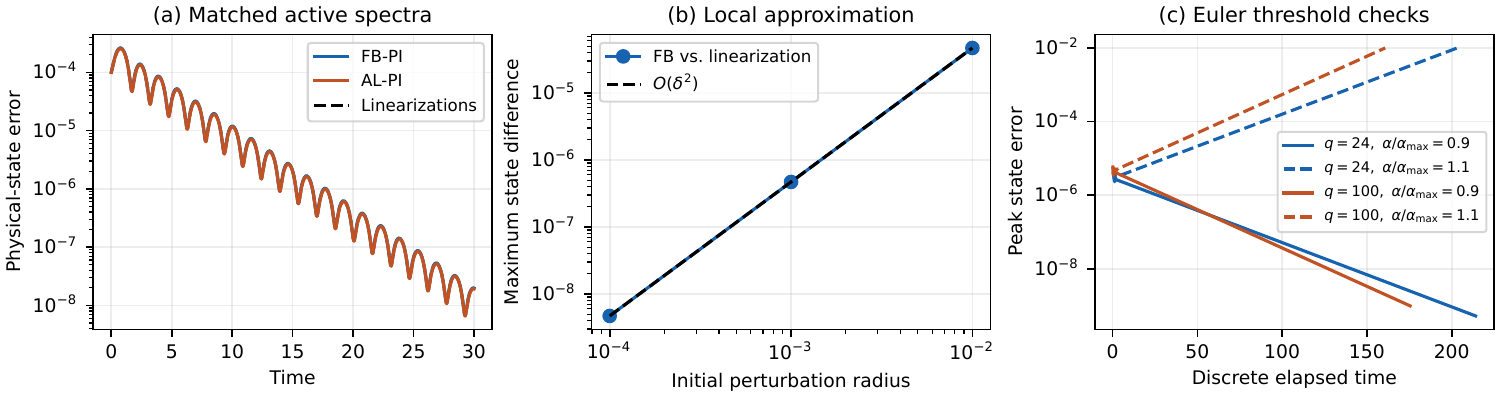}
\caption{Local controls: (a) matched active spectra, (b) FB linearization error, and (c) Euler trajectories below and above the local step limit.}
\label{fig:app-local-controls}
\end{figure*}

Figure~\ref{fig:app-local-controls} collects the matched-spectrum comparison, the local approximation errors, and the Euler tests. The Euler panel connects successive local maxima of the state error to display the decay or growth of the oscillations; the archived trajectories retain every iterate.

\subsection{Exact Discrete Parameter Selection}
\label{app:exp-discrete}
\label{app:euler-comparison}

Consider
\[
\begin{aligned}
\min_{x\in\mathbb R^3}\;&
\frac12x_1^2+x_2^2+\frac12x_3^2-\frac25x_1-x_2,\\
\text{s.t.}\;&\frac25x_1\le0,\;x_2\le0,\;x_3-1\le0.
\end{aligned}
\]
The reference pair is $x^\star=0$, $\mu^\star=(1,1,0)^\top$. Each method uses shared scalar gains and a uniform Euler step in physical coordinates $(x,\mu)$. Put $K\mathrel{:=}k_p^{\rm FB}$ for FB-PI and $K\mathrel{:=}k_p^{\rm AL}+\rho$ for AL-PI. The admissible gain classes are $K\ge0$, $k_i>0$ for FB-PI and $0<\rho\le K$, $k_i>0$ for AL-PI.

\paragraph{Local spectra and Schur conditions}
For an active channel with curvature $h$ and constraint coefficient $\sigma$, the FB linearization is
\[
\delta\dot x=-h\delta x-\sigma\delta\mu,~~
\delta\dot\mu=k_i\sigma\delta x+K\sigma\delta\dot x.
\]
For AL-PI, the linear coordinate change $\delta\widehat\mu\mathrel{:=}\delta\mu+\rho\sigma\delta x$ gives the same equations with $\widehat\mu$ in place of $\mu$. Thus, with $a\mathrel{:=}\alpha$, $p\mathrel{:=}\alpha K$, and $q\mathrel{:=}\alpha^2k_i$, the active Euler characteristic polynomials are
\begin{equation}
 \begin{gathered}
 z^2-(2-B_j)z+1-B_j+d_j,\\
 B_j\mathrel{:=}ah_j+p\sigma_j^2,~~
 d_j\mathrel{:=}q\sigma_j^2,
 \end{gathered}
 \label{eq:app-discrete-poly}
\end{equation}
where $(h_1,\sigma_1)=(1,2/5)$ and $(h_2,\sigma_2)=(2,1)$. The inactive primal eigenvalue is $1-a$ for both methods. The inactive multiplier eigenvalue is $1-q/(a+p)$ for FB-PI and $1-q/t$ for AL-PI, where $t\mathrel{:=}\alpha\rho\in(0,p]$.

Fix a target radius $r\in(0,1)$ and set $u\mathrel{:=}1-r$, $v\mathrel{:=}1+r$. The quadratic Schur conditions for both roots of \eqref{eq:app-discrete-poly} to lie in $|z|\le r$ are
\begin{equation}
 d_j-B_j\le-uv,~~
 uB_j-d_j\le u^2,~~
 vB_j-d_j\le v^2.
 \label{eq:app-discrete-schur}
\end{equation}
These follow by applying the real quadratic Jury criterion after replacing $z$ by $rz$. The inactive primal mode requires $u\le a\le v$. The inactive multiplier condition is
\begin{equation}
\begin{aligned}
 &u(a+p)\le q\le v(a+p),&&\text{FB-PI},\\
 &q\le vp,&&\text{AL-PI}.
\end{aligned}
\label{eq:app-discrete-inactive}
\end{equation}
For AL-PI, the latter condition is equivalent to the existence of
$t\in[q/v,\min\{q/u,p\}]$. At fixed $r$, all these conditions are linear in $(a,p,q)$.

\paragraph{Attaining the optimal factors}
Let $c\mathrel{:=}(2/5)^2=4/25$. If $z_{j,1},z_{j,2}$ are the active Euler roots, then
$d_j=(1-z_{j,1})(1-z_{j,2})$. Hence $|z_{j,k}|\le r$ implies $u^2\le d_j\le v^2$. In particular, $cq\ge u^2$ and $q\le v^2$, giving
\[
 r\ge\frac{1-2/5}{1+2/5}=\frac37.
\]
For FB-PI, the choice
\begin{equation}
 \alpha=\frac{120}{119},~~
 k_p=\frac56,~~ k_i=\frac{289}{144}
 \label{eq:app-fb-euler-optimum}
\end{equation}
gives $(a,p,q)=(120/119,100/119,100/49)$ and the six Euler eigenvalues
\[
 \left\{\frac37,\frac37,-\frac37,-\frac37,
 -\frac1{119},-\frac8{77}\right\}.
\]
The lower bound is therefore attained.

For AL-PI, define $\xi\mathrel{:=}vp-q\ge0$, and abbreviate the first active coefficients by $B\mathrel{:=}a+cp$ and $d\mathrel{:=}cq$. The first two inequalities in \eqref{eq:app-discrete-schur} imply $uB-u^2\le d\le B-uv$, and hence $B\ge2u$. Since $d=vB-va-c\xi$, they also imply
\[
 va+c\xi\ge rB+uv\ge4u-3u^2.
\]
For the second active channel, the third Schur inequality gives $2va+\xi\le v^2$. Since $c\le1/2$, these inequalities yield
\[
 4u-3u^2\le va+c\xi\le\frac{2va+\xi}{2}\le\frac{v^2}{2}.
\]
Equivalently, $7r^2-2r-1\ge0$, so $r\ge(1+2\sqrt2)/7$. At equality, choose
\begin{equation}
\begin{aligned}
 r&=\frac{1+2\sqrt2}{7},~
 \alpha=\frac{1+r}{2},~ k_p^{\rm AL}=0,\\
 \rho&=\frac{2(1-r)^2}{c(1+r)^2},~~ k_i=2\rho.
\end{aligned}
\label{eq:app-al-euler-optimum}
\end{equation}
Then $a=v/2$, $q=u^2/c$, and $p=t=q/v$. The first active block has the repeated root $r$, the second has roots $-r$ and $1-p$, and the remaining roots are $u/2$ and $-r$. Direct substitution gives $|1-p|<r$ and $u/2<r$. Consequently,
\begin{equation}
 r_{\rm FB}^\star=\frac37\approx0.428571,
 ~~ r_{\rm AL}^\star=\frac{1+2\sqrt2}{7}\approx0.546918.
 \label{eq:app-discrete-optima}
\end{equation}
These are minima over the stated shared-gain, uniform-step classes. If a strictly positive AL proportional gain is required, its displayed value is the infimum, approached as $k_p^{\rm AL}\downarrow0$.

The inactive channel determines the separation. At the FB-optimal $(a,p,q)$, the smallest inactive Euler modulus available to AL-PI is
\[
 \min_{0<t\le p}|1-q/t|=q/p-1=\frac{10}{7},
\]
whereas the FB inactive modulus is $8/77$. The FB inactive time scale therefore permits the active blocks to attain their joint lower bound without violating the inactive stability restriction. The repeated active roots can introduce polynomial factors in finite-iteration errors; \eqref{eq:app-discrete-optima} concerns the local asymptotic spectral factors.

\begin{figure*}[t]
\centering
\includegraphics[width=0.8\textwidth]{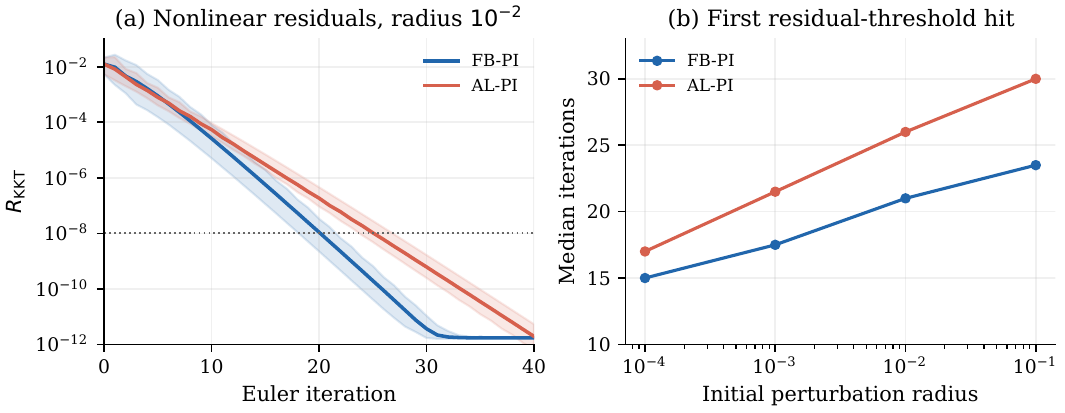}
\caption{Optimized Euler implementations: residuals at $a_0=10^{-2}$ and median iteration counts over $20$ shared directions.}
\label{fig:app-euler-comparison}
\end{figure*}

\paragraph{Nonlinear iterations}
The nonlinear FB implementation uses $\varepsilon=10^{-6}$ and the parameters in \eqref{eq:app-fb-euler-optimum}. AL-PI uses \eqref{eq:app-al-euler-optimum}, numerically
\[
\begin{aligned}
\alpha&\approx0.7734590803,~~ k_p=0,\\
k_i&\approx2.1446609407,~~ \rho\approx1.0723304703.
\end{aligned}
\]
For each $a_0\in\{10^{-4},10^{-3},10^{-2},10^{-1}\}$, Gaussian directions $v\in\mathbb R^6$ are generated with NumPy seeds $32000$--$32019$. Both methods start at $z_0=z^\star+a_0v/\|v\|$, with no clipping of multiplier components. The recorded count is the first index satisfying $R_{\rm KKT}\le10^{-8}$, under a budget of $1000$ iterations. All $160$ trajectories meet this tolerance. The median iteration counts are
\begin{center}
\begin{tabular}{c|rrrr}
$a_0$&$10^{-4}$&$10^{-3}$&$10^{-2}$&$10^{-1}$\\\hline
FB-PI&15&17.5&21&23.5\\
AL-PI&17&21.5&26&30
\end{tabular}
\end{center}
Figure~\ref{fig:app-euler-comparison} shows the corresponding residual curves and iteration counts. Thus, the spectral comparison also describes the ordering observed over these finite initial perturbations.
The residual curves show pointwise medians and the full range over the $20$ directions. They are continued to iteration $40$ for display, while the recorded counts retain the first threshold hit.

\subsection{Modal Parameter Design for Sparse Problems}
\label{app:exp-modal-design}

\paragraph{Training model}
For the sparse family, the active constraint gradients are coordinate vectors. Partition the Lagrangian Hessian into normal and tangent blocks, $H_{AA}$ and $H_{TT}$. The design model retains these two blocks and omits their cross block. Each normal eigenvalue $h$ then gives \eqref{eq:app-discrete-poly} with $\sigma=1$, while each tangent eigenvalue $h_T$ contributes $1-ah_T$. The inactive multiplier conditions remain \eqref{eq:app-discrete-inactive}. Since the inequalities are affine in $h$ and $h_T$, enforcing them at the endpoints of prescribed curvature intervals covers every eigenvalue in those intervals.

The intervals are obtained from three $n=100$ strictly complementary instances with seeds $75001$, $75002$, and $75003$. For the convex family they are
\[
 \begin{aligned}
 h&\in[1.2129468367126781,3.2290477974008334],\\
 h_T&\in[1.1271954626032692,3.2320038380316].
\end{aligned}
\]
At each trial $r$, linear feasibility is tested in $(a,p,q)$. Bisection uses $42$ steps on $[0,1-10^{-7}]$, with HiGHS primal and dual feasibility tolerances $10^{-9}$. A final radius $10^{-6}$ above the numerical feasibility boundary is used, and $p+q$ is minimized over its feasible set. The recovered gains are $\alpha=a$, $k_i=q/a^2$, and $K=p/a$. For AL-PI, $t$ is the projection of $q$ onto $[q/(1+r),\min\{q/(1-r),p\}]$, after which $\rho=t/a$ and $k_p^{\rm AL}=K-\rho$.

The numerical feasibility boundary is approximately $0.5989106924$ for both methods. The selected model radius is $0.5989116924$, with
\begin{equation}
\begin{aligned}
 \alpha&=0.3558285327862788,\\
 k_i&=5.065018582200931,\\
 k_p^{\rm FB}&=2.391617298679508,\\
 k_p^{\rm AL}&=0.5893391680397124,\\
 \rho^{\rm AL}&=1.8022781306397955.
\end{aligned}
\label{eq:app-sparse-modal-gains}
\end{equation}
The designs are computed separately; they yield the same $(\alpha,K,k_i)$ in this family. The selected AL inactive Euler eigenvalue is zero, while the FB inactive modulus is approximately $0.468608$. These parameters are frozen before evaluating the coupled nonlinear test dynamics.

\subsection{Constraint Switching and Integration Accuracy}
\label{app-exp-switching}

\paragraph{Problem construction and parameter matching}
The switching experiment uses
\[
\begin{aligned}
f(x)&\mathrel{:=}\tfrac12x^\top Hx-(C^\top\mu^\star)^\top x,\\
g(x)&\mathrel{:=}Cx-d,~~ C\mathrel{:=}[I_m\;0],
\end{aligned}
\]
with $n=1000$, $m=500$, and $x^\star=0$. Write
\[
H\mathrel{:=}D+0.3L_w,~~
L_w\mathrel{:=}\sum_{j=1}^n w_j(e_j-e_{j+1})(e_j-e_{j+1})^\top,
\]
where $e_{n+1}=e_1$, $D_{jj}\sim\mathcal U[1,2]$, and $w_j\sim\mathcal U[0.5,1.5]$. The first $250$ constraints have $d_j=0$ and $\mu_j^\star\sim\mathcal U[0.8,1.2]$; the remaining constraints have $d_j\sim\mathcal U[0.5,1]$ and $\mu_j^\star=0$. These quantities are generated in the order $D$, $w$, positive reference multipliers, and positive slacks, using independent draws from NumPy's default random-number generator. Problem seeds $95001$, $95002$, and $95003$ define instances $1$, $2$, and $3$. Their initial states are
\[
x(0)=0.8\xi,~ \xi\sim\mathcal N(0,I_n),~~
\mu(0)=0.2\mathbf1,
\]
with initialization seeds $117$, $118$, and $119$, respectively. Both methods use the same physical initial state on each instance.

FB-PI uses $k_p=4$, $k_i=5$, and $\varepsilon=10^{-6}$. AL-PI uses $k_i=5$ and $k_p=4-\rho$ for $\rho\in\{0.25,1,2,4\}$. Thus, the ideal active linearizations coincide in their respective integral coordinates. For each instance and tolerance, the AL parameter with the smallest right-hand-side evaluation count is retained. This selection concerns the four specified matching choices; the nonlinear and inactive-channel dynamics remain different.

Both systems are integrated in physical coordinates $(x,\mu)$ using DOP853 on $[0,60]$, with relative tolerance $\tau\in\{10^{-7},10^{-9}\}$ and absolute tolerance $10^{-2}\tau$. The stopping event is the downward crossing of $R_{\mathrm{KKT}}=10^{-6}$. Performance runs use the solver's event localization with \texttt{dense\_output=False}; the reported evaluation count is \texttt{nfev}, including evaluations associated with rejected trial steps. All runs reach the stopping event before the time limit.

\paragraph{Integration work and switching diagnostics}
Table~\ref{tab:app-switch-candidates} gives the evaluation counts before selecting the AL parameters. The selected values are $\rho=2,1,1$ for instances $1,2,3$, respectively, at both tolerances. At $\tau=10^{-9}$, the median counts are $1265$ for FB-PI and $7145$ for AL-PI. Table~\ref{tab:app-switch-terminal} reports the corresponding stopping times and terminal state errors. The difference in evaluation counts occurs at comparable continuous stopping times.

\begin{table*}[t]
\centering
\caption{Function evaluations for the switching experiment.}
\label{tab:app-switch-candidates}
\small
\setlength{\tabcolsep}{6pt}
\begin{tabular}{ccrrrrr}
\hline
$\tau$ & Instance & FB-PI & AL, $\rho=0.25$ & AL, $\rho=1$ & AL, $\rho=2$ & AL, $\rho=4$\\
\hline
$10^{-7}$ & 1 & 905 & 4397 & 2081 & 1709 & 4625\\
           & 2 & 845 & 4565 & 1625 & 1769 & 4649\\
           & 3 & 905 & 4637 & 1541 & 1553 & 4985\\
$10^{-9}$ & 1 & 1337 & 14273 & 8441 & 8417 & 14837\\
           & 2 & 1253 & 13637 & 7145 & 7745 & 13337\\
           & 3 & 1265 & 12533 & 6389 & 8033 & 15125\\
\hline
\end{tabular}
\end{table*}

\begin{table*}[t]
\centering
\caption{Terminal results at relative tolerance $10^{-9}$.}
\label{tab:app-switch-terminal}
\small
\setlength{\tabcolsep}{7pt}
\begin{tabular}{clrrrrr}
\hline
Instance & Method & Evaluations & Accepted steps & $t_{\rm stop}$ & $E(t_{\rm stop})$ & $\Delta z_{\rm stop}$\\
\hline
1 & FB-PI & 1337 & 95  & 19.265828 & $1.6441\times10^{-6}$ & $2.1130\times10^{-10}$\\
1 & AL-PI & 8417 & 303 & 18.520837 & $2.5192\times10^{-6}$ & $1.3847\times10^{-10}$\\
2 & FB-PI & 1253 & 88  & 19.107265 & $1.5646\times10^{-6}$ & $3.0108\times10^{-9}$\\
2 & AL-PI & 7145 & 257 & 18.392508 & $2.5675\times10^{-6}$ & $4.1028\times10^{-9}$\\
3 & FB-PI & 1265 & 94  & 18.700238 & $1.6234\times10^{-6}$ & $3.0932\times10^{-9}$\\
3 & AL-PI & 6389 & 236 & 18.378079 & $2.6752\times10^{-6}$ & $8.9954\times10^{-10}$\\
\hline
\end{tabular}
\end{table*}

The first instance supplies the representative trajectories. An AL branch change is recorded when the sign of $\mu_j+\rho g_j(x)$ differs between consecutive accepted integration nodes. The resulting $170$ componentwise changes are sampled diagnostics; multiple crossings inside one integration interval need not be resolved. Over $t\leq0.5$, FB-PI and AL-PI take $42$ and $262$ accepted steps, respectively. The final event-shortened interval is omitted from step-size plots because its length is determined by event localization. The minimum sampled value of $\sqrt{\mu_j^2+g_j(x)^2}$ along the FB trajectory is approximately $5.63\times10^{-3}$. The concentration of small AL steps near branch changes, together with the observed separation of FB complementarity pairs from the singular point, is consistent with feedback smoothness reducing adaptive integration work in this experiment.

\paragraph{Reference trajectories and terminal accuracy}
For each selected run, a separate reference trajectory uses relative tolerance $10^{-12}$ and absolute tolerance $10^{-14}$. Its horizon exceeds the larger of the two original stopping times for that method and instance by one time unit. Dense output evaluates this reference at the original accepted nodes and stopping time. With $z\mathrel{:=}(x,\mu)$, the last column of Table~\ref{tab:app-switch-terminal} is
\[
\Delta z_{\rm stop}\mathrel{:=}
\|z_\tau(t_{\rm stop})-z_{\rm ref}(t_{\rm stop})\|.
\]
At $\tau=10^{-9}$, its maximum over the three instances is $3.09\times10^{-9}$ for FB-PI and $4.10\times10^{-9}$ for AL-PI. The reference KKT residuals evaluated at these same stopping times lie in $[0.99715,1.00178]\times10^{-6}$ for FB-PI and $[0.999993,1.000055]\times10^{-6}$ for AL-PI. The maximum absolute difference between reported and reference first-crossing times is $3.50\times10^{-3}$ for FB-PI and $6.67\times10^{-5}$ for AL-PI.

\begin{figure*}[t]
\centering
\includegraphics[width=0.8\textwidth]{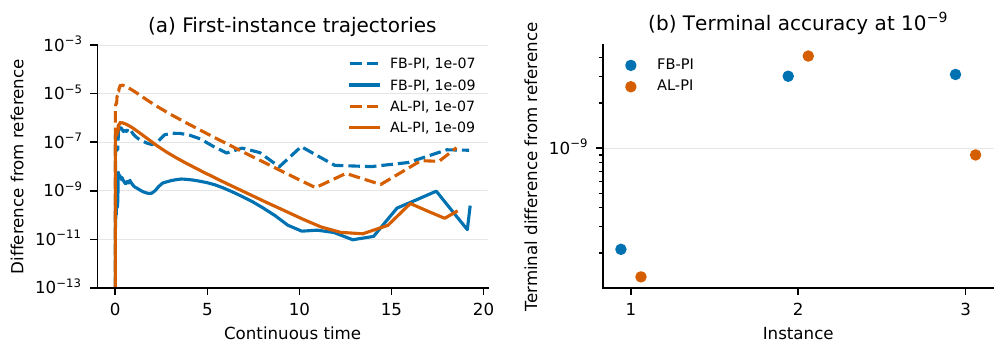}
\caption{Integration accuracy: (a) first-instance discrepancies from tighter reference trajectories; (b) terminal discrepancies for all three instances.}
\label{fig:app-switch-accuracy}
\end{figure*}

For the first instance, a further reference calculation uses relative tolerance $10^{-13}$ and absolute tolerance $10^{-15}$. At the original $10^{-9}$ stopping times, the differences between the two reference solutions are $2.78\times10^{-13}$ for FB-PI and $1.10\times10^{-12}$ for AL-PI. These tighter-solve comparisons assess numerical consistency; they are not certified error bounds relative to an exact trajectory. Figure~\ref{fig:app-switch-accuracy} shows trajectory and terminal discrepancies. Dense-output reference evaluations are excluded from the performance counts in Table~\ref{tab:app-switch-candidates}.

\paragraph{Event localization at the looser tolerance}
The $10^{-7}$ runs give median evaluation counts of $905$ and $1625$ for FB-PI and AL-PI. Their maximum terminal discrepancies from the tighter references are $3.29\times10^{-7}$ and $3.69\times10^{-7}$, respectively. At this tolerance, the residual can cross below and then above the threshold within one accepted step. An event detector based on endpoint signs can consequently report a later downward crossing.

To distinguish solver termination from the first crossing of its computed trajectory, separate dense-output solves were examined on grids of spacing $0.01$ and $0.005$, augmented by the accepted nodes. Downward crossings were bracketed and localized with Brent's method. The two grids gave first-crossing times differing by at most $1.45\times10^{-11}$ across all tested runs and references. Earlier crossings were found in three $10^{-7}$ runs: instance $2$ for FB-PI, and instances $2$ and $3$ for AL-PI. The corresponding delays between the first detected crossing and the original solver termination were $0.68443$, $0.76235$, and $0.31032$. Thus, the $10^{-7}$ counts describe the implemented stopping rule. No earlier crossing was detected for any $10^{-9}$ run under this refinement check; these runs provide the main-text comparison.

\paragraph{Effect of smoothing on the switching trajectory}
For the first instance, ideal FB-PI and FB-PI with $\varepsilon=10^{-6}$ were also integrated at relative tolerance $10^{-13}$ and absolute tolerance $10^{-15}$. Their maximum sampled physical-state difference is $1.91\times10^{-10}$, and their difference at the original smoothed stopping time is $3.99\times10^{-11}$. The minimum sampled complementarity-pair distance along the ideal trajectory is $5.63336\times10^{-3}$. Both reference integrations require $4577$ function evaluations with dense-output construction enabled. This separate comparison indicates that the smoothing perturbation is small along the sampled switching trajectory.

At a smoothed equilibrium, stationarity and perturbed complementarity imply
\[
R_{\mathrm{KKT}}(x^\varepsilon,\mu^\varepsilon)
=\sqrt m\,\varepsilon^2
=2.23607\times10^{-11},
\]
because primal and dual feasibility are strict and each complementarity product has magnitude $\varepsilon^2$. This equilibrium residual is well below the stopping threshold in the switching experiment. The trajectory comparison and equilibrium calculation concern different effects: the former measures the finite-time perturbation of the dynamics, whereas the latter identifies the residual remaining at a smoothed equilibrium.

\subsection{Weak Complementarity and Squared-Slack Dynamics}
\label{app:exp-weak}

\paragraph{Scalar family and common initialization}
Consider
\[
\begin{gathered}
\min_{x\le0}\; f_\delta(x)\mathrel{:=}\tfrac12(x-\delta)^2,\\
(x^\star,\mu^\star)=(\min\{\delta,0\},\max\{\delta,0\}).
\end{gathered}
\]
The tested values are $\delta=0$ and $\delta=\pm10^{-k}$ for $k=0,1,2,3$. They pass from an inactive constraint through weak complementarity to a strictly complementary active constraint. All methods start from $x(0)=-1$ and $\mu(0)=(1+\delta)/5$. FL-CMO uses $s(0)=1$, whose algebraically recovered multiplier agrees with this physical initialization. FB-PI uses the ideal residual in integral coordinates with $k_p=k_i=1$. AL-PI uses $k_p=\rho=1/2$, $k_i=1$; the PD dual gain and FL-CMO output gain are both one.

The integrations use DOP853 with relative tolerance $10^{-10}$, absolute tolerance $10^{-12}$, and horizon $1000$. Sampling combines $1600$ logarithmically spaced positive times with $800$ uniformly spaced times on $[0,40]$. PD is evaluated by joining its positive-multiplier linear arcs and zero-multiplier boundary arcs. The reported scalar threshold times are the first qualifying sample, rather than continuously localized events.

\paragraph{Internal decay on the feasible manifold}
For FL-CMO, $h=x+s^2$ is initially zero and satisfies $\dot h=-h$. Setting $w\mathrel{:=}s^2$ therefore gives
\begin{equation}
 x=-w,~~
 \mu=\frac{w+\delta}{1+4w},~~
 \dot w=-\frac{4w(w+\delta)}{1+4w}.
 \label{eq:app-scalar-slack}
\end{equation}
For $\delta>0$, the linear decay rate of $w$ at zero is $4\delta$. For $\delta<0$, linearization at $w=-\delta$ gives rate $4|\delta|/(1+4|\delta|)$. Both rates vanish as $\delta\to0$.

At weak complementarity, $\delta=0$, separation of variables yields
\begin{equation}
 t=\frac14\left(\frac1{w(t)}-\frac1{w(0)}\right)
   +\log\frac{w(0)}{w(t)},~~ w(t)\sim\frac1{4t}.
 \label{eq:app-slack-exact-solution}
\end{equation}
Since $x+\mu=-4w^2/(1+4w)$, primal and dual feasibility hold and
\begin{equation}
\begin{aligned}
 E&=\left(w^2+\frac{w^2}{(1+4w)^2}\right)^{1/2}
       \sim\frac{\sqrt2}{4t},\\
 R_{\mathrm{KKT}}&=\frac{\sqrt{17}\,w^2}{1+4w}
       \sim\frac{\sqrt{17}}{16t^2}.
\end{aligned}
\label{eq:app-slack-decay-orders}
\end{equation}
The Gram matrix $1+4s^2$ remains bounded below by one. The algebraic tail comes from the internal dynamics on $h=0$, even though the constraint output is exactly regulated. It also produces different decay orders for state error and KKT residual.

\begin{figure*}[t]
\centering
\includegraphics[width=0.8\textwidth]{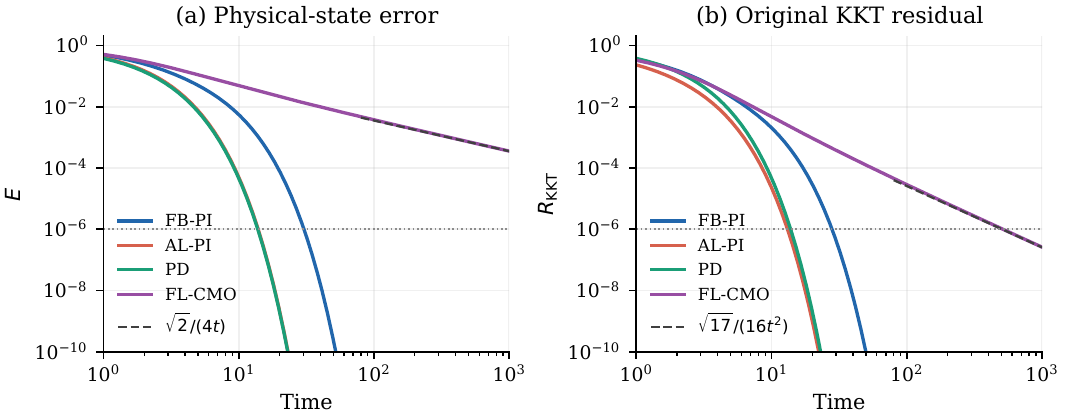}
\caption{Scalar weak complementarity: state errors and KKT residuals, with the FL-CMO asymptotes.}
\label{fig:app-weak-scalar}
\end{figure*}

Figure~\ref{fig:app-weak-scalar} shows the $\delta=0$ trajectories. FL-CMO follows the two algebraic asymptotes, whereas the other three methods show substantially faster decay at the selected gains. Their first sampled times with $R_{\mathrm{KKT}}\le10^{-6}$ are $28.04$, $13.17$, $13.87$, and $519.33$ for FB-PI, AL-PI, PD, and FL-CMO, respectively. The corresponding first sampled times with $E\le10^{-6}$ are $30.39$, $13.97$, and $13.87$ for the first three methods. FL-CMO has $E(1000)=3.56\times10^{-4}$ and $R_{\mathrm{KKT}}(1000)=2.62\times10^{-7}$. At $\delta=-10^{-3}$ and $10^{-3}$, its state errors at the horizon are $2.75\times10^{-5}$ and $2.71\times10^{-5}$, respectively, consistent with slow internal decay near the weakly complementary case.

Evaluating \eqref{eq:app-slack-exact-solution} independently gives a maximum primal-trajectory discrepancy of $4.83\times10^{-10}$. For $\delta\in\{-10^{-3},0,10^{-3}\}$, independent Radau integrations differ from the FB and AL integrated coordinates by at most $1.06\times10^{-10}$. Step-halved projected Euler checks approach the event-resolved PD trajectory with discrepancy-reduction factors approximately $1.98$. The plotted curves are truncated at their first sampled value below $10^{-10}$ to exclude the subsequent numerical-precision tail.

\subsection{Large-Scale Complementarity Transfer and Smoothing Error}
\label{app:exp-large-weak}
\label{subsec:app-four-method-protocol}

\paragraph{Problem family and initialization}
Use the sparse convex family from Appendix~\ref{app-exp-switching}, now with $n=5000$, $m=2500$, and $n/4=1250$ zero-slack constraints. The problem seeds remain $95001$, $95002$, and $95003$. Replace the positive reference multipliers by $\theta\mu^\star$ and change the objective linear term by the same factor, for $\theta\in\{1,10^{-2},0\}$. The Hessian, constraints, and primal solution remain unchanged.

For each instance, a Gaussian direction $v\in\mathbb R^{n+m}$ is generated after the problem data using the same random-number stream. Define the nominal perturbed state by
\[
 (\widetilde x_0,\widetilde\mu_0)
 =(0,\theta\mu^\star)+a_0v/\|v\|,
 ~~ a_0=0.1.
\]
FB-PI, AL-PI, and PD use $x_0=\widetilde x_0$, $\mu_0=[\widetilde\mu_0]_+$. Thus, $a_0$ is the perturbation radius before multiplier projection. FL-CMO shares $x_0$ and uses
\(
s_{0,j}=\sqrt{[-g_j(x_0)]_++10^{-4}},
\)
with its multiplier recovered from \eqref{eq:app-fl-feedback}. The direction and problem data are retained as $\theta$ changes.

Both PI methods use the frozen modal parameters in \eqref{eq:app-sparse-modal-gains}. PD uses $k_i=0.8939481282632452$, $\alpha=0.5498357712120769$, while FL-CMO uses $\kappa=2.607161285838228$, $\alpha=0.4735569851899427$. The latter two choices come from the independent $n=100$ training instances under the finite-search protocol described below. All gains and step sizes are fixed before testing. FB-PI uses $\varepsilon=10^{-6}$. Each run has a $5000$-iteration budget and terminates at the first iterate satisfying $R_{\mathrm{KKT}}\le10^{-6}$.

\paragraph{Finite-search training protocol}
The trajectory-based parameter search uses $n=100$, $a_0=0.1$, and seeds $75001$, $75002$, and $75003$, separately for the convex and nonconvex families. For each method and family, it evaluates $256$ scrambled Sobol candidates followed by $16$ local candidates around each of six leading choices, for $352$ candidates in total. The ranges are $k_i,\kappa\in[0.25,40]$ and $\alpha\in[0.005,0.8]$. FB-PI uses $k_p\in[0.2,16]$; AL-PI uses $k_p+\rho\in[0.2,16]$ and $\rho/(k_p+\rho)\in[0.05,1]$. Positive scale parameters are sampled logarithmically, and the AL fraction is sampled linearly. Local refinements multiply gains and step sizes by factors in $[0.65,1.55]$, perturb the AL fraction by up to $0.2$, and clip to the stated ranges.

Screening uses the first training instance and a $500$-iteration budget. The best $12$ candidates are then evaluated on all three training instances with a $5000$-iteration budget. Selection is lexicographic in the number of reference-recovery failures, number of KKT failures, worst iteration count, median iteration count, and worst terminal residual; the candidate index breaks ties. Reference recovery means reaching the residual threshold with $E\le10^{-4}$. The selected parameters are stored before any held-out evaluation. The main comparison uses this protocol for PD and FL-CMO, while the alternative-parameter study below also uses its FB-PI and AL-PI choices.

\paragraph{Iteration counts and physical-state accuracy}
All $36$ runs reach the residual threshold. Table~\ref{tab:app-complementarity-transfer} reports median iterations and maximum terminal state errors across the three instances. The PI methods retain comparable iteration counts across the three complementarity regimes. FL-CMO is fastest in iterations at $\theta=1$, but its iteration count increases as the active multipliers shrink, in agreement with the loss of slack damping seen in the scalar model. The terminal state errors show why residual and state accuracy are recorded separately.

\begin{table*}[t]
\centering\small
\caption{Large-scale complementarity transfer: median iterations and maximum terminal state errors.}
\label{tab:app-complementarity-transfer}
\begin{tabular}{ccrrrr}
\hline
$\theta$ & Measure & FB-PI & AL-PI & PD & FL-CMO\\
\hline
1 & Iterations & 21 & 20 & 63 & 18\\
  & $\max E$ & $6.34\times10^{-7}$ & $5.89\times10^{-7}$ & $3.86\times10^{-7}$ & $6.07\times10^{-7}$\\
$10^{-2}$ & Iterations & 21 & 20 & 63 & 382\\
  & $\max E$ & $6.47\times10^{-7}$ & $5.89\times10^{-7}$ & $3.86\times10^{-7}$ & $6.14\times10^{-5}$\\
0 & Iterations & 19 & 18 & 63 & 3891\\
  & $\max E$ & $5.23\times10^{-5}$ & $9.23\times10^{-7}$ & $2.81\times10^{-7}$ & $4.16\times10^{-3}$\\
\hline
\end{tabular}
\end{table*}

\paragraph{Equilibrium displacement under smoothing}
At $\theta=0$, write $z^\varepsilon\mathrel{:=}(x^\varepsilon,\mu^\varepsilon)$ for the smoothed equilibrium. For each of the $1250$ active constraints, $g_j(x)=x_j$ and $x_j^\star=\mu_j^\star=0$. Perturbed complementarity gives $-\mu_j^\varepsilon x_j^\varepsilon=\varepsilon^2$, hence
\begin{equation}
 \|z^\varepsilon-z^\star\|^2
 \ge\sum_{j\in\mathcal A}
 \bigl[(x_j^\varepsilon)^2+(\mu_j^\varepsilon)^2\bigr]
 \ge2|\mathcal A|\varepsilon^2.
 \label{eq:app-smoothing-lower-bound}
\end{equation}
For $\varepsilon=10^{-6}$, this yields the lower bound $5\times10^{-5}$.

The actual displacement can be computed independently of the dynamic trajectory. Eliminating the multipliers gives
\begin{equation}
 Hx^\varepsilon+C^\top\mu^\varepsilon=0,~~
 \mu_j^\varepsilon=\frac{\varepsilon^2}{d_j-[Cx^\varepsilon]_j},
 ~ d-Cx^\varepsilon>0.
 \label{eq:app-central-path-system}
\end{equation}
These are the first-order conditions of the strictly convex function
\(
\tfrac12x^\top Hx-\varepsilon^2\sum_j\log(d_j-[Cx]_j),
\)
which has a unique minimizer. For numerical evaluation, set $x=\varepsilon u$ and solve
\[
 Hu+C^\top\frac{\varepsilon\mathbf1}{d-\varepsilon Cu}=0
\]
by sparse Newton iteration with feasibility-preserving backtracking. Division is componentwise. Its Jacobian is
\(
H+C^\top\operatorname{diag}\bigl(\varepsilon^2/(d_j-\varepsilon[Cu]_j)^2\bigr)C.
\)
The initial active coordinates are $u_j=-D_{jj}^{-1/2}$ and all remaining coordinates are zero. Newton steps use an Armijo factor $10^{-4}$, repeated step halving, and stopping tolerance $2\times10^{-13}$ on the scaled stationarity norm. Four updates suffice on each instance.

Table~\ref{tab:app-smoothing-decomposition} separates the equilibrium displacement from the remaining terminal transient. Re-evaluation gives stationarity norms below $9.88\times10^{-21}$ and perturbed-complementarity norms below $6.06\times10^{-28}$. The original KKT residual at each computed equilibrium is $\sqrt m\,\varepsilon^2=5\times10^{-11}$.

\begin{table}[t]
\centering\small
\caption{Weak-complementarity FB errors and the independently computed smoothed equilibrium.}
\label{tab:app-smoothing-decomposition}
\setlength{\tabcolsep}{4pt}
\begin{tabular}{crrr}
\hline
Seed & $\|z^\varepsilon-z^\star\|$ & $\|z_{\rm stop}-z^\star\|$ & $\|z_{\rm stop}-z^\varepsilon\|$\\
\hline
95001 & $5.22744\times10^{-5}$ & $5.23219\times10^{-5}$ & $6.37679\times10^{-7}$\\
95002 & $5.22396\times10^{-5}$ & $5.22904\times10^{-5}$ & $6.83055\times10^{-7}$\\
95003 & $5.22589\times10^{-5}$ & $5.23031\times10^{-5}$ & $6.79576\times10^{-7}$\\
\hline
\end{tabular}
\end{table}

The terminal state lies much closer to the smoothed equilibrium than to the original KKT pair. Thus, the observed FB state-error scale in this weakly complementary test is primarily the smoothing displacement. The transient vector and equilibrium-displacement vector sum to the total error vector; their norms are reported separately rather than added as scalar contributions.

\subsection{Additional Scales, Initial Conditions, and Parameter Choices}
\label{app:exp-robustness}

\paragraph{Alternative parameter selections}
Two further choices examine dependence on the design information. A finite trajectory search on the same training instances gives the ``trajectory'' parameters below. A second modal design uses the construction bounds $h,h_T\in[1,3.8]$: since $H=D+0.3L_w$, $1\le D_{jj}\le2$, and the ring weights lie in $[0.5,1.5]$, the weighted Laplacian satisfies $0\preceq L_w\preceq6I$. The resulting ``family'' design uses the same feasibility procedure and has model radius $0.6747737084$.
The trajectory-search protocol is described in Section~\ref{subsec:app-four-method-protocol}. Table~\ref{tab:app-trajectory-gains} records its frozen choices for both sparse families; the nonconvex family is defined below. The PD and FL-CMO parameters used in the main comparison are the convex entries of this table, whereas the main FB-PI and AL-PI parameters use the modal design in Section~\ref{app:exp-modal-design}.
\begin{table}[ht]
\centering\small
\caption{Frozen trajectory-search parameters.}
\label{tab:app-trajectory-gains}
\setlength{\tabcolsep}{2pt}
\begin{tabular}{llrrrr}
Family&Method&$\alpha$&$k_p$&$k_i$ or $\kappa$&$\rho$\\\hline
Convex&FB-PI&0.277901&3.257643&7.388956&--\\
&AL-PI&0.302416&1.407052&5.847672&1.764899\\
&PD&0.549836&--&0.893948&--\\
&FL-CMO&0.473557&--&2.607161&--\\
Nonconvex&FB-PI&0.395501&4.911755&4.467161&--\\
&AL-PI&0.614685&2.537329&2.270075&1.213186\\
&PD&0.401201&--&13.572169&--\\
&FL-CMO&0.577437&--&1.559955&--\\\hline
\end{tabular}
\end{table}
For the convex family-bound design, the gains are
\[
\begin{array}{c|rrrr}
&\alpha&k_p&k_i&\rho\\\hline
\text{FB-PI}&0.325226&2.349539&5.149539&--\\
\text{AL-PI}&0.325226&0.674774&5.149539&1.674766
\end{array}
\]

The full-precision choices are retained in the parameter files. The modal model isolates the channel time scales used in parameter selection; all nonlinear tests retain the original Hessian coupling.

\paragraph{Scale and initialization tests}
The held-out seeds are $95001$, $95002$, and $95003$. Table~\ref{tab:app-modal-scale} reports median first-hit iteration counts for $R_{\rm KKT}\le10^{-6}$ at nominal initial perturbation radii $a_0=0.1$ and $1$. Initialization and multiplier projection follow the sparse-test protocol in Section~\ref{subsec:app-four-method-protocol}. The convex strict-complementarity results are
\begin{table*}[t]
\centering\small
\caption{Convex sparse tests: median iterations for three parameter selections.}
\label{tab:app-modal-scale}
\begin{tabular}{rr|rrr|rrr}
&&\multicolumn{3}{c|}{FB-PI}&\multicolumn{3}{c}{AL-PI}\\
$n$&$a_0$&Trajectory&Modal&Family&Trajectory&Modal&Family\\\hline
100&0.1&25&19&22&22&18&20\\
100&1&30&23&26&27&22&24\\
1000&0.1&25&21&22&22&21&21\\
1000&1&30&25&26&27&26&26\\
5000&0.1&25&21&22&22&20&21\\
5000&1&30&26&26&27&26&25
\end{tabular}
\end{table*}
For $n=5000$, reducing the reference active multipliers by $\theta$ gives the additional results in Table~\ref{tab:app-modal-weak}. The associated linear term of the objective is changed by the same factor, keeping $x^\star=0$. The parameters are unchanged.
\begin{table*}[t]
\centering\small
\caption{Complementarity transfer at $n=5000$: median iterations.}
\label{tab:app-modal-weak}
\begin{tabular}{rr|rrr|rrr}
&&\multicolumn{3}{c|}{FB-PI}&\multicolumn{3}{c}{AL-PI}\\
$\theta$&$a_0$&Trajectory&Modal&Family&Trajectory&Modal&Family\\\hline
$10^{-2}$&0.1&25&21&22&22&20&21\\
$10^{-2}$&1&30&25&27&27&25&25\\
0&0.1&25&19&21&22&18&20\\
0&1&30&23&26&27&22&25
\end{tabular}
\end{table*}
The modal parameters maintain comparable iteration counts as complementarity weakens. Their derivation uses strictly complementary training points; the weak-complementarity tests examine transfer of the fixed parameters.

\paragraph{Full coupled spectra}
On the three strictly complementary $n=100$ test instances, dense eigenvalue calculations give full Euler spectral radii $0.574606$, $0.600048$, and $0.590878$ for both modal selections. Sparse eigenvalue calculations at $n=5000$ give estimated radii $0.695986$, $0.676462$, and $0.663432$, again equal for the two methods to the reported precision; the maximum computed eigenpair residual is below $5.8\times10^{-11}$. The full model radii differ from the training-model radius because the held-out curvature intervals and the retained coupling need not match the design model. For example, the normal-curvature upper endpoints at $n=5000$ lie between $3.3667$ and $3.4374$, above the training upper endpoint $3.2290$. The family-bound design gives corresponding estimated radii $0.642800$, $0.641006$, and $0.639131$. The computed full-system spectra quantify the effect of transferring the modal parameters to the coupled test instances.

\paragraph{Negative normal curvature}
The scalar problem $\min_{x\le0}(-x-x^2/2+x^4/4)$ has reference pair $(0,1)$ and normal curvature $-1$. Set $x_0=-10^{-3}$, $s_0=\sqrt{-x_0}$, and $\mu_0=(1+x_0-x_0^3)/(1+4s_0^2)$, so that the four representations share their initial physical state. FB-I, FB-PI, AL-PI, and PD use $k_i=5$; their proportional gains are $0$, $4$, $2$, and $0$, respectively, with $\rho=2$ for AL-PI. FL-CMO uses output gain $3$. The ideal FB flows are integrated in integral coordinates. At $t=5$, the error ratios $E(t)/E(0)$ are $6.089$, $4.76\times10^{-4}$, $3.71\times10^{-4}$, $6.097$, and $2.08\times10^{-9}$ for FB-I, FB-PI, AL-PI, PD, and FL-CMO, respectively. DOP853 and Radau trajectories agree to within $1.17\times10^{-11}$ in their integrated coordinates. This example complements the positive-normal-curvature test: here proportional feedback supplies the curvature needed for local stabilization, whereas pure integral feedback and PD depart from the reference point over the displayed interval.

\paragraph{Nonconvex transfer}
An additional family uses $H=D+0.15L_w$, with the first $n/4$ diagonal entries of $D$ replaced by independent samples from $-\mathcal U[1,1.5]$. Write $m=n/2$, $c=0.2$, and $q_4=0.05$, and set
\[
\begin{aligned}
 g_j(x)&=x_j+\frac c2x_{m+j}^2-d_j,\\
 f(x)&=\frac12x^\top Hx-(\mu^\star)^\top x_{1:m}\\
 &~-\frac c2\sum_{j=1}^m\mu_j^\star x_{m+j}^2
 +\frac{q_4}{4}\sum_{j=1}^n x_j^4.
\end{aligned}
\]
The reference Lagrangian Hessian is exactly $H$. Its tangent block is positive definite and its normal block is negative definite. The same training seeds give normal and tangent intervals $[-1.400470,-0.410258]$ and $[1.088297,2.516747]$, respectively. Separate modal designs yield
\[
\begin{array}{c|rrrr}
 &\alpha&k_p&k_i&\rho\\\hline
\text{FB-PI}&0.405321&6.330604&7.977872&--\\
\text{AL-PI}&0.473817&3.015879&3.409335&1.615402
\end{array}
\]
The corresponding model radii are $0.5588904946$ and $0.4843460864$. The trajectory-selected comparison parameters are given in Table~\ref{tab:app-trajectory-gains}.
\begin{table}[ht]
\centering\small
\caption{Nonconvex sparse tests: median iterations.}
\label{tab:app-modal-nonconvex}
\begin{tabular}{rr|rr|rr}
&&\multicolumn{2}{c|}{FB-PI}&\multicolumn{2}{c}{AL-PI}\\
$n$&$a_0$&Trajectory&Modal&Trajectory&Modal\\\hline
100&0.1&31&21&22&14\\
100&1&37&25&26&17\\
1000&0.1&30&21&23&15\\
1000&1&36&27&28&18\\
5000&0.1&30&21&23&14\\
5000&1&37&25&28&17
\end{tabular}
\end{table}
Across the modal-design study, all $396$ runs, including the convex complementarity-transfer cases at each tested size, meet $R_{\rm KKT}\le10^{-6}$ within $5000$ iterations and finish within $10^{-4}$ of the reference pair. Compared with the finite trajectory selections, modal design reduces iteration counts in $69$ of $72$ FB-PI runs and $67$ of $72$ AL-PI runs; it leaves $2$ and $1$ unchanged and increases the remaining counts by at most $1$ and $2$, respectively. The family-bound design improves all $54$ corresponding FB-PI runs and $53$ of $54$ AL-PI runs, with the last count unchanged. The gain-selection results depend on the curvature and channel geometry: the two methods have different optimum model factors for the nonconvex family, whereas they share the same factor for the convex training model.

\paragraph{Step-size refinement with four methods}
The trajectory-selected gains in Table~\ref{tab:app-trajectory-gains} are also tested at $\alpha/4$, with the iteration budget increased from $5000$ to $20000$. Thus the nominal integration horizon remains unchanged. The stopping threshold remains $R_{\mathrm{KKT}}\le10^{-6}$, and reference recovery is recorded separately as $E\le10^{-4}$. A run is stopped for growth if $R_{\mathrm{KKT}}>10^6\max\{R_{\mathrm{KKT}}(0),1\}$ or $E>10^6$; nonfinite states and exhausted budgets are recorded separately.

Table~\ref{tab:app-step-refinement} gives the convex results at nominal perturbation radius $a_0=0.1$. Here both PI methods use the trajectory-selected parameters, so their counts differ from the modal-parameter results in the main text. Each displayed count is the median over the same three test seeds, and every run reaches the KKT threshold.
\begin{table}[ht]
\centering\small
\caption{Median iterations at $\alpha$ / $\alpha/4$ with trajectory-selected gains.}
\label{tab:app-step-refinement}
\setlength{\tabcolsep}{3pt}
\begin{tabular}{rr|rrrr}
$\theta$&$n$&FB-PI&AL-PI&PD&FL-CMO\\\hline
$1$&100&25/118&22/102&46/190&18/82\\
&1000&25/116&22/102&47/198&17/80\\
&5000&25/116&22/103&63/195&18/81\\
$10^{-2}$&100&25/118&21/100&44/178&307/1233\\
&1000&25/116&22/101&44/188&367/1474\\
&5000&25/116&22/102&63/184&382/1532\\
$0$&100&25/119&21/100&41/127&1417/5676\\
&1000&25/116&22/101&37/124&2581/10329\\
&5000&25/116&22/102&63/127&3891/15568\\\hline
\end{tabular}
\end{table}
The slow FL-CMO behavior under weak complementarity remains after refinement. At $n=5000$ and $\theta=0$, its largest terminal state error changes from $4.1644\times10^{-3}$ to $4.1642\times10^{-3}$; the corresponding FB-PI errors are $5.2297\times10^{-5}$ and $5.2268\times10^{-5}$. These errors reflect the state-versus-residual distinction and the smoothing displacement examined above. At the refined step size, the largest terminal errors of AL-PI and PD are $1.1329\times10^{-6}$ and $4.5740\times10^{-6}$, respectively. For $\theta=1$ or $10^{-2}$, all four methods recover the reference pair under the stated $10^{-4}$ criterion. For $\theta=0$, FB-PI, AL-PI, and PD meet this criterion, while FL-CMO does not, at all three sizes and both step sizes.

The strict-complementarity tests also use $a_0=1$. At $n=5000$, the convex median counts at $\alpha$ / $\alpha/4$ are $30/141$, $27/127$, $79/245$, and $22/100$ for FB-PI, AL-PI, PD, and FL-CMO. Each method recovers the reference pair in all $18$ strict-complementarity convex tests combining three sizes, three seeds, and two initial radii, at each step size.

For the nonconvex family, FB-PI, AL-PI, and FL-CMO recover the prescribed reference pair in all $18$ tests at both step sizes. At $a_0=0.1$, their refined-step median iteration counts for $n=100,1000,5000$ are $(140,135,136)$, $(60,62,62)$, and $(72,68,69)$, respectively. PD reaches another KKT point in $6$ tests and triggers the growth condition in $12$ tests at the nominal step size; after refinement it reaches another KKT point in $15$ tests and exhausts the budget in $3$. The other KKT points have objective values below the reference value $f(0)=0$. These outcomes distinguish local attraction to the specified constrained minimum from termination at another stationary point.

\subsection{Implementation Checks and Reproducibility}
\label{app:experiment-consistency}

The physical-coordinate FB velocity was checked against the complete PI identity, and its multiplier recovery was checked against the integral coordinate transformation. Across the tested gains and smoothing values, the largest absolute discrepancies were $2.11\times10^{-15}$ and $4.44\times10^{-16}$, respectively. Finite-difference Jacobians for the sparse convex and nonconvex families agreed with the implemented analytic Jacobians to within $3.30\times10^{-10}$. For the matched local comparison, transforming the two active Jacobians to their respective integral coordinates produced identical matrices; their physical-coordinate matrices need not coincide.

The archived four-method tests contain 432 terminal states, covering 288 strict-complementarity runs and 144 complementarity-transfer runs. Independent evaluation through weighted ring-edge differences verified the gradients, original KKT residuals, and reference errors. The largest relative residual discrepancy was $2.70\times10^{-16}$. Scalar FL-CMO trajectories were checked against the implicit exact solution, FB and AL trajectories against an independent Radau integration, and projected trajectories against step-halved projected Euler integration. The specific error and event-location checks for the new local and switching tests are reported in their corresponding subsections.

The recorded numerical environment uses Python~3.12.14, NumPy~2.3.5, SciPy~1.17.0, and Matplotlib~3.10.8. Sparse operations and linear algebra use the supplied problem implementations; the experiment drivers fix random seeds and restrict BLAS to one thread. The accompanying files include the frozen parameters, trial-level results, stored terminal states, figure data, and scripts for reproduction. Parameter values displayed in tables are rounded, while the supplied JSON files retain their full numerical precision.

\subsection{Constrained Rosenbrock and Griewank Problems}
\label{app:exp-original-benchmarks}

We also examine two nonquadratic benchmark families. The chained Rosenbrock function~\cite{Rosenbrock1960} and the Griewank function~\cite{Griewank1981} are
\[
\begin{aligned}
F_R(x)&\mathrel{:=}\sum_{i=1}^{n-1}\bigl[100(x_{i+1}-x_i^2)^2+(1-x_i)^2\bigr],\\
F_G(x)&\mathrel{:=}1+\frac{\|x\|^2}{4000}-\prod_{i=1}^n\cos\!\left(\frac{x_i}{\sqrt i}\right).
\end{aligned}
\]
Their curved and oscillatory landscapes supplement the quadratic examples above. We form the constrained variants
\[
\min_x\; f(x)\mathrel{:=}F(x)+q^\top x,
\qquad Ax=b,\quad Cx-d\leq0,
\]
with $n=25$, five affine equalities, and respectively $m=10$ and $m=12$ inequalities. The reference points are $x^\star=\mathbf1$ for Rosenbrock and $x^\star=0$ for Griewank. Entries of $A$ and $C$ are sampled independently from the standard normal distribution, using seed $200$ for Rosenbrock and $202$ for Griewank. We set $b=Ax^\star$ and $d_j=C_jx^\star$ for the first four inequalities; the remaining entries have independent slack values in $[0.5,1.3]$. The reference equality multipliers have independent standard normal entries, while the first four inequality multipliers are sampled independently from $\mathcal U[0.2,0.7]$ and the others are zero. Finally,
\(
q\mathrel{:=}-A^\top\lambda^\star-C^\top\mu^\star
\)
makes the prescribed triple satisfy the KKT equations. For feasible $x$,
\[
f(x)-f(x^\star)=F(x)-F(x^\star)
-(\mu^\star)^\top C(x-x^\star)\geq0.
\]
Since each base function has a unique zero at its reference point, $x^\star$ is the unique global minimizer of the corresponding constrained variant.

\paragraph{Multiplier dynamics and initialization}
These problems include equality constraints. Their physical equality multipliers satisfy
\(\dot\lambda=k_i(Ax-b)+k_eA\dot x\),
and the primal channel includes $-A^\top\lambda$. We compare FB-PI and FB-I with AL-PI, projected primal--dual dynamics (PD), proximal augmented-Lagrangian dynamics (Prox-AL), and squared-slack PI dynamics (Slack-PI). For FB-PI and Slack-PI, $k_e=k_p$; FB-I and PD use $k_e=0$. AL-PI uses $k_e=k_p+\rho$. Prox-AL is the case $k_p=0$ of this augmented-Lagrangian implementation and has $k_e=\rho$. The inequality channels of FB-PI and AL-PI are those specified in Appendix~\ref{app:experimental-implementations}, and FB-I sets $k_p=0$. All FB tests here use $\varepsilon=10^{-5}$.

For Slack-PI, let $c\mathrel{:=}Cx-d+s\odot s$ and let $\eta,\xi$ denote integral states. Its physical multipliers and velocities are
\[
\begin{aligned}
\lambda&=\eta+k_p(Ax-b),\qquad \mu=\xi+k_pc,\\
\dot x&=-\nabla f(x)-A^\top\lambda-C^\top\mu,\\
\dot s&=-2\mu\odot s,\\
\dot\eta&=k_i(Ax-b),\qquad \dot\xi=k_ic.
\end{aligned}
\]
This PI controller for the squared-slack equations differs from the feedback-linearization implementation FL-CMO tested in Appendix~\ref{app:exp-weak}.

Every method starts from the same physical triple, generated by
\[
\begin{aligned}
x(0)&=x^\star+0.05\,u_x,\\
\lambda(0)&=\lambda^\star+0.04\,u_\lambda,\\
\mu(0)&=[\mu^\star+0.04\,u_\mu]_++0.08\mathbf1,
\end{aligned}
\]
where the entries of the three vectors are independent standard normal samples. Slack-PI additionally takes
\(s_j(0)=\sqrt{\max\{-g_j(x(0)),10^{-6}\}}\)
and initializes its integral states through the physical-multiplier identities above. This common physical initialization need not satisfy $\nu(0)\geq0$ for FB-PI; the dual-feasibility term records the resulting multiplier sign changes.

\paragraph{Gain selection and accuracy}
For each fixed problem, one independent initialization with seed $9999$ is used to select gains by minimizing the original KKT residual at $t=2$. The common grid has total proportional gains $p\in\{3,5,8\}$ and integral gains $k_i\in\{8,16,24\}$. FB-PI and Slack-PI take $k_p=p$, AL-PI takes $\rho=1$ and $k_p=p-1$, and Prox-AL takes $\rho=p$ and $k_p=0$. PD selects $k_i\in\{4,8,12,16\}$. FB-I uses the selected FB-PI integral gain. All parameters are frozen before the three test initializations, with seeds $5000$, $5001$, and $5002$. Both problems select $(k_p,k_i)=(3,24)$ for FB-PI, $(k_p,k_i,\rho)=(2,24,1)$ for AL-PI, and $(k_i,\rho)=(24,3)$ for Prox-AL. Slack-PI selects $(k_p,k_i)=(3,24)$ for Rosenbrock and $(3,8)$ for Griewank; PD selects $k_i=4$ for both.

The mixed-constraint residual and physical-state error are
\[
\begin{aligned}
R_{\mathrm{mix}}^2\mathrel{:=}\;&
\|\nabla f(x)+A^\top\lambda+C^\top\mu\|^2
+\|Ax-b\|^2\\
&+\|[Cx-d]_+\|^2+\|[-\mu]_+\|^2
+\|\mu\odot(Cx-d)\|^2,\\
E_{\mathrm{mix}}\mathrel{:=}\;&
\|(x,\lambda,\mu)-(x^\star,\lambda^\star,\mu^\star)\|.
\end{aligned}
\]
The first four time units reproduce the original observation horizon, and integrations continue to $t=40$ to distinguish slower transients. The five unprojected flows use BDF with relative tolerance $10^{-9}$ and absolute tolerance $10^{-11}$. PD uses projected Heun with step sizes $2.5\times10^{-4}$ for Rosenbrock and $10^{-3}$ for Griewank, followed by a step-halving check. Solver status and the separate residual components are recorded for every trajectory.

\begin{figure*}[t]
\centering
\includegraphics[width=0.8\textwidth]{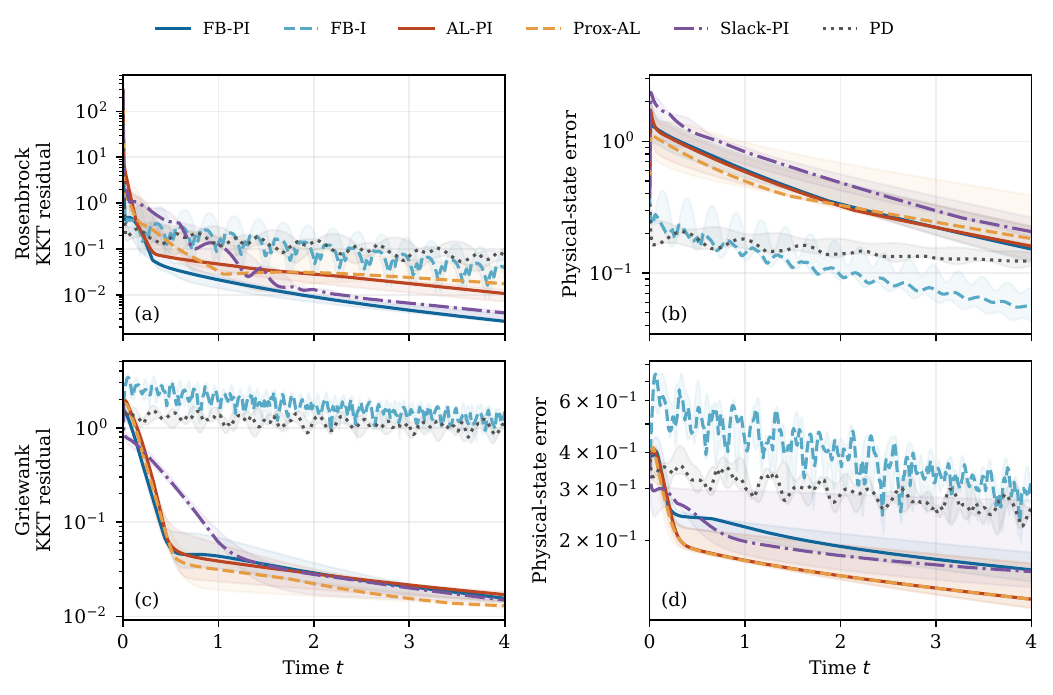}
\caption{Rosenbrock (top) and Griewank (bottom) over the original four-unit horizon: KKT residuals (left) and state errors (right). Lines show medians; bands span three initializations.}
\label{fig:app-nonlinear-transients}
\end{figure*}

\begin{figure*}[t]
\centering
\includegraphics[width=0.8\textwidth]{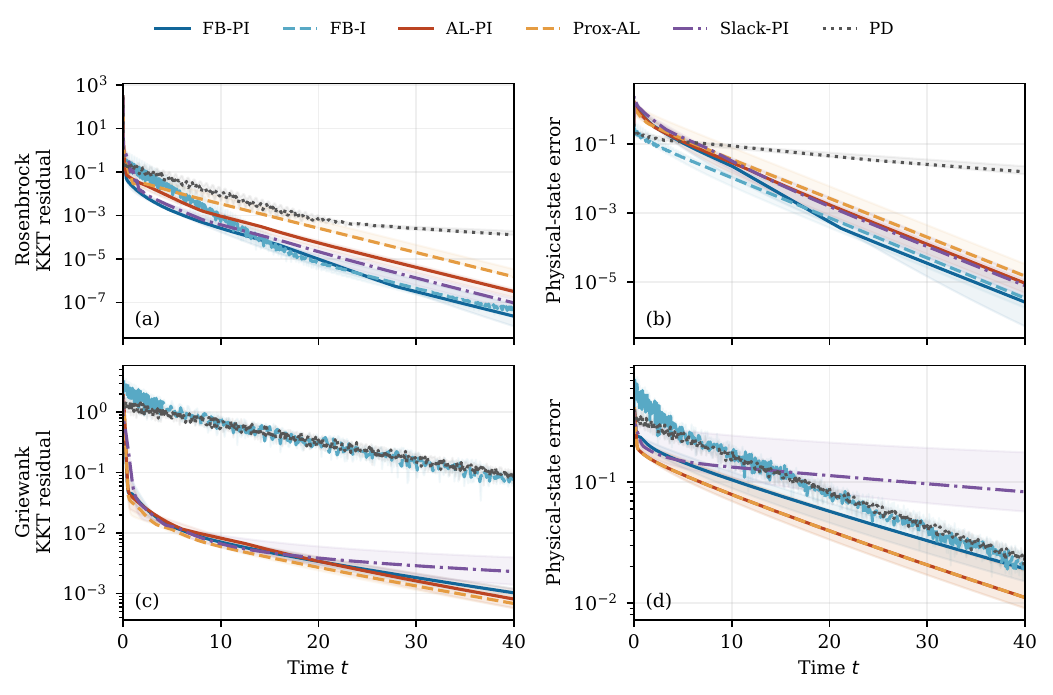}
\caption{The same nonlinear benchmark trajectories continued to $t=40$, with medians and ranges over the three initializations.}
\label{fig:app-nonlinear-extended}
\end{figure*}

\paragraph{Nonlinear transients}
Figures~\ref{fig:app-nonlinear-transients} and~\ref{fig:app-nonlinear-extended} show the median trajectories and the ranges over the three initializations. Table~\ref{tab:app-nonlinear-benchmarks} gives the corresponding residuals and state errors. All $36$ integrations reach $t=40$.

On Rosenbrock, FB-PI gives a median residual of $2.66\times10^{-3}$ at $t=4$, compared with $3.86\times10^{-2}$ for FB-I and $1.07\times10^{-2}$ for AL-PI. The state errors distinguish the remaining transient: at the same time, the FB-PI and FB-I medians are $0.152$ and $0.0554$, respectively. By $t=40$, their residuals are below $10^{-7}$ and their median state errors are below $4\times10^{-6}$.

On Griewank, the FB-PI residual falls to $1.55\times10^{-2}$ at $t=4$, while FB-I gives $1.49$. FB-PI, AL-PI, and Prox-AL retain smaller residuals over the extended horizon, with median values between $6.84\times10^{-4}$ and $1.03\times10^{-3}$ at $t=40$. The slower state decay is accompanied by smaller tangent curvature: the minimum tangent Hessian eigenvalues are approximately $252$ for Rosenbrock and $0.0466$ for Griewank. These examples show how the feedback and local curvature jointly shape transients beyond quadratic objectives.

\begin{table*}[t]
\centering\small
\caption{Nonlinear benchmarks: median residuals and physical-state errors.}
\label{tab:app-nonlinear-benchmarks}
\setlength{\tabcolsep}{7pt}
\begin{tabular}{lrrr|rrr}
\hline
&\multicolumn{3}{c|}{Rosenbrock}&\multicolumn{3}{c}{Griewank}\\
Method&$R_{\rm mix}(4)$&$R_{\rm mix}(40)$&$E_{\rm mix}(40)$
&$R_{\rm mix}(4)$&$R_{\rm mix}(40)$&$E_{\rm mix}(40)$\\\hline
FB-PI&$2.66\times10^{-3}$&$2.37\times10^{-8}$&$2.63\times10^{-6}$&$1.55\times10^{-2}$&$1.03\times10^{-3}$&$1.91\times10^{-2}$\\
FB-I&$3.86\times10^{-2}$&$6.26\times10^{-8}$&$3.43\times10^{-6}$&$1.49$&$8.61\times10^{-2}$&$2.04\times10^{-2}$\\
AL-PI&$1.07\times10^{-2}$&$3.24\times10^{-7}$&$9.29\times10^{-6}$&$1.70\times10^{-2}$&$8.09\times10^{-4}$&$1.11\times10^{-2}$\\
Prox-AL&$1.76\times10^{-2}$&$1.53\times10^{-6}$&$1.49\times10^{-5}$&$1.29\times10^{-2}$&$6.84\times10^{-4}$&$1.11\times10^{-2}$\\
Slack-PI&$4.06\times10^{-3}$&$9.73\times10^{-8}$&$7.92\times10^{-6}$&$1.48\times10^{-2}$&$2.31\times10^{-3}$&$8.33\times10^{-2}$\\
PD&$7.28\times10^{-2}$&$1.30\times10^{-4}$&$1.56\times10^{-2}$&$1.16$&$8.49\times10^{-2}$&$2.52\times10^{-2}$\\
\hline
\end{tabular}
\end{table*}

\paragraph{Numerical verification}
The reference KKT residuals are below $7\times10^{-16}$, and the combined equality and active-inequality Jacobian has rank $9$ in both problems. Finite-difference gradient checks have relative errors below $8.40\times10^{-10}$. Independent Radau integrations of FB-PI and AL-PI from the first test initialization change the sampled residuals by at most $1.68\times10^{-7}$ and the state errors by at most $3.63\times10^{-9}$. Halving the PD step changes the sampled state-error curves by at most $1.01\times10^{-4}$ and the terminal residuals by at most $3.70\times10^{-4}$. The archived results include the frozen gains, individual trajectories, five residual components, and solver status for every run.

\clearpage
\raggedbottom
\end{document}